\documentclass[11pt]{article}

\usepackage{amsmath,amssymb,amsfonts,amsthm,enumitem,accents,mathrsfs}
\usepackage{xy} 
\xyoption{all}

\usepackage[dvipsnames]{xcolor}
\usepackage{soul} 
\usepackage{pstricks,pstricks-add,pst-node}
\usepackage{makeidx}
\usepackage{graphicx}
\usepackage{setspace}
\usepackage{extarrows}
\usepackage{stmaryrd}

\usepackage{hyperref}
\hypersetup{pdftitle={Relative SW invariants and a sum formula}, verbose,
            pdfpagemode=UseThumbs,       
            colorlinks=true,
            linkcolor=black, 
            citecolor=black, 
            urlcolor=black, 
            menucolor=black, 
            anchorcolor=black,
            naturalnames=true,
            pdfstartview=XYZ
            }

\numberwithin{equation}{section}

\usepackage{tocloft}

\usepackage[Symbolsmallscale]{upgreek} 
\usepackage{bbm}        

\newcommand*{\defeq}{\mathrel{\vcenter{\baselineskip0.5ex \lineskiplimit0pt
                     \hbox{\scriptsize.}\hbox{\scriptsize.}}}%
                     =}

\newtheorem{theorem}{Theorem}[section]
\newtheorem{thmx}{Theorem}

\newtheorem{lemma}[theorem]{Lemma}
\newtheorem{conjlemma}[theorem]{Conjectural Lemma}
\newtheorem{corollary}[theorem]{Corollary}
\newtheorem{proposition}[theorem]{Proposition}
\newtheorem{conjecture}[theorem]{Conjecture}

\theoremstyle{definition}
\newtheorem{definition}[theorem]{Definition}
\newtheorem{remark}[theorem]{Remark}
\newtheorem{example}[theorem]{Example}

\def\ov#1{\overline{#1}}
\def\tn#1{\textnormal{#1}}
\def\mf#1{\mathfrak{#1}}
\def\wt#1{\widetilde{#1}}
\def\wh#1{\widehat{#1}}

\def\ll{\left\langle}
\def\rr{\right\rangle}
\def\mc{\mathcal}
\def\lra{\longrightarrow}
\def\dbar{\bar\partial}

\def\ve{\varepsilon}

\newcommand{\abs}[1]{\left\vert #1 \right\vert}

\def\bEq#1{\begin{equation}\label{#1}}
\def\eEq{\end{equation}}

\def\bsEq{\begin{equation*}}
\def\esEq{\end{equation*}}

\def\bDf#1{\begin{definition}\label{#1}}
\def\eDf{\end{definition}}

\def\bTh#1{\begin{theorem}\label{#1}}
\def\eTh{\end{theorem}}

\def\bCn#1{\begin{conjecture}\label{#1}}
\def\eCn{\end{conjecture}}

\def\bLm#1{\begin{lemma}\label{#1}}
\def\eLm{\end{lemma}}

\def\bCLm#1{\begin{conjlemma}\label{#1}}
\def\eCLm{\end{conjlemma}}

\def\bRm#1{\begin{remark}\label{#1}}
\def\eRm{\end{remark}}

\def\bEx#1{\begin{example}\label{#1}}
\def\eEx{\end{example}}

\def\bPr#1{\begin{proposition}\label{#1}}
\def\ePr{\end{proposition}}

\def\bCr#1{\begin{corollary}\label{#1}}
\def\eCr{\end{corollary}}

\def\bFg#1{\begin{figure}\label{#1}}
\def\eFg{\end{figure}}

\def\bPf{\begin{proof}}
\def\ePf{\end{proof}}

\def\bIt{\begin{itemize}[leftmargin=*]}
\def\eIt{\end{itemize}}

\def\bEn{\begin{enumerate}[label=$(\arabic*)$,leftmargin=*]}
\def\eEn{\end{enumerate}}

\def\bEnalph{\begin{enumerate}[label=$(\alph*)$,leftmargin=*]}
\def\eEnalph{\end{enumerate}}

\def\coker{\tn{coker}}

\def\nd{\tn{d}}

\def\ev{\tn{ev}}
\def\GW{\tn{GW}}
\def\SW{\tn{SW}}

\def\ker{\tn{ker}}

\def\sss{\scriptscriptstyle}

\def\cB{\mc{B}}
\def\cC{\mc{C}}
\def\cD{\mc{D}}
\def\cJ{\mc{J}}
\def\cM{\mc{M}}
\def\cO{\mc{O}}
\def\cE{\mc{E}}

\def\cP{\mc{P}}

\def\cN{\mc{N}}
\def\cR{\mc{R}}

\def\cC{\mc{C}}
\def\cH{\mc{H}}

\def\cT{\mc{T}}

\def\R{\mathbb R}
\def\C{\mathbb C}
\def\Z{\mathbb Z}

\def\T{\mathbb T}

\def\mfi{\mf{i}}
\def\mfj{\mf{j}}

\def\mfs{\mf{s}}
\def\mft{\mf{t}}

\def\la{\lambda}

\def\De{\Delta}

\def\om{\omega}
\def\Om{\Omega}
\def\si{\sigma}
\def\Si{\Sigma}
\def\al{\alpha}
\def\be{\beta}

\def\Ga{\Gamma}
\def\ze{\zeta}
\def\na{\nabla}

\def\cancel#1#2{\ooalign{$\hfil#1\mkern1mu/\hfil$\crcr$#1#2$}}
\def\dirac{\mathpalette\cancel\partial}             
\def\G{\ifmmode{\cal G}\else{${\cal G}$}\fi}        
\def\spin{\ifmmode{\mathrm{spin}}\else{spin}\fi}
\def\Spin{\ifmmode{\mathrm{Spin}}\else{Spin}\fi}
\def\spinc{\ifmmode{\mathrm{spin}^c}\else{spin$^c$ }\fi}
\def\Spinc{\ifmmode{\mathrm{Spin}^c}\else{Spin$^c$ }\fi}
\def\SW{\ifmmode{\tn{SW}}\else{SW}\fi}

\def\mend{\mathrm{End}}
\def\Im{\mathrm{Im}}
\def\det{\tn{det}}
\def\bbC{\mbox{${\mathbb C}$}}

\def\bb1{\mbox{${\mathbbm 1}$}}         

\def\negsp{\mkern-9mu}                  
\def\mquadt{\kern-1em}                  
\def\mqquadt{\kern-2em}                 
\def\clm{\rho}                      
\def\sb{\mathbb{S}}                 
\def\ub{\mathbb{U}} 
\def\scan{\sb_{\tn{can}}}
\def\ucan{\ub_{\tn{can}}}
\def\clmcan{\clm_{\tn{can}}}
\def\lcan{L_{\tn{can}}}
\def\kcan{K_{\sss \Si}}             
\def\piysi{\pi}                     
\def\sycan{\scan}                   
\def\clmycan{\clmcan}
\def\gmetric{\ifmmode{\mathrm{g}}\else{$\mathrm{g}$ }\fi} 
\def\dsw{\ifmmode{\mathsf{d}}\else{${\mathsf{d}}$ }\fi} 
\def\tcyend{\mathsf{t}}
\def\etacyend{\hat{\eta}}           
\def\sfrak{\ifmmode{\mathfrak s}\else{${\mathfrak s}$ }\fi}  
\def\tfrak{\ifmmode{\mathfrak t}\else{${\mathfrak t}$ }\fi}

\def\mfsy{\mfs_{\sss Y}}            
\def\mfsyeclass{\mfs_{\sss Y,[E]}}  
\def\mfsye{\mfs_{\sss Y,E}}         
\def\mfssie{\mfs_{\sss \Si,E}}      
\def\mfsn{\mfs_{\sss\cN}}           

\def\syfrak{\ifmmode{\sfrak_{\sss Y}}\else{$\sfrak_{\sss Y}$ }\fi} 
\def\smfrak{\ifmmode{\sfrak_{\sss M}}\else{$\sfrak_{\sss M}$ }\fi}
\def\mfsm{\mfs_{\sss M}}            
\def\sz{\ifmmode{{\tn\ss}}\else{\ss}\fi} 
\def\sbold{\ifmmode{\mbox{\boldmath$s$\unboldmath}}\else{\boldmath$s$\unboldmath}\fi}
\def\tbold{\ifmmode{\mbox{\boldmath$t$\unboldmath}}\else{\boldmath$t$\unboldmath}\fi} 
\makeindex
\begin{document}
\title{Relative Seiberg-Witten invariants\\ and a sum formula}
\author{Mohammad Farajzadeh-Tehrani\footnote{Supported by NSF grant DMS-2003340}~~and Pedram Safari}
\date{\today}
\maketitle
\begin{abstract}
We study relative Seiberg-Witten moduli spaces and define relative invariants for a pair $(X,\Sigma)$ consisting of a smooth, closed, oriented 4-manifold $X$ and a smooth, closed, oriented 2-dimensional submanifold $\Sigma\!\subset\!X$ with positive genus. 
These relative Seiberg-Witten invariants are meant to be the counterparts of relative Gromov-Witten invariants. 
We also obtain a sum formula (aka a product formula) that relates the SW invariants of a sum $X$ of two closed oriented 4-manifolds $X_1$ and $X_2$ along a common oriented surface $\Sigma$ with dual self-intersections to the relative SW invariants of $(X_1,\Sigma)$ and $(X_2,\Sigma)$. 
Our formula generalizes Morgan-Szab{\'o}-Taubes' product formula.
\end{abstract}
\tableofcontents
\section{Introduction}\label{intro}
Relative moduli spaces associated to a pair $(X,D)$ of a symplectic manifold/complex projective variety $X$ and a (symplectic) divisor $D$ have many applications in both symplectic and complex algebraic geometry. 
In particular, we use moduli spaces of  stable relative maps to define relative Gromov-Witten (or GW) invariants. 
We also use them to derive a GW sum formula when a smooth one-parameter family of symplectic manifolds/varieties $\{Y_\la\}_{\la\in \C}$ degenerates to a simple normal crossing variety $Y_0= \bigcup_{i\in \mc{I}} X_i$.

\medskip
\par\noindent
In real dimension four, Seiberg-Witten (or SW) moduli spaces have been an excellent tool for exploring smooth closed orientable 4-manifolds.  
When $X$ is symplectic, by a celebrated work of Taubes, we know that the SW invariants of $X$ are equal to certain (possibly disconnected) GW invariants \cite{t3}. 
This correspondence has been a great tool for classifying symplectic $4$-manifolds; e.g. see \cite{L}.
In light of SW-GW correspondence, it is natural to seek a construction of relative SW invariants for every pair $(X,\Si)$ of a closed oriented 4-manifold $X$ and a closed, oriented, possibly non-connected, 2-dimensional submanifold $\Si\!\subset\!X$. 
It is also natural to expect a sum formula that expresses the SW invariants of a sum $X$ of two oriented 4-manifolds $X_1$ and $X_2$ along a common oriented surface $\Si$ with dual self-intersections in terms of the relative SW invariants of $(X_1,\Si)$ and $(X_2,\Si)$. It is expected that in the symplectic case and under Taubes' GW-SW correspondence, such an SW sum formula should correspond to the well-known GW sum formula. 
$$
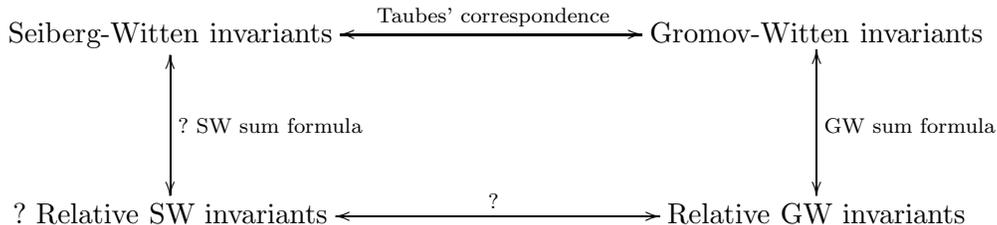
\xymatrix{
\tn{Seiberg-Witten invariants} \ar@{<->}[rrrr]^{\small{\tn{Taubes' correspondence}}} \ar@{<->}[dd]^{\small{\tn{?~SW sum formula}}} &&&& \tn{Gromov-Witten invariants} \ar@{<->}[dd]^{\small{\tn{GW sum formula}}}\\
&&&&\\
\tn{?~Relative SW invariants} \ar@{<->}[rrrr]^{?}&&&& \tn{Relative GW invariants}\\
}
$$

\medskip
\par\noindent
If $\Si$ has genus $g\!>\!0$ and trivial normal bundles in $X_1$ and $X_2$, then Morgan-Szab{\'o}-Taubes' product formula \cite[Thm.~3.1]{mst}, when the restriction of the characteristic line bundle to $\Si$ has degree $2g-2,$ is a particular case of this picture. 
Applications of their formula include the proof of Thom conjecture \cite{mst} --- also proved independently in \cite{KM1} using a vanishing argument --- as well as construction of non-symplectic 4-manifolds with non-trivial SW invariants \cite{FS}. 
Also, the vanishing result of \cite[Thm.~2.1]{OS} is another particular case of the SW sum formula envisioned above.

\medskip
\par\noindent
In this paper, we introduce a setup, in two formats, for constructing relative SW moduli spaces and invariants without any major restriction on the topology of $(X,\Si)$. 
In one format, as formulated in Theorem~\ref{Relaive_thm2}, we obtain relative SW invariants $\wt\SW{}^{\sss X,\Si}$ and a sum formula that generalizes Morgan-Szab{\'o}-Taubes' product formula, but its relation to the GW side is not clear to us. 
In the other format, formulated in Theorem~\ref{Relaive_thm}, we obtain relative SW invariants $\SW^{\sss X,\Si}$ and a sum formula, subject to the regularity of the tunneling spaces, which should be equivalent to the GW side of the diagram above under Taubes' GW-SW correspondence. 
We believe, moreover, that for certain choice of $\Om$ in (\ref{SWI-Relative_e}), the relative invariants provided by these two approaches are equivalent. 

\medskip
\par\noindent
Given a smooth, closed, oriented 4-manifold $X$ and a transverse union $\Si=\bigcup_{i\in\mc{I}}\Si_i$ of positively intersecting, closed, oriented, 2-dimensional sub-manifolds $\Si_i\!\subset\!X$, we can define a notion of logarithmic tangent bundle $TX(-\log\Si)$ that coincides with the corresponding notions in algebraic geometry/symplectic topology whenever $(X,\Si)$ is complex/symplectic; see \cite{FMZ}. There is a vector bundle homomorphism
$$
\iota \colon TX(-\log \Si) \lra TX,
$$
covering $\tn{id}_X$, which is an isomorphism away from $\Si$. 
In this paper, for simplicity, we assume that $\Si$ is smooth. 
We will then have 
$$
TX(-\log\Si)|_{\Si}\cong T\Si\oplus\C_\Si,
$$
where $\C_{\Si}$ is the trivial complex line bundle $\Si\times\C$.
In order to define relative SW invariants, our idea is to use the logarithmic tangent bundle $TX(-\log\Si)$ instead of the classical tangent bundle $TX,$ as is traditionally used in SW theory. The advantage of this approach is two-fold.  
First of all, it is amenable to working with cylindrical-end manifolds, for which a good deal of literature is available at our disposal \cite{mmr,tl2,mst,moy,N,N1,gluingbook,KM}. We will review the essential material in Section~\ref{setup_s} for the sake of coherence of the presentation and establishing the notation. 
In this context, we will work over the cylindrical-end manifold $X\!-\!\Si$ and simply use the logarithmic tangent bundle to distinguish different topological components of particular SW moduli spaces in dimensions 3 and 4.
\par\noindent
From another perspective, using the logarithmic tangent bundle enables us to directly work over the closed manifold $X$ and generalize the results to arbitrary normal crossing case, where it is not feasible to study $X\!-\!\Si$ as a manifold with one cylindrical end. This approach is outlined in Section~\ref{LSW_s}, where we use a so-called logarithmic connection on $TX(-\log \Si)$ and logarithmic connections on the related spinor bundles to derive a logarithmic version of the SW equations. We will treat the details of this direct approach and the general normal crossing case in a separate paper.
\medskip
\par\noindent
Here is an outline of the rest of the paper. 
In Section~\ref{LogStr_ss}, we define the notion of a relatively canonical \spinc structure for the logarithmic tangent bundle. 
Every other \spinc structure will then be obtained by tensoring a relatively canonical \spinc structure with a hermitian line bundle $E$ on $X$. 
Let $\sfrak$ be a \spinc structure on $TX(-\log\Si)$ with spinor bundle $\sb\!=\!\sb^+\oplus\sb^-$ and characteristic line bundle 
$$L=L_\mfs=\tn{det}(\sb^+)=\tn{det}(\sb^-).$$
Define the degree of $\sfrak$ along $\Si$ to be\footnote{For simplicity, we are assuming that $\Si$ is connected. If $\Si$ is not connected, some of the numbers and structures in the following description, such as the integer $d(\sfrak),$ should be defined component-wise.} the integer
\bEq{DegS_e}
d(\sfrak)\defeq (g-1)- \abs{m(\mfs)}, \qquad m(\mfs)=\frac{1}{2}\deg(L_{\sss \Si}),
\eEq
where $g$ is the genus of $\Si$ and $L_{\sss \Si}\!\defeq\!L_\mfs|_\Si.$ 
The motivation for this definition comes from complex geometry: if $X$ is a symplectic $4$-manifold, then $TX(-\log\Si)$ can be equipped with a complex structure. 
In this case, if $\mfs_{\tn{can}}$ is the canonical \spinc structure on $TX(-\log\Si)$ and $\mfs_{\sss E}\!=\!\mfs_{\tn{can}}\otimes E$ is the canonical \spinc structure twisted by a complex line bundle $E$, then $d(\mfs_{\sss E})$ in (\ref{DegS_e}) is simply the degree of $E_{\sss \Si}\!\defeq\! E|_\Si$ or of the Serre dual $\kcan\otimes E^*_{\sss \Si},$ depending on the sign of $\deg(L_{\sss \Si})$ being non-positive or positive. If we set $\sz\!=\!\tn{PD}(c_1(E))\!\in\!H_2(X,\Z),$ then the \spinc structure $\mfs_{\sss E}$ is determined by $\sz$ and $d(\mfs_{\sss E})$ is equal to the product of homology classes $\Si\cdot\sz$ (or $\Si\cdot(\kcan-\sz)$). 
\medskip
\par\noindent
We focus next on $X\!-\!\Si$. 
The restriction of a suitable metric on $TX(-\log \Si)$ to $X\!-\!\Si$ gives $X\!-\!\Si$ the structure of a manifold with cylindrical end $[0,\infty)\times Y,$ where $Y$ is a circle bundle over $\Si.$
We will then proceed in two steps, which we call {\it gliding} and {\it descent}, as described below. 
Let $\cM(X\!-\!\Si,\sfrak_{\sss X-\Si})$ denote the moduli space of monopoles on $X\!-\!\Si$ with finite energy on the end with respect to the induced $\spinc$ structure $\sfrak_{\sss X-\Si}\defeq\sfrak|_{\sss X-\Si}$ on $X\!-\!\Si.$ 
Different $\spinc$ structures $\mfs$ and $\mfs'$ on $TX(-\log \Si)$ may have the same restriction on $X\!-\!\Si$. 
Therefore, for such $\mfs$ and $\mfs'$,  $\cM(X\!-\!\Si,{\sfrak}_{\sss X-\Si})$ will be the same as $\cM(X\!-\!\Si,\sfrak'_{\sss X-\Si})$. 
In the following, by looking at the limits of the monopoles over the cylindrical end, we will choose a component of $\cM(X\!-\!\Si,{\sfrak}_{\sss X-\Si})$ corresponding to $\mfs$.

\medskip
\par\noindent
In temporal gauge, a monopole on $X\!-\!\Si$ with finite energy on the end is the gradient flow line of the Chern-Simons-Dirac function $\tn{CSD}$ on the cylindrical end (see Section~\ref{SW4cyend_ss}), where a stationary solution on the end corresponds to a monopole on $Y$. 
Therefore, we have a {\it limiting} map
\bEq{Limit-Map_e}
\partial:\cM(X\!-\!\Si,\sfrak_{\sss X-\Si})\lra\cM(Y,\sfrak_{\sss Y}),
\eEq
where $\sfrak_{\sss Y}$ is the induced $\spinc$ structure on $Y$. 
Meanwhile, as we will see in Section~\ref{SW3_ss}, an irreducible monopole on $Y$ in this setup actually {\it descends} to a monopole on $\Si,$ i.e., it will be the pullback of a solution of $\SW^2$ on $\Si$ via the bundle projection $Y\!\lra\!\Si.$ 
Because of this, the moduli space $\cM(Y,\sfrak_{\sss Y})$ 
decomposes into a disjoint union of components
$$\cM(Y,\sfrak_{\sss Y})=\mc{J}\cup \bigcup_{m} \cM(Y,\sfrak_{\sss Y})_{m},$$
where $\mc{J}\cong \tn{Map}(\Si,S^1)\cong \T^{2g}$ is the subspace of reducible solutions\footnote{This component is empty if $\Si\cdot\Si=0$ and $\tn{deg}(L_{\sss \Si})\neq 2-2g,$ and is equal to $\T^{2g+1}$ if $\Si\cdot\Si=0$ and $\tn{deg}(L_{\sss \Si})\!=\! 2-2g$. Since the case $\Si\cdot\Si=0$ is treated extensively in \cite{mst}, throughout the paper, we often implicitly or explicitly assume that $\Si\cdot\Si\neq 0$.} and the second union is over the set of integers\footnote{For $m\!\not\in\! [1\!-\!g, g\!-\!1]$, the component $\cM(Y,\sfrak_{\sss Y})_{m}$ is empty, and for $m=0,$ we get the reducible component $\mc{J}$ instead. Different $m$'s among non-trivial components differ by a multiple of $\Si\cdot\Si$.} $m=m(\mfs)$ for all \spinc structures $\mfs$ on $TX(-\log \Si)$ that restrict to $\mfs_{\sss X-\Si}$ on $X\!-\!\Si$. 
The integer $m(\mfs)$ and the restriction $\mfs_{\sss X-\Si}$ uniquely determine $\mfs;$ therefore, for each $\mfs,$ we define
$$
\cM(X\!-\!\Si,\sfrak)\defeq\partial^{-1} \big(\cM(Y,\sfrak_{\sss Y})_{m(\mfs)}\big) \subset \cM(X\!-\!\Si,\sfrak_{\sss X-\Si}).
$$
By definition, all the monopoles in $\cM(X\!-\!\Si,\sfrak)$ are automatically irreducible.
We also define 
$$
\cM(X\!-\!\Si,\sfrak_{\sss X-\Si},\cJ)=\partial^{-1} \big(\cJ\big).
$$
Thus, we have a decomposition 
$$
\cM(X\!-\!\Si,\sfrak_{\sss X-\Si})=\cM(X\!-\!\Si,\sfrak_{\sss X-\Si},\mc{J})\cup\bigcup _{\mfs} \cM(X\!-\!\Si,\sfrak),
$$
where the second union runs over all \spinc structures $\mfs$ on $TX(-\log\Si)$ that restrict to the fixed \spinc structure $\sfrak_{\sss X-\Si}$ on $X\!-\!\Si.$ 

\medskip
\par\noindent
A monopole in $\cM(Y,\sfrak_{\sss Y})_{m(\mfs)}$ can be identified with an effective divisor of degree $d(\mfs)>0$ on $\Si,$ or a single point if $d(\mfs)=0.$
Depending on whether $m(\mfs)\!<\!0$ or $m(\mfs)\!>\!0$, this divisor is the zero set of a non-zero holomorphic section of $E_{\sss \Si}$ or $\kcan\otimes E_{\sss \Si}^*$, respectively. For $d>0,$ let
$$
\tn{Div}_{d}(\Si)\cong \tn{Sym}^{d}(\Si)
$$
denote the space of effective divisors of degree $d$ on $\Si;$
if $d=0,$ this space is taken to be a single point.
By the argument above, we obtain a {\it landing} map 
$$
\flat:\cM(Y,\sfrak_{\sss Y})_{m(\mfs)}\lra \tn{Div}_{d(\mfs)}(\Si)
$$ 
which assigns to each monopole on $Y$ the corresponding effective divisor on $\Si$ of degree $d(\mfs)$ on which the non-zero spinor lands; this will be the ``empty" divisor if $d(\mfs)\!=\!0$. The landing map is a diffeomorphism. Combining the two steps, and using the nomenclature of Gromov-Witten theory, we obtain an {\it evaluation} map 
\bEq{Ev-Map_e}
\tn{ev}\defeq\flat\circ\partial:\cM(X\!-\!\Si,\sfrak)\lra\tn{Div}_{d(\mfs)}(\Si).
\eEq

\medskip
\par\noindent
In order to make sure that the moduli space is cut transversely and contains no reducible solutions, we need to perturb the SW equations on $X\!-\!\Si$ in a way that (\ref{Ev-Map_e}) is still defined. 
This can be achieved by using either a {\it compact perturbation} or an {\it adapted perturbation}; see Section~\ref{SW4cyend_ss}. 
In the first case, the perturbation term is a compactly supported self-adjoint $2$-form $\eta_o$. 
In the second case, restricted to the neck $[0,\infty)\!\times\! Y,$ the perturbation term is additionally equal to the self-adjoint part of the pullback of a non-trivial 2-form $\eta$ on $Y$.  
Furthermore, following \cite{mst}, we consider a special type of adapted perturbations $\eta_\nu,$ where $\eta$ is obtained from a holomorphic $1$-form $\nu$ on $\Si$; see Lemma~\ref{Nic-lm}. 
We denote the resulting moduli spaces by $\cM_{\eta_o}(X\!-\!\Si,\sfrak)$ and $\cM_{\eta_{\nu}}(X\!-\!\Si,\sfrak),$ respectively.  
The discussion leading to (\ref{Ev-Map_e}) readily generalizes to the compactly-perturbed moduli spaces $\cM_{\eta_o}(X\!-\!\Si,\sfrak)$. 
In the second case, $\tn{CSD}$ is circle-valued, but (\ref{Ev-Map_e}) is still defined and we get a major restriction on its image. Given a non-trivial holomorphic $1$-form $\nu$, vanishing on 
$$
\tn{Div}(\nu)=\{p_1,\ldots,p_{2g-2}\}\!\subset\! \Si
$$ 
(counted with multiplicities), let $S_d(\nu)$ ($\cong S_{2g-2-d}(\nu)$)
denote the set of subsets of $\tn{Div}(\nu)$ of size~$d$.
In the case of $\cM_{\eta_{\nu}}(X\!-\!\Si,\sfrak)$, $\cJ$ will be empty and (\ref{Ev-Map_e}) takes values in $S_{d(\mfs)}(\nu)$. Therefore, 
$$
\cM_{\eta_\nu}(X\!-\!\Si,\sfrak_{\sss X-\Si})=\bigcup _{\mfs} \cM_{\eta_\nu}(X\!-\!\Si,\sfrak).
$$ 

\begin{remark}
We have assumed $g\!>\!0$; otherwise, there is no such $\nu$. Furthermore, if $g\!=\!0$, we have $\cM_{\eta_o}(X\!-\!\Si,\sfrak_{\sss X-\Si})=\cM_{\eta_o}(X\!-\!\Si,\sfrak_{\sss X-\Si},\cJ)$. 
\end{remark}

\begin{thmx}\label{Relaive_thm}
Let $X$ be a smooth, closed, oriented $4$-manifold and $\Si\!\subset\!M$ a smooth, closed, oriented surface with positive genus. 
If $b^{+}_{\sss X-\Si}>0,$ for a generic compact perturbation $\eta_o$, $\cM_{\eta_o}(X\!-\!\Si,\sfrak)$ is a smooth orientable (but not necessarily compact) manifold of real dimension
\bEq{DimFormula_e}
\dsw=\frac{(c_1(L)+\Si)^2-2\chi(X)-3\si(X)}{4}.
\eEq
Furthermore, the evaluation map (\ref{Ev-Map_e}) is a smooth submersion. 
\end{thmx}

\medskip
\par\noindent
Unlike the classic case in SW theory, the moduli spaces $\cM_{\eta_o}(X\!-\!\Si,\sfrak)$ are not necessarily compact. 
Because of the tunneling phenomenon (see \cite[Sec.~16]{KM} or \cite[Sec.~4.4.2]{N}), a sequence of monopoles in $\cM_{\eta_o}(X\!-\!\Si,\sfrak)$ will, after passing to a sub-sequence, ``converge" to a finite ordered set of monopoles (called a \textit{broken trajectory}\textnormal{}~in \cite{KM}), where the first one is a monopole on $X\!-\!\Si$ and the rest are non-trivial (i.e., non-stationary) monopoles on the cylinder $\R\!\times\! Y.$ 
Furthermore, 
\bEn
\item the limit at $+\infty$ of the monopole on $X\!-\!\Si$ defined in (\ref{Limit-Map_e}) coincides with the limit at $-\infty$ of the monopole on the first copy of $\R\!\times\! Y$, and 
\item the limit at $+\infty$ of the monopole on the $i$-th copy of  $\R\!\times\! Y$ coincides with the limit at $-\infty$ of the monopole on the $(i+1)$-th copy of $\R\!\times\! Y$. 
\eEn
In order to obtain a compact moduli space without boundary, we will consider the monopoles on the cylinder $\R\!\times\! Y$ up to the natural $\C^*$-action generated by translation in the $\R$-direction and rotation in the $Y$-direction. 
We denote the space containing such limits by $\ov\cM_{\eta_o}(X-\Si,\sfrak)$. 
With the exception of taking quotients by $\C^*\!=\!\R\times S^1$ instead of just by $\R$, this is the same compactification considered in \cite{KM,OS}. 
The limiting map (\ref{Limit-Map_e}) on 
$$
\ov\cM_{\eta_o}(X\!-\!\Si,\sfrak)\subset \ov\cM_{\eta_o}(X\!-\!\Si,\sfrak_{\sss X-\Si})$$ 
is taken at $+\infty$ of the last monopole. 
This compactification is the direct analogue of the relative compactification of the moduli space of pseudo-holomorphic curves relative to $\Si$ (with the contact orders all equal to one). 

\medskip
\par\noindent
While we show that the complement
$$
\ov\cM_{\eta_o}(X\!-\!\Si,\mfs)-\cM_{\eta_o}(X\!-\!\Si,\mfs)
$$
is a finite union of strata of expected codimension at least 2, it is not clear to us if the tunneling spaces are always manifolds of the expected dimension; see Section~\ref{Tunneling_ss}. Therefore, without further restrictions (e.g. as in \cite{OS}), it may not always be the case that, for generic $\eta_o$, $\ov\cM_{\eta_o}(X-\Si,\mfs)$ is a $C^0$-manifold of the expected dimension (\ref{DimFormula_e}). Whenever the latter happens, the relative Seiberg-Witten invariants of $(X,\Si)$ in the class of a logarithmic $\spinc$ structure $\mfs$ are defined by integration on $\ov\cM_{\eta_o}(X-\Si,\mfs)$ in the following way. 

\medskip
\par\noindent
Let 
$$
c\in H^2(\ov\cM_{\eta_o}(X\!-\!\Si,\sfrak),\Z)
$$ 
denote the first Chern class of the natural circle bundle on the moduli space as in the classical case (see \cite[Sec.~6.7]{m} or \cite[p.~249, (7.24)]{DS}). For $\Om\!\in\!H^*(\tn{Div}_{\sss d(\mfs)}(\Si))$ satisfying 
$$
\tn{deg}(\Om)+2r=\dsw,
$$
define 
\bEq{SWI-Relative_e}
\SW^{\sss X,\Si}_{\eta_o}(\mfs;\Om)
\defeq\int_{\ov\cM_{\eta_o}(X-\Si,\mfs)} c^{r}\wedge\tn{ev}^*\Om
=\int_{\cM_{\eta_o}(X-\Si,\mfs)} c^{r}\wedge\tn{ev}^*\Om \, .
\eEq
Under the assumption that the tunneling spaces have their expected dimensions, the inclusion $\cM_{\eta_o}(X\!-\!\Si,\mfs)\subset \ov\cM_{\eta_o}(X\!-\!\Si,\mfs)$ defines a pseudo-cycle in the sense of \cite{Z}; this justifies the second equality in (\ref{SWI-Relative_e}). 
In this situation, unlike in the definition of the classical SW invariants, the condition $b^{+}_{\sss X-\Si}\!>\!0$ is sufficient to conclude that $\SW^{\sss X,\Si}_{\eta_o}(\mfs;-)$ is independent of the choice of generic $\eta_o$ and the cylindrical metric, since all the monopoles in $\cM_{\eta_o}(X\!-\!\Si,\mfs)$ are automatically irreducible by definition. We will elaborate more on this point in Section~\ref{Compact_ss}. Therefore, we drop $\eta_o$ from the notation and denote the invariants by $\SW^{\sss X,\Si}(\mfs;-).$ When $X$ is symplectic, we believe that a special case of (\ref{SWI-Relative_e}), described below, is equal to certain relative Gromov-Witten invariant in the sense of Taubes' correspondence.

\medskip
\par\noindent
Given a $\spinc$ structure $\sfrak$ on $TX(-\log\Si)$ such that
$d(\mfs)\!\geq\!0$ (otherwise the moduli space is empty), let $\mft\!=\!\{t_1,\dots,t_k\}$ be a (possibly empty) partition of $d(\mfs)$ into $k=k(\mft)$ positive integers, i.e.,
\bEq{DivisorData_e}
d(\mfs)\!=\!t_1+\dots+t_k.
\eEq
Let 
$$
\tn{Div}_{\mft}(\Si)=\Big\{\sum_{i=1}^k t_i \tn{p}_i\colon~ \tn{p}_i\!\in\!\Si~~~\forall ~i=1,\ldots,k\Big\} \subset \tn{Div}_{d(\mfs)}(\Si)
$$ 
denote the subspace of effective divisors that can be written as a sum of $k$ (not necessarily distinct) points $\{\tn{p}_i\}_{i=1}^k$ with multiplicities $\{t_i\}_{i=1}^k$.  If we sort the numbers in (\ref{DivisorData_e}) so that 
$$
t_1=\cdots=t_{i_1}<t_{i_1+1}=\cdots=t_{i_1+i_2}<\cdots<t_{i_1+\cdots+i_{n-1}+1}=\cdots=t_{i_1+\cdots +i_n},\quad k=i_1+\cdots+i_n,
$$
then 
$$
\tn{Div}_{\mft}(\Si)\cong \tn{Sym}^{i_1}(\Si)\times \cdots\times \tn{Sym}^{i_n}(\Si).
$$
For $\Omega=\tn{PD}(\tn{Div}_{\mft}(\Si))$ in (\ref{SWI-Relative_e}), we define a particular type of relative SW invariants by
$$\SW^{\sss X,\Si}(\mfs;\mft)\defeq
\SW^{\sss X,\Si}\big(\mfs;\tn{PD}(\tn{Div}_{\mft}(\Si))\big).$$
If X is symplectic, as mentioned earlier, we can identify the set of \spinc structures $\mfs$ on $TX(-\log \Si)$ with the set of homology classes $\sz\in\!H_2(X,\Z)$. 
Under this identification, (\ref{DimFormula_e}) simplifies to 
$$
\dsw=\sz\cdot\sz-K_X\cdot\sz
$$
and $\mft$ is a partition of $d(\mfs)=\Si\cdot\sz$ (or $\Si\cdot(\kcan-\sz)$) into a sum $t_1\!+\cdots+t_k$ of positive integers.
Then we believe that, similar to Taubes' GW-SW correspondence theorem, 
$\SW^{\sss X,\Si}(\sz;\mft)$ is equal to certain relative Gromov-Witten invariant $\GW^{\sss X,\Si}(\sz;\mft)$; it is a count of (possibly disconnected) $J$-holomorphic curves of degree $\sz$ that intersect $\Si$ in $k$ points with tangency orders $\mft\!=\!\{t_1,\dots,t_k\}$ and pass through $r$ generic points in $X$.

\medskip
\par\noindent
In order to resolve the issues above regarding regularity of tunneling spaces and non-compactness, we may consider adapted perturbations $\eta_\nu$ corresponding to holomorphic $1$-forms $\nu$ on $\Si$. In this case, we show that the tunneling spaces are empty. 

\begin{thmx}\label{Relaive_thm2}
Let $X$ be a smooth, closed, oriented 4-manifold and $\Si\!\subset\!M$ a smooth, closed, oriented surface with positive genus. If $b^{+}_{\sss X-\Si}>0,$ for every $\nu\!\neq\! 0$ and generic adapted perturbation $\eta_\nu,$ $\cM_{\eta_\nu}(X\!-\!\Si,\sfrak)$ is a smooth closed manifold of real dimension
\bEq{DimFormula_e2}
\wt{\dsw}=\frac{(c_1(L)+\Si)^2-2\chi(X)-3\si(X)}{4}-d(\mfs).
\eEq
For different such $\eta_{\nu}$ and $\eta_{\nu'}$, $\cM_{\eta_\nu}(X\!-\!\Si,\sfrak)$ and $\cM_{\eta_{\nu'}}(X\!-\!\Si,\sfrak)$ are smoothly cobordant.
\end{thmx}

\noindent
In this case, for each point $q\in S_{d(\mfs)}(\nu)$, if $\wt{\dsw}$ is even, we define 
\bEq{SWI-Relative_e2}
\wt\SW{}_{\nu}^{\sss X,\Si}(\mfs;q) \defeq \int_{\ev^{-1}(q)\subset\cM_{\eta_\nu}(X-\Si,\mfs)} c^{\wt{\dsw}/2}
\eEq
By the second statement of Theorem~\ref{Relaive_thm2}, the right-hand side (\ref{SWI-Relative_e2}) is independent of the choice of the compactly-supported part of $\eta_\nu$, which justifies the notation on the left-hand side. Recall that, by definition, the elements of $\cM_{\eta_\nu}(X-\Si,\mfs)$ are automatically irreducible. Furthermore, the finite sum
$$
\wt\SW{}^{\sss X,\Si}(\mfs)=
\sum_{q\in S_{d(\mfs)}(\nu)}\wt\SW{}_{\nu}^{\sss X,\Si}(\mfs;q)
$$
is independent of the choice of $\nu.$

\medskip
\par\noindent
The dimension formula (\ref{DimFormula_e2}) differs from (\ref{DimFormula_e}) by a term $\nd(\mfs)$. Particularly, if $X$ is symplectic, identifying the set of $\spinc$ structures $\mfs$ on $TX(-\log \Si)$ with the set of homology classes $\sz\in\!H_2(X,\Z)$ again,  (\ref{DimFormula_e2}) simplifies to 
$$
\wt\dsw=\sz\cdot\sz-(K_X+\Si)\cdot\sz.
$$
Note that $K_X+\Si$ is the logarithmic canonical line bundle of $(X,\Si)$. The integer $\wt\dsw$ is even iff $d(\mfs)=\Si\cdot\sz$ is even. It is not clear to us if $\wt\SW$ is equal to certain count of $J$-holomorphic curves.

\medskip
\par\noindent
Next, we prove that there is a SW sum formula that relates the SW invariant of a sum $X\!=\!X_1\#_\Si X_2$ of two closed oriented 4-manifolds $X_1$ and $X_2$ along a common oriented surface $\Si$ with dual self-intersections to the relative SW invariants of $(X_1,\Si)$ and $(X_2,\Si).$ Here, we only need to consider $\spinc$ structures on $X$ that are of ``pullback type"  on the separating hypersurface $Y$ of this sum manifold. Recall that the SW invariants of the connected sum $X$ of two closed oriented $4$-manifolds $X_1$ and $X_2$ with $b_2^+(X_i)>1$ are zero.

\medskip
\par\noindent
Suppose $X\!=\!X_1\#_\Si X_2$ and $\mfs,\mfs'$ are two $\spinc$ structures on $X$. We say $\mfs$ is equivalent to $\mfs'$ if both restrict to the same $\spinc$ structures on $X_1-\Si$ and $X_2-\Si$. We denote the equivalence class of $\mfs$ by $[\mfs]$.  

\begin{thmx}\label{SW-Sum_thm2}
Suppose $X\!=\!X_1\#_\Si X_2$ and $\mfs$ is a $\spinc$ structure on $X$ that is of pullback type on the separating hypersurface $Y$. For each holomorphic $1$-form $\nu\neq0$ we have the sum formula
\bEq{SW-Sum_e2}
\sum_{\mfs'\in [\mfs]}\SW^{\sss X}(\mfs')=
\sum_{[\mfs]=\mfs_1\#\mfs_2}\; \sum_{q\in S_{d(\mfs)}(\nu)} \ve_{\mfs_1,\mfs_2,q}\,\,
\wt\SW{}^{\sss X_1,\Si}_\nu(\mfs_1;q)\cdot 
\wt\SW{}^{\sss X_2,\Si}_{\nu}(\mfs_2;q).
\eEq
\end{thmx}

\noindent
In (\ref{SW-Sum_e2}), $\ve_{\mfs_1,\mfs_2,q}\!\in\!\{\pm1\}$ depend on the choice of orientations.
Whenever the relative invariants (\ref{SWI-Relative_e}) are defined, we get another sum formula
\bEq{SW-Sum_e}
\sum_{\mfs'\in [\mfs]}\SW^{\sss X}(\mfs')=
\sum_{[\mfs]=\mfs_1\#\mfs_2}\ve_{\mfs_1,\mfs_2}\sum_{\Om}\, 
\SW^{\sss X_1,\Si}(\mfs_1;\Om)\cdot
\SW^{\sss X_2,\Si}(\mfs_2;\Om^*)
\eEq
where $\Om$ runs over the terms in the K\"unneth decomposition $\sum \Om \otimes \Om^*$ of the diagonal in $\tn{Div}(\Si)\!\times\! \tn{Div}(\Si)$. Under the
GW-SW correspondence, the SW sum formula (\ref{SW-Sum_e}) should correspond to the well-known GW sum formula, with $X$ being the symplectic sum of $X_1$ and $X_2$ along $\Si$.

\begin{remark}
If $\mfs$ is a $\spinc$ structure on $X$ that is not of pullback type on the separating hypersurface $Y$, then it follows from Fact 1 in \cite[p.~94]{N1} and the same convergence argument used in the proof of Theorem~\ref{SW-Sum_thm2} that all the terms in (\ref{SW-Sum_e}) are zero.
\end{remark}

\begin{remark}
The issue that the sum $\mfs_1\#\mfs_2$ of two $\spinc$ structures on $X_1$ and $X_2$ is only well-defined up to the equivalence relation defined before Theorem~\ref{SW-Sum_thm2} also appears in the GW sum formula: two homology classes $A_1\!\in\!H_2(X_1,\Z)$ and $A_2\!\in\!H_2(X_2,\Z)$ that have the same intersection number with $\Si$ can be glued together to produce a homology class $A\!\in\!H_2(X,\Z)$; however, $A$ is only well-defined up to addition with homology classes associated to certain elements of $H_1(\Si,\Z)$, known as rim tori. A natural question is whether we can refine the sum formulas (\ref{SW-Sum_e2}) and (\ref{SW-Sum_e}) to express each individual $\SW^{\sss X}(\mfs')$ in terms of some ``refined relative invariants" of $X_1$ and $X_2$. In the context of GW sum formula, this problem is extensively studied in \cite{FZ1,FZ2}. 
\end{remark}

\section{A tour of Seiberg-Witten theory}\label{setup_s}
\subsection{\Spinc structures}
Given an oriented riemannian rank $2n$ (real) vector bundle $V\!\lra\! X$ over a smooth manifold $X,$ a \spinc structure $\sfrak\!=\!(\sb,\clm)$ 
on $V$ consists of a rank $2^n$ (complex) hermitian vector bundle $\sb\!\lra\!X,$ called the spinor bundle, and a Clifford multiplication, which is a linear map of vector bundles
$$
\clm: V\lra\mend(\sb),\qquad V_x\!\ni v \lra \clm(v)\in\mend(\sb|_{x}),
$$ 
satisfying 
\bEq{SpincMap_e}
\clm(v)^*+\clm(v)=0, \qquad \clm(v)^*\clm(v)=|v|^2~\tn{id}, \qquad \forall~v\!\in\!V.
\eEq
We sometimes denote the action of $\clm(v)$ on $\psi\in\sb$ by $v\cdot\psi$ when there is no chance of confusion. 
Every such $\sb$ admits a canonical splitting into rank $2^{n-1}$ hermitian vector bundles, $\sb\!=\!\sb^+\oplus \sb^-$, satisfying
$$\clm(v)\sb^{\pm}|_{x}=\sb^{\mp}|_{x}, \qquad \forall~x\!\in\!X,~v\!\in\!V_x.$$
In other words, we have
$$\clm(v)=
 \left[ 
  {\begin{array}{cc}
   0 & \gamma(v) \\
   -\gamma(v)^* & 0 \\
  \end{array}} 
 \right],\qquad 
\gamma: V\lra\tn{Hom}(\sb^-,\sb^+),$$
where 
$$\gamma(v)^*\gamma(v)=|v|^2~\tn{id}.$$
A \spinc isomorphism from $(\sb_0, \clm_0)$ to $(\sb_1, \clm_1)$ is a unitary bundle isomorphism $f:\sb_0\lra\sb_1$ which induces a bundle isomorphism $\mend(f):\mend(\sb_0)\lra\mend(\sb_1)$ such that $\mend(f)\circ\clm_0=\clm_1.$
Denote by $\mc{S}^c(V)$ the set of isomorphism classes of $\tn{spin}^c$ structures on $V$. Given a riemannian $2n$-manifold $X,$ let $\mc{S}^c(X)$ be the set of isomorphism classes of $\tn{spin}^c$ structures on $TX.$  This is in fact the set of principal \spinc bundles lifting the principal tangent bundle of $X$ up to bundle isomorphism. 

\medskip
\par\noindent
If $(\sb,\clm)$ is a \spinc structure on $V$ and $n\!>\! 1,$ then 
$$\det(\sb^+)=\det(\sb^-)=L_{\mfs}^{\otimes 2^{n-2}}$$
for a unique complex line bundle $L_{\mfs},$ called the {\it characteristic line bundle} of the \spinc structure, 
for which $c_1(L_{\mfs})$ is an integral lift of $w_2(V).$ 
In fact, $L_\sfrak$ is the determinant line bundle of the principal $\tn{Spin}^c(2n)$-bundle which lifts the principal $\tn{SO}(2n)$-bundle associated to $V\!\lra\! X.$
We define $c_1(\mfs)\defeq c_1(L_{\mfs}).$ Similarly, we say that the \spinc structure $\mfs$ is torsion if $c_1(L_{\mfs})$ is torsion. 
When there is no chance of confusion, we will drop the subscript $\mfs$ and simply write $L.$

\medskip
\par\noindent
Given a $\tn{spin}^c$ structure $\mfs_0\!=\!(\sb_0, \clm_0)$, every other \spinc structure $\mfs\!=\!(\sb,\clm)\!\in\!\mc{S}^c(V)$ has the form 
$$
\sb\cong \sb_0\otimes E,\qquad \clm=\clm_0\otimes \tn{id},
$$ 
where $E$ is a hermitian line bundle. We express the above relation between $\sfrak$ and $\sfrak_0$ by writing $\sfrak=\sfrak_0\otimes E.$
Note that $L_{\mfs}=L_{\mfs_0}\otimes E^{\otimes{2}}$ and the two \spinc structures $\sfrak$ and $\sfrak_0$ are isomorphic if and only if $c_1(E)\!=\! 0.$

\medskip
\par\noindent
Every complex vector bundle $(V,J)$ equipped with a hermitian metric $\ll\,,\,\rr$ admits a canonical \spinc structure $\sfrak_{\tn{can}}=(\sb_{\tn{can}},\clm_{\tn{can}})$ with
\bEq{CanW_e}
\aligned
& 
\sb_{\tn{can}}^+=\Lambda^{0,\tn{even}}V^* \qquad \sb_{\tn{can}}^-=\Lambda^{0,\tn{odd}}V^*, \\
&
\clm_{\tn{can}}(v)\alpha = -\sqrt{2}\iota_v\alpha + \frac{1}{\sqrt{2}} \ll\cdot,v\rr \wedge\alpha.
\endaligned
\eEq
As a result, the correspondence $\sfrak_{\tn{can}}\otimes E \leftrightarrow E$ determines a canonical bijection between isomorphism classes of \spinc structures and isomorphism classes of complex line bundles over $X.$

\medskip
\par\noindent
Specifically, every almost-complex manifold $(X^{2n},J)$ admits a canonical \spinc structure $\sfrak_{\tn{can}}$ with 
$$
\aligned
&\sb_{\tn{can}}^+=\Lambda^{0,\tn{even}}T^*X, \quad \sb_{\tn{can}}^-=\Lambda^{0,\tn{odd}} T^*X, \quad L_{\sfrak_{\tn{can}}}=K_X^*,\\
&\clm_{\tn{can}}(v)\alpha=-\sqrt{2}\iota_v\alpha+\frac{1}{\sqrt{2}} \ll\cdot,v\rr \wedge\alpha.
\endaligned
$$
In particular, if $(X,\om_X)$ is a symplectic manifold, we can choose an $\om_X$-compatible (or $\om_X$-tame) almost-complex structure $J$ and define $\sfrak_{\tn{can}}$ as above. In this case, $K_{X}=\Lambda^{n,0}T^*X$ is the canonical bundle of $X$ with respect to $\om_X$ and $\ll\cdot,v\rr=\om_X(\cdot,\mfi v + Jv)$ is the $(0,1)$-form dual to $\mfi v\!+\!Jv\in T^{1,0}X.$
If $J$ is $\om_X$-tame but not $\om_X$-compatible, then one has to replace $\om_X$ in the definition of $\ll\cdot ,v\rr$ with
$$\wt\om_X (u,v)=\frac{1}{2}\big(\om_X(u,v)+\om_X(Ju,Jv)\big).$$

\medskip
\par\noindent
Let \gmetric be a riemannian metric on $X^{2n},$ $\nabla$ the corresponding Levi-Civita connection, and $\sfrak\!\in\!\mc{S}^c(X).$ 
A hermitian connection $\wt\nabla$ on $\sb$ is called a \spinc connection if it is compatible with $\nabla$ and the Clifford multiplication
$$\wt\nabla_v(w\cdot\Phi)=w\cdot\wt\nabla_v\Phi+(\nabla_v w)\cdot\Phi, \qquad \forall~\Phi\in\Ga(\sb),~v,w\!\in\Ga(TX).$$ 
Every $\spinc$ connection preserves $\sb^\pm.$ Every two such connections differ by an imaginary-valued $1$-form $\alpha\!\in\!\Om^1(X,\mfi\R).$ Moreover, the gauge group $\G=\tn{Maps}(X,S^1)$ acts on the space of connections by 
$$(u^*\wt\nabla)\Phi=u^{-1}\wt\nabla(u\Phi)=\frac{\nd u}{u}\otimes\Phi+\wt\nabla\Phi$$
for $\Phi\in\Ga(\sb)$ and $u\!\in\!\tn{Maps}(X,S^1).$

\medskip
\par\noindent
Every \spinc connection $\wt\nabla$ on the spinor bundle $\sb$ is uniquely determined by $\nabla$ and a connection $A$ on the characteristic line bundle $L$ of the \spinc structure. 
This is essentially due to the fact that the underlying principal bundle of $\sb$ (with fiber $\Spinc(2n)$) lifts the principal tangent bundle of $X$ (with fiber SO($2n$)) as a circle-bundle extension (corresponding to the line bundle $L$).
This corresponds to the lifting of $w_2(V)$ to the integral class $c_1(L);$ see \cite[Chap.~3]{m}.   
The space of connections on the characteristic line bundle $L,$ denoted by $\mc{A}(L),$ is an affine space with tangent space $\Om^1(X,\mfi\R).$ 
If $A\!\in\!\mc{A}(L)$ and $a\!\in\!\Om^1(X,\mfi\R),$ then the curvature $F_A\!\in\!\Om^2(X,\mfi\R)$ and $F_{A+a}\!=\!F_A+\nd a.$
The connection $A$ on $L$ induces a \spinc connection $\wt\nabla=\nabla_A$ on $\sb$ 
and their curvatures are related by \cite[Sec.~6.1]{DS}
$$F_A(v,w)=\frac{1}{2^n}\tn{tr}\big(F_{\nabla_A}(v,w)\big).$$
\par\noindent
The gauge group $\G=\tn{Maps}(X,S^1)$ acts on $\mc{A}(L)$ by 
$$u^* A= \frac{\nd u}{u}+ A$$
and leaves $F_A$ invariant. The operator
$$
\dirac_A\colon \Ga(\sb^\pm)\lra \Ga(\sb^\mp), \qquad 
\dirac_A\Phi=\sum_{{i=1}}^{2n}e_i\cdot\wt\nabla_{e_i}\Phi,
$$
where $e_1,\dots,e_{2n}$ form an orthonormal frame for $TX,$ is the Dirac operator associated to $A\!\in\!\mc{A}(L)$. It is a self-adjoint operator independent of the particular choice of $e_1,\dots,e_{2n}\,.$

\medskip
\par\noindent
If $(X,\om_X,J)$ is K\"ahler, then 
$$\Ga(\sb_{\tn{can}}^+\otimes E)=\Om^{0,\tn{even}}(X,E), \qquad \Ga(\sb_{\tn{can}}^-\otimes E)=\Om^{0,\tn{odd}}(X,E),$$
and 
\bEq{SympDirac_e}
\frac{1}{\sqrt{2}}\dirac_A=\dbar_{A_E} +\dbar^*_{A_E} :\Om^{0,\tn{even}}(X,E)\lra  \Om^{0,\tn{odd}}(X,E),
\eEq
where $A_E$ is a hermitian connection on $E,$ $\dbar_{A_E}=\nabla_{A_E}^{0,1}$ is the $\dbar$-operator associated to $A_0,$ and $A$ is the induced connection on the characteristic line bundle $L$ of the \spinc structure.
For example, in dimension 4, $\sb_{\tn{can}}^-=\Lambda^{0,1} T^*X,$ $\det(\sb_{\tn{can}}^-)=K_X^*$ and 
$L=\det(\sb^-)=\det(\sb_{\tn{can}}^-\otimes E)=K_X^* \otimes E^2.$
Therefore, any hermitian connection $A$ on $L$ is equivalent to a hermitian connection $A_E$ on $E$ and we have $A=A_{\tn{can}}\otimes A_E^2,$ where $A_{\tn{can}}$ is the holomorphic hermitian connection on $K_X^*.$

\medskip
\par\noindent
The identity (\ref{SympDirac_e}) continues to hold in the symplectic case \cite[Thm.~6.17]{DS}. 
If $(X,\om_X)$ is a symplectic manifold, $J$ is compatible with $\om_X$ (resp. tames $\om_X$), and $\nabla$ is the Levi-Civita connection of the metric  $\gmetric\!=\!\om_X(\cdot, J\cdot)$ (resp. $\gmetric\!=\!\wt\om_X(\cdot, J\cdot)$), then $\nabla$ does not preserve $\Om^{0,k}(X)$ unless $(X,\om_X,J)$ is K\"ahler, i.e., $\nabla J=0$. However, there is a canonical hermitian connection on $TX,$ defined by 
$$\wt\nabla_v w\defeq\nabla_v w -\frac{1}{2} J(\nabla_vJ) w,$$
which gives rise to a hermitian connection on $\sb_{\tn{can}}=\Lambda^{0,*}T^*X,$ compatible with the Clifford multiplication but not with $\nabla.$ 
However, it is possible to modify $\wt\nabla$ further to produce a \spinc connection $\wh\nabla,$ compatible with both the Clifford multiplication and the Levi-Civita connection $\nabla$
$$\wh\nabla_v \Phi\defeq\wt\nabla_v \Phi +\frac{1}{2} \mu(J\nabla_v J) \Phi,$$
where $\Phi\in\Ga(\sb_{\tn{can}})=\Om^{0,*}(X)$ and 
$\mu\colon\mf{so}(TX)\lra\tn{End}(\sb_{\tn{can}})$ is characterized by
$$[\mu(A),\clm_{\tn{can}}(v)]=\clm_{\tn{can}}(Av).$$
When $X$ is K\"ahler, the connections $\wh\nabla$ and  $\wt\nabla$ coincide with the Levi-Civita connection $\nabla$ on forms \cite[pp.~198--199]{DS}.
%
\subsection{SW equations in dimension 4}\label{SW4_ss}
\par
Let $X$ be a smooth closed connected oriented riemannian $4$-manifold (with metric $\gmetric$) and $\sfrak=(\sb,\clm)$ a $\tn{spin}^c$ structure on $X.$ The (unperturbed) Seiberg-Witten monopole equations are a system of first order differential equations for a pair $(A,\Phi)$ in the configuration space $\cC(X,\sfrak)\!=\!\mc{A}(L)\!\times\!\Gamma(\sb^+),$ where $A$ is a connection on the characteristic line bundle of $\sfrak$ and $\Phi$ is a plus-spinor. 
The spaces $\mc{A}(L)$ and $\Ga(\sb^+)$ are completed with respect to appropriate Sobolev norms, so that we will be working in the context of Banach spaces; see \cite{m}. 
The Seiberg-Witten equations read
$$F_A^+=(\Phi\Phi^*)_0, \qquad \dirac_A\Phi=0, \eqno{(\tn{SW}^4)}$$
where $(\Phi\Phi^*)_0\in\Ga(\tn{End}_0(\sb^+)),$ defined by
$$(\Phi\Phi^*)_0(w)=\ll\Phi,w\rr \Phi-\frac{1}{2}|\Phi|^2 w \qquad \forall~w\!\in\!\Ga(\sb^+),$$
is the trace-less part of $\Phi\Phi^*\!\in\!\Ga(\tn{End}(\sb^+))$ and the two sides of the curvature equation in ${\tn{SW}^4}$ are identified via the bundle isomorphism
\bEq{CM4_e}
\Lambda^{2,+}T^*X\otimes\C \lra \tn{End}_0(\sb^+), \qquad 
\sum_{i<j} c_{ij}e_i^*\wedge e_j^* \lra \sum_{i<j} c_{ij}\clm(e_i)\clm(e_j).
\eEq
We will occasionally denote the self-dual 2-form representing the quadratic term $(\Phi\Phi^*)_0$ by $q(\Phi).$

\medskip
\par\noindent
For an imaginary-valued self-dual $2$-form $\eta\!\in\!\Om^{2,+}(X,\mfi\R),$ the perturbed Seiberg-Witten monopole equations are
$$
F_A^+-\eta=(\Phi\Phi^*)_0, \qquad \dirac_A\Phi=0. 
\label{SW4eta}\eqno{(\tn{SW}^4_\eta)}
$$
We call a solution to $\tn{SW}^4_\eta$ an $\eta$-monopole. 
The set of $\eta$-monopoles is invariant under the action of $\mc{G}=\tn{Maps}(X,S^1),$ which sends $(A,\Phi)$ to $(u^*A,u^{-1}\Phi).$ Here, $\mc{G}$ is of course completed with an appropriate Sobolev norm consistent with the completions of $\mc{A}(L)$ and $\Ga(\sb^+).$ 
Let
$$
\Om^{2,+}(X,\mfi \R)_0=\{\eta\!\in\!\Om^{2,+}(X,\mfi \R)\colon~\exists A\!\in\!\mc{A}(L)~\tn{s.t.}~F_A^+-\eta=0\}.
$$
This is an affine subspace of codimension $b^+$ in $\Om^{2,+}(X,\mfi \R).$ 
An $\eta$-monopole for which the spinor $\Phi$ is identically zero is called reducible, in which case $\eta\!\in\!\Om^{2,+}(X,\mfi \R)_0,$ and irreducible otherwise.
The Seiberg-Witten moduli space is the quotient
$$
\cM(X,\sfrak)=\cM_{\gmetric,\eta}(X,\sfrak)\defeq
\{(A,\Phi)\in\mc{A}(L)\times\Ga(\sb^+)~\tn{satisfying}~\tn{SW}^4_\eta\}/\mc{G}. 
$$
If $b^+>0,$ we can avoid reducible solutions of $\tn{SW}^4_\eta$ with a generic perturbation $\eta$ and obtain a smooth closed orientable moduli space of real dimension
\bEq{DimClassic4_eq}
\dsw=\frac{c_1(L)^2-2\chi(X)-3\si(X)}{4}.
\eEq
The deformation-obstruction complex of $\cM(X,\sfrak)$ at any $\eta$-monopole $(A,\Phi)$ is given by 
$$
0\longrightarrow
\Om^0(X;\mfi\R)
\stackrel{\cD^0}{\longrightarrow}
\Om^1(X;\mfi\R)\oplus\Gamma(\sb^+)
\stackrel{\cD^1}{\longrightarrow}
\Om^{2,+}(X;\mfi\R)\oplus\Gamma(\sb^-)
\longrightarrow 0.
\eqno{(\cE_X(A,\Phi))}
$$
The Banach spaces in $\cE_X(A,\Phi)$ are tangent spaces to $\mc{G},$ $\mc{A}(L)\times\Gamma(\sb^+)$ and $\Om^{2,+}(\Si;\mfi\R)\times\Gamma(\sb^-)$ with appropriate completions, respectively, and the maps 
$$\cD^0(f)=(2\,\nd f, -f\Phi) \quad , \quad
\cD^1(a,\phi)=\Big(\nd^+ a-D_\Phi q(\phi),\frac{1}{2} a\cdot\Phi+\dirac_A\phi\Big)$$
are, respectively, the linearizations of the gauge group action and the map 
$$\tn{SW}_{\eta}^4(A,\Phi)=\Big(F_A^+-q(\Phi)-\eta,\, \dirac_A\Phi\Big).$$ 
We have used small letters $a$ and $\phi$ to distinguish tangent space variables at $(A,\Phi)$ from the point itself; we may occasionally adopt this convention later on.
\\
The complex $\cE_X(A,\Phi)$ is an elliptic complex and its index is $-\dsw.$ In fact, $H^0(\cE_X(A,\Phi))=\ker({\cD^0})$ is zero exactly when $\Phi\neq0$ at least at some point, which is the definition of $(A,\Phi)$ being irreducible. The first cohomology $H^1(\cE_X(A,\Phi))=\ker({\cD^1})/\Im({\cD^0})$ is the Zariski tangent space of the moduli space $\cM(X,\sfrak)$ and $H^2(\cE_X(A,\Phi))=\coker({\cD^1})$ is its obstruction space, which is zero when $(A,\Phi)$ is a regular point of the moduli space.

\medskip
\par\noindent
The Seiberg-Witten invariants $\tn{SW}(X,\sfrak)$ are defined as the integral of a canonical cohomology class (a suitable power of the $2$-form $c,$ as in (\ref{SWI-Relative_e})) over $\cM(X,\sfrak),$ if this moduli space is even-dimensional and non-empty, and zero otherwise. The invariants are well-defined (independent of the choice of the metric $\gmetric$ and the generic perturbation term $\eta$) when $b^+>1.$ 
If $b^+=1,$ the formal dimension $\dsw$ of the SW moduli space in (\ref{DimClassic4_eq}) is even exactly when $b_1$ is even, in which case 
the SW invariants $\tn{SW}^{\pm}(X,\sfrak)$ depend on the chamber where the metric $\gmetric$ resides and satisfy a wall-crossing formula
$$\tn{SW}^{+}(X,\sfrak)-\tn{SW}^{-}(X,\sfrak)=-(-1)^{\dsw/2}.$$
\par\noindent
If $(X,\om_X,J)$ is K\"ahler, using (\ref{SympDirac_e}), the $\tn{SW}^4_\eta$ equations for the pair $(A,\Phi)$ 
take the form
$$
2(2F_{A_E}-\eta)^{0,2}=\ov{\Phi_0}\Phi_2, \quad 
4\mfi(F_{A_{\tn{can}}}+2F^+_{A_E}-\eta)^{1,1}=(|\Phi_2|^2-|\Phi_0|^2)\om_X, \quad 
\dbar_{A_E}\Phi_0+\dbar_{A_E}^*\Phi_2=0,
$$
where 
$$
A=A_{\tn{can}}\otimes A_E^2, \qquad
\Phi=(\Phi_0,\Phi_2)\in \Om^{0,0}(X,E)\times \Om^{0,2}(X,E).
$$
If $\eta\!\in\!\Om^{1,1}\cap\Om^{2,+}$, then either $\Phi_0=0$ or $\Phi_2=0$. The latter will happen if 
$$
c_1(E)\cdot \om_X <\frac{1}{2} c_1(K_X)\cdot \om_X,
$$
or $\eta$ has a large (positive) multiple of $-\mfi\om_X$; see \cite[Sec.~12.2]{DS}. 
Recall that $c_1(K_X)\cdot\om_X\geq 0$ for every K\"ahler surface with $b^+>1.$
By adapting the arguments in \cite[Sec.~7.2]{m}, 
we obtain a holomorphic description of the moduli space of $\eta$-monopoles in terms of holomorphic structures on $E$ (or, equivalently, on the characteristic line bundle $L$ of the \spinc structure) and non-zero holomorphic sections of $E$ or $K_X \otimes E^*$ (up to constant scalar multiples), depending on the case. 

\medskip
\par\noindent
If $\Phi_2=0$, the unperturbed SW equations reduce to the vortex equations 
$$
F_{A_E}^{0,2}=0,\quad  
4\mfi(F_{A_{\tn{can}}}+2F_{A_E})^{1,1}=-|\Phi_0|^2\om_X,\quad 
\dbar_{A_E}\Phi_0=0. 
$$ 
Similar results hold in the symplectic case. In the almost-complex case, the $\tn{SW}^4_\eta$ equations for the pair $(A_{\tn{can}}\otimes A_E^2, \Phi)$ take the form
$$
(F_{A_{\tn{can}}}+2F_{A_E})^+-\eta=q(\Phi),\qquad
\dbar_{A_E}\Phi_0+\dbar_{A_E}^*\Phi_2=0.
$$

\subsection{SW equations in dimension 2}\label{SW2_ss}
\par
In this and the following sections, suppose $(\Si,\mfj,\om)$ is a closed Riemann surface of genus $g,$ equipped with a complex structure $\mfj$ and a K\"ahler form $\om.$ As a K\"ahler manifold of complex dimension one (real dimension two), $\Si$ carries a canonical \spinc structure $(\ucan,\clmcan),$ where $\ucan$ is a hermitian vector bundle of rank 2, 
and any other \spinc structure $(\ub,\clm)$ on $\Si$ is obtained by twisting this canonical \spinc structure by a hermitian line bundle $E.$ The characteristic line bundle is $L=\det(\ub)$ and the spinor bundle $\ub=\ub^+\oplus\ub^-$ splits as a direct sum of two complex hermitian line bundles. We have
$$\ucan^+=\bbC_{\sss \Si},\qquad \ucan^-=\kcan^*, \qquad\tn{and}~~~\lcan=\kcan^*,$$
where $\kcan=\Lambda^{1,0}T^*\Si$ is the canonical line bundle of $\Si$ and $\bbC_{\sss \Si}$ is the trivial line bundle, and
$$\ub^+=E,\qquad \ub^-\!=\!\kcan^*\otimes E, 
\quad\tn{and}\quad L=\kcan^*\otimes E^2.$$
When there is a chance of confusion, we may add a subscript and, for example, denote the spinor bundles $\ub^\pm$ or the line bundle $L$ by $\ub^\pm_{\sss E}$ or $L_{\sss E},$ respectively.
\begin{remark}\label{DwN_rmk}
It is also possible to define a canonical {\it spin} structure on $\Si$ by simply lifting its principal $\tn{SO}(2)$-tangent-bundle to a $\tn{Spin}(2)$-bundle as a
non-trivial double cover. One can then tensor with $\C$ to obtain a \spinc structure on $\Si.$ This is the approach taken in \cite[p. 93]{N1}. The spinor bundles obtained in this case are $\ub^-=\kcan^{-1/2}$ and $\ub^+=\kcan^{1/2}$ and the characteristic line bundle is trivial. This complexification of the canonical spin structure is obviously different from our canonical \spinc structure above, which is induced by the (almost-) complex structure on $\Si,$ but it can easily be offset by an additional twist by the line bundle $\kcan^{-1/2}=\sqrt{\kcan^*}.$
We have chosen here to work with the canonical \spinc structure, rather than the canonical spin structure, because of its compatibility with the overall framework of complex structures. 
\qed
\end{remark}

\medskip
\par\noindent
Now, consider the set of pairs $(A,\Psi)\!\in\!\mc{A}(L)\times\Gamma(\ub)$ satisfying the equations
$$F_A=\frac{|\Psi_+|^2-|\Psi_-|^2}{2}\,\mfi\om,\qquad \dirac_A\Psi=0, \eqno{(\tn{SW}^2)}$$
where $\Psi=(\Psi_+,\Psi_-)$ is a section of $\ub=\ub^+\oplus\ub^-.$ 
The set of solutions to these equations is invariant under the action of the gauge group $\mc{G}=\tn{Maps}(\Si,S^1)$ on $\mc{A}(L)\times\Gamma(\ub).$
As before, we call a solution reducible if the spinor $\Psi$ is identically zero.
\medskip
\par\noindent
If $A_{\tn{can}}$ is the holomorphic hermitian connection on $\lcan\!=\!\kcan^*,$ then any connection on $L\!=\!\kcan^*\otimes E^2$ is of the form $A\!=\!A_{\tn{can}}\otimes A_E^2,$ where $A_E\!\in\!\mc{A}(E),$ and the corresponding Dirac operator $\dirac_A$ on sections of $\ub=\ub^+\oplus\ub^-$ takes the form
$$
\frac{1}{\sqrt{2}}\dirac_A=
\left[ 
  {\begin{array}{cc}
   0 & \dbar_{A_E}^* \\
   \dbar_{A_E} & 0 \\
  \end{array}} 
\right].
$$
The Dirac equation then splits into two equations $\dbar_{A_E}\Psi_+ =0$ and $\dbar_{A_E}^*\Psi_-=0$,
and the curvature equation turns into 
$$F_{A_{\tn{can}}}+2F_{A_E}=\frac{|\Psi_+|^2-|\Psi_-|^2}{2}\,\mfi\om.$$ 
Since $\dbar_{A_E}\dbar_{A_E}=0$ and $\dbar_{A_E}\Psi_+=0,$ we conclude that $A_E$ induces a holomorphic structure on $E$ for which $\Psi_+$ is a holomorphic section. Similarly, using the involution on the \spinc structures induced by complex conjugation, which sends the characteristic line bundle to its inverse, we can see that the equation $\dbar_{A_E}^*\Psi_-=0$ implies that $\overline{\Psi}_-$ is a holomorphic section of $\kcan\otimes E^*,$ the Serre dual of~$E.$ Since line bundles of negative degree can not have non-zero holomorphic sections, we conclude that one of $\Psi_+$ or $\Psi_-$ is identically zero, unless $0\leq d=\deg(E)\leq 2g-2.$

\medskip
\par\noindent
We will also be interested in another variant of the Seiberg-Witten equations
$$F_A=\frac{|\Psi_+|^2-|\Psi_-|^2}{2}\,\mfi\om,\qquad \dirac_A\Psi=0, \quad\tn{and}\quad \Psi_+\ov\Psi_-=\mfi\nu, \eqno{(\tn{SW}^2_\nu)}$$
where $\nu$ is a fixed holomorphic 1-form. Notice that $\Psi_+\ov\Psi_-$, as a section of $\ub^+\otimes(\ub^-)^*=\kcan$, identifies with a holomorphic 1-form on $\Si.$ This variant of the SW equations is closely related to the perturbed SW equations in higher dimensions; see Lemma \ref{Nic-lm}.

\subsection{SW equations on a circle bundle over a Riemann surface}\label{SW3_ss}
Similarly to $2$-manifolds, for a closed $3$-manifold $Y,$ a \spinc structure $\mfs=(\sb,\clm)$ consists of a rank 2 hermitian bundle $\sb$ and a Clifford multiplication map $\clm\colon TY\lra\mend(\sb)$ satisfying (\ref{SpincMap_e}), which result from a principal $\tn {Spin}^c(3)$-lifting of the principal frame bundle of $TY.$ Note that, unlike dimensions 2 and 4, the spinor bundle $\sb$ does not decompose into plus- and minus-spinor bundles. 
\par\noindent
We are interested in the case where $Y$ is a circle bundle over a Riemann surface $\Si.$ In this subsection, to simplify notations, we will temporarily use the shorthand $\piysi$ to denote the bundle projection
$$\piysi=\pi_{\sss Y,\Si}\colon Y\!\lra\!\Si.$$ 
Viewing $Y$ as a principal $U(1)$-bundle, a choice of a principal $U(1)$-connection gives rise to a decomposition
\bEq{Splitting_e}
TY\cong T^{\tn{ver}}Y \oplus T^{\tn{hor}}Y,
\eEq
where the vertical summand is tangent to the $S^1$ fibers $\ker(\nd\piysi)$ and 
the horizontal summand $T^{\tn{hor}}Y\!\cong\!\piysi^*T\Si$ is a lifting of the tangent space of $\Si$ which is equivariant under the $U(1)$ action.
We will denote this connection by $\mfi\al,$ where $\al$ is a $U(1)$-invariant real $1$-form on $Y$ whose restriction to each fiber is $\nd\theta.$

\medskip
\par\noindent
Under such a decomposition, $Y$ admits a canonical $\spinc$ structure $(\sycan,\clmycan),$ where
$$\sycan=\bbC_{\sss Y}\oplus\piysi^*(\kcan^*)$$
is the pullback of the canonical spinor bundle $\ucan=\bbC_{\sss \Si}\oplus\kcan^*$ on $\Si,$ the Clifford map $\clmycan$ on 
$T^{\tn{ver}}Y$ sends the generator $\partial_\theta$ to 
$$
\begin{bmatrix} 
\mfi & 0 \\
0 & -\mfi 
\end{bmatrix},
$$
and $\clmycan$ on
$T^{\tn{hor}}Y$ is the pullback of the canonical $\clm_{\sss \Si,\tn{can}}$ as in (\ref{CanW_e}).
The isomorphism class of this \spinc structure is independent of the choice of the principal $U(1)$-connection and lifts the canonical Spin$^c(2)$ structure on $\Si$ to the canonical Spin$^c(3)$ structure on $Y.$    
Every other \spinc structure on $Y$ is obtained by tensoring with a complex line bundle on $Y.$ 
\begin{remark}\label{Gysin_rmk}
By the Gysin long exact sequence
$$H^2(\Si,\Z)\lra H^2(Y,\Z)\lra H^1(\Si,\Z) \lra 0,$$ 
we see that, unless $\Si=S^2$, there are complex line bundles on $Y$ that are not the pullback of a complex line bundle on~$\Si$. 
In fact, the Gysin exact sequence implies that if the Euler number $\ell$ of the circle-bundle $Y\!\lra\!\Si$ is non-zero, then 
$$H^2(Y,\Z)\cong\Z^{2g}\oplus\big(\Z/\ell\Z\big),$$
where the torsion part is generated by the pullback $\piysi^*\om$ of the volume form on $\Si.$
Therefore, the SW theory on $Y$ involves a larger class of \spinc structures than those on $\Si.$ In the setup that we will be developing for relative SW theory in Section~\ref{LogStr_ss}, we will be working with complex line bundles that are defined over the entire $X$. Therefore we will be dealing only with complex line bundles on $Y$ that are pullbacks of those on $\Si.$ Unless the circle bundle $Y\!\lra\!\Si$ is trivial, these line bundles on $Y$ will always be torsion. 
Moreover, it follows from \cite[p.~94, Fact~1]{N1} that $\spinc$ structures on $Y$ that are not of pullback type result in trivial relative invariants, as we will see in Section~\ref{LogStr_ss}.\qed
\end{remark}
\begin{remark}\label{PullBackSpinc_rmk}
If the degree $\ell=\deg(Y)$ of the $U(1)$-bundle $Y\!\lra\!\Si$ is not zero and $E$ and $E'$ are two complex line bundles on $\Si$, then 
$$
\piysi^*(E)\cong\piysi^*(E')\quad \tn{if and only if}\quad
\deg(E)\equiv\deg(E')\quad \tn{modulo}~\ell.
$$
Unless $\deg(E)\!=\!\deg(E')$, such an isomorphism does not extend to the disk bundle $\cD$ over $\Si$ which has $Y$ as its boundary, where $\cD$ is defined using the same $U(1)$-cocycles of $Y\!\lra\!\Si,$ only with circle fibers replaced by unit disks. \qed
\end{remark}

\medskip
\par\noindent
The general setup of Seiberg-Witten equations in dimension 3 is in most respects analogous to those in dimensions 2 and 4. For a closed, oriented, riemannian 3-manifold $Y,$ equipped with a \spinc structure $\sfrak\!=\!(\sb,\clm),$ consider the following (perturbed) Seiberg-Witten equations for a pair $(B,\Psi),$ consisting of a connection on the characteristic line bundle $L=\det(\sb)$ and a section of the spinor bundle $\sb$ 
$$
F_B-\eta=(\Psi\Psi^*)_0, \qquad \dirac_B\Psi=0,
\label{SW3eta}\eqno{(\tn{SW}^3_\eta)}
$$
where $\eta\!\in\!\Om^{2}(Y,\mfi\R)$ is a closed $2$-form and the two sides of the curvature equation are identified via the isomorphism
\bEq{CM_e}
\Lambda^{2}T^*Y\otimes\C \stackrel{\cong}{\lra} \tn{End}_0(\sb), \qquad 
\sum_{i<j} c_{ij}e_i^*\wedge e_j^* \lra \sum_{i<j} c_{ij}\clm(e_i)\clm(e_j),
\eEq
which is the adjoint of the Clifford multiplication.

\medskip
\par\noindent
As in Sections~\ref{SW4_ss} and \ref{SW2_ss}, we consider the following deformation-obstruction complex, as described in \cite{mst}, to understand the infinitesimal structure of the $\tn{SW}^3_\eta$ moduli space at each solution $(B,\Psi),$
\glossary{$\cE_Y(B,\Psi)$}%
$$
0\lra\Om^0_2(Y;\mfi\R)\stackrel{\cD^0}{\lra}
\Om^1_1(Y;\mfi\R)\oplus\Gamma_1(\sb)\stackrel{\cD^1}{\lra}
\Om^1_0(Y;\mfi\R)\oplus\Gamma_0(\sb)/\cD^0(\Om^0_1(Y;\mfi\R))\lra 0,
\eqno{(\cE_Y(B,\Psi))}
$$
where a subscript $k$ indicates that the space has been completed using the $L^2_k$ norm, $\cD^0$ is the linearization of the gauge group action, and $\cD^1$ is induced by the linearization of
$$\tn{SW}^3_{\eta}(B,\Psi)=\Big(\!*\!\big(F_B-\eta-q(\Psi)\big),\dirac_B\Psi\Big).$$ 
The cohomologies of this complex, as before, determine irreducibility, the Zariski tangent space and the obstruction space at a solution $(B,\Psi),$ respectively. 
Moreover, as is the case for any elliptic complex on an odd-dimensional manifold, the index of $\cE_Y(B,\Psi)$ is zero.

\medskip
\par\noindent
In what follows, we are going to describe how solutions of the $\tn{SW}^3_\eta$ equations on a circle bundle $\piysi\colon Y\!\lra\!\Si$ are related to those of the $\tn{SW}^2_\nu$ equations on $\Si,$ if the \spinc structure $\syfrak$ on $Y$ is of pullback type and $\eta$ is obtained from $\nu$ in a particular way. Such a \spinc structure is specified by 
tensoring the canonical \spinc structure on $Y$ with
a line bundle $\piysi^*E,$ where $E$ is a hermitian line bundle of degree $d$ on $\Si.$ However, as we have seen in remark \ref{PullBackSpinc_rmk}, any other line bundle on $\Si$ whose degree is $d$ modulo $\ell\!=\!\deg(Y)$ will also pull back to the same $\piysi^*E$ on $Y.$
Therefore, given a fixed \spinc structure $\syfrak$ on $Y,$ the choice of $E$ on $\Si$ is not unique. We will indicate this ambiguity by a subscript $[E]$ in our notation; for example, $\mfsy=\mfsyeclass$ is meant to emphasize that the \spinc structure depends on the equivalence class $[E]$ under the relation given in remark~\ref{PullBackSpinc_rmk}.

\medskip
\par\noindent
With this notation in place, fix the \spinc structures $\mfssie$ on $\Si$ and $\mfsyeclass$ on $Y$ as above. These are obtained by twisting the canonical \spinc structures on $\Si$ and $Y$ by $E$ and $\piysi^*E,$ respectively.
Then the spinor bundle $\ub=\ub_{\sss E}$ on $\Si$ pulls back to the spinor bundle $\sb=\sb_{{\sss [E]}}$ on $Y$ and we have 
\bEq{U3_e}
\begin{array}{rccccrl}
\ub=&\negsp\ub^+\oplus\ub^-&\negsp=&\negsp E\oplus(\kcan^*\otimes E) 
\qquad & \tn{and} & \qquad L=L_{\sss \Si,E}&\negsp=\kcan^*\otimes E^{\otimes 2}
\\
\sb=&\negsp\piysi^*\ub^+\oplus\piysi^*\ub^-&\negsp=&\negsp \piysi^*E\oplus\piysi^*(\kcan^*\otimes E)
\qquad & \tn{and} & \quad L_{\sss Y,[E]}&\negsp=\piysi^*L,
\end{array}
\eEq
where $L$ and $\piysi^*L$ are the characteristic line bundles of the \spinc structures on $\Si$ and $Y,$ respectively.

\medskip
\par\noindent
We will follow the construction laid out in \cite{N1} to explicitly write out the Dirac operator and the Seiberg-Witten equations on $Y.$ 
Start from a constant-curvature conformal metric on $\Si$ 
and use the splitting (\ref{Splitting_e}) and the connection 1-form $\al$ on $Y$ 
as the first element of an orthonormal frame 
to equip $Y$ with a riemannian metric $\gmetric.$ Now fix a $\gmetric$-compatible connection on the tangent bundle of $Y$ and couple that with any unitary connection $B=B_{\tn{can}}\otimes B_E^2$ on the characteristic line bundle $\piysi^*L$ on $Y$ to define the Dirac operator $\dirac_B$ on $Y.$ According to the decomposition of $\sb$ in (\ref{U3_e}) above, any spinor on $Y$ splits as a section $\Psi=(\Psi_+,\Psi_-)$ of $\piysi^*\ub^+\oplus\piysi^*\ub^-\lra Y.$
With respect to this splitting, the Dirac operator $\dirac_B\colon\Ga(\sb)\lra\Ga(\sb)$ on $Y$ takes the form
\bEq{DiracOnS1bundle}
\dirac_B=
\left[ 
  {\begin{array}{cc}
   \nabla_{-\mfi\partial_\theta} & \sqrt{2}\dbar_{B_E}^* \medskip \\
   \sqrt{2}\dbar_{B_E} & \nabla_{\mfi\partial_\theta}  \\
  \end{array}} 
\right]+\la\bb1,
\eEq
where $\la$ is some constant depending on the choice of the $\gmetric$-compatible connection on $Y,$ 
$\nabla_{\partial_\theta}$ is the covariant derivative in the fiber direction, and $\dbar_{B_E}$ is defined with respect to the lift of the complex structure on $\Si$ to the horizontal tangent space, so that it is zero in the fiber direction \cite[(2.7)]{N1}. Moreover, if we choose our $\gmetric$-compatible connection on $Y$ to be the {\it adiabatic} connection, that is, one which acts trivially in the fiber direction $T^{\tn{ver}}Y,$ then $\la=0.$ Using this (non-Levi-Civita) connection on $Y,$ there is an explicit description of the solutions of the (perturbed) $\SW_\eta^3$ equations, which we state below. This is essentially based on 
\cite[Sec.~5, pp.~719--726]{mst}, \cite[Sec.~3.2, pp.~92--97]{N1}, 
\cite[Sec.~3, pp.~356--372]{N2}, and \cite[Sec.~5, pp.~701--722]{moy} 
and continues to hold for perturbation terms of the form
\bEq{om-pert_e}
\eta=\mfi (*n),
\eEq
where $n=\piysi^*(\nu+\ov{\nu})$ is the pullback of the real 1-form associated to a holomorphic 1-form $\nu.$

\bLm{Nic-lm}
With notation as above, and using the adiabatic connection ($\la\!=\!0$), for every solution $(B,\Psi)$ of $\SW^3_\eta$ 
we have
$$
\Psi_+\ov{\Psi}_-=\mfi\pi^*\nu,
$$
both $\Psi_+$ or $\Psi_-$ are covariantly constant in the vertical direction, $\dbar_{B_E}^*(\Psi_-)=0$ and $\dbar_{B_E}(\Psi_+)=0.$ The curvature equation of $\SW^3_\eta$ simplifies to
\bEq{CEPullBack_e}
F_B=\frac{|\Psi_+|^2-|\Psi_-|^2}{2}\,\mfi\piysi^*\om,
\eEq
where $\om$ is a K\"ahler form on $\Si$ with constant curvature. 
\eLm
\bPf
In the following, to keep the notation simple, we write $\dbar_B$ instead of $\dbar_{B_E}$.
By (\ref{DiracOnS1bundle}), for $\la\!=\!0$, the Dirac equation in $\tn{SW}^3_\eta$ becomes
\bEq{SW3Dirac_e}
\aligned
-\mfi\nabla_{\partial_\theta}\Psi_++\sqrt{2}\dbar_{B}^*\Psi_- & =0 \\
\sqrt{2}\dbar_{B}\Psi_++\mfi\nabla_{\partial_\theta}\Psi_-& =0
\endaligned
\eEq
Following the local calculations in \cite[pp.~722--723]{mst} and \cite[Appendix D]{N1}, fix a local orthonormal coframe $\al_1$ and $\al_2$ on $\Si$ such that $\nd\al_1=\kappa \al_1\wedge \al_2$ and $\nd \al_2=0$, where $\kappa$ is a non-negative constant. Note that $\om=\al_1\wedge \al_2$. By pulling back to $Y$, we can think of $\al_1$ and $\al_2$ as $1$-forms on $Y$. Let 
\bEq{alhol_dfn}
\al_{1,0}=\al_1+\mfi \al_2\qquad\tn{and}\qquad \al_{0,1}=\al_1-\mfi \al_2
\eEq
denote the corresponding $(1,0)$ and $(0,1)$ forms, respectively.
The triple $(\al_1,\al_2,\al)$ is an orthonormal coframe for $Y$ with dual frame $(\ze_1,\ze_2,\partial_\theta)$. In this notation we have 
\bEq{dbar_dfn}
\dbar_{B}=\al_{0,1} \otimes \frac{1}{2}\big(\nabla_{\ze_1}+\mfi\nabla_{\ze_2} \big).
\eEq
The curvature equation of $\tn{SW}^3_\eta$ reads:
\bEq{CurvEq_e}
\aligned
-\mfi(F_{\al_1,\al_2}-\eta_{\al_1,\al_2}) & = \frac{1}{2}(|\Psi_+|^2-|\Psi_-|^2) \\ 
\frac{1}{2}(F_{\al_2,\al}-\eta_{\al_2,\al}-\mfi F_{\al_1,\al}+\mfi\eta_{\al_1,\al})\al_{0,1} & =  \overline{\Psi}_+\Psi_-
\endaligned
\eEq
where 
$$
\aligned
F_B&= F_{\al_1,\al_2}~\al_1\wedge \al_2 + F_{\al_1,\al}~\al_1 \wedge \alpha+ F_{\al_2,\al}~\al_2 \wedge \alpha,\\
\eta&= \eta_{\al_1,\al_2}~\al_1\wedge \al_2 + \eta_{\al_1,\al}~\al_1 \wedge \alpha+ \eta_{\al_2,\al}~\al_2 \wedge \alpha.
\endaligned
$$

\medskip
\par\noindent
Applying $\dbar_{B}$ to the first equation in (\ref{SW3Dirac_e}) we get
\bEq{Start_Le}
-\mfi\dbar_{B}\nabla_{\partial_\theta}\Psi_++\sqrt{2}\dbar_{B}\dbar_{B}^*\Psi_-  =0.
\eEq
For $\eta$ as in (\ref{om-pert_e}), by the curvature equation   
$$
\dbar_{B} ~ \nabla_{\partial_\theta}-\nabla_{\partial_\theta}~ \dbar_{B}= (F_{\al_1,\al}+\mfi F_{\al_2,\al} )\al_{0,1},
$$
and using the second equation in (\ref{CurvEq_e}), we can rewrite (\ref{Start_Le}) as  
$$
-\mfi \nabla_{\partial_\theta}\,\dbar_{B}\Psi_+ + 2 (\ov{\Psi}_+\Psi_-+ \mfi\pi^*\ov\nu) \Psi_++ \sqrt{2}\dbar_{B}\dbar_{B}^*\Psi_-  =0.
$$
Applying the second equation in (\ref{SW3Dirac_e}) to the first term we get
$$
-\frac{1}{\sqrt{2}}\nabla_{\partial_\theta}\,\nabla_{\partial_\theta}\Psi_- + 2(\ov{\Psi}_+\Psi_-+ \mfi\pi^*\ov\nu) \Psi_+ + \sqrt{2} \dbar_{B}\dbar_{B}^*(\Psi_-)  =0.
$$
Since rotation in the circle direction acts by isometries on $TY$, we have
$\nabla_{\partial_\theta}^*=-\nabla_{\partial_\theta}$.
Taking the inner product with $\Psi_-$ then gives
\bEq{Norms_e}
\frac{1}{\sqrt{2}}\|\nabla_{\partial_\theta}\Psi_-\|^2 + 2 \|\ov{\Psi}_+\Psi_-+\mfi\pi^*\ov\nu \|^2+ \sqrt{2} \|\dbar_{B}^*(\Psi_-)\|^2 =0.
\eEq
Here, to get the middle term, we have used the fact that 
\bEq{VanishingPart_e}
\ll \mfi\pi^*\ov\nu, \ov{\Psi}_+\Psi_-+\mfi\pi^*\ov \nu\rr
=\int_{Y} \mfi\pi^*\ov\nu\wedge F_B=0,
\eEq
because $\ov\nu$ is a closed $1$-form on $\Si$ and $L$ is pulled back from $\Si$; see \cite[p.~726]{mst}.
Lemma~\ref{Nic-lm} now follows from (\ref{Norms_e}) and the first equation in (\ref{CurvEq_e}).
\ePf

\medskip
\par\noindent
Lemma~\ref{Nic-lm} implies that an irreducible solution $(B,\Psi)$ of $\SW^3_\eta$ is indeed pulled back from a pair of counterparts $(A,\Psi)$ on $\Si$ satisfying $\SW^2_\nu,$ with some abuse of notation on the spinor component. Such a solution $(A,\Psi)$ on $\Si$ defines a holomorphic structure on $\ub^+\!=\!E$ via $\dbar_{A_E}$, for which $\Psi_+$ is a holomorphic section. Analogous statements hold for the Serre dual bundle $\kcan\otimes E^*$ and $\ov{\Psi}_-.$ 
Therefore, to understand irreducible $\eta$-monopoles on $Y,$ we need to understand irreducible solutions of $\tn{SW}^2_\nu$ on $\Si.$

\medskip
\par\noindent
A word of caution here is in order though. In the discussion below, for the sake of specificity, we are going to fix the \spinc structure $\mfssie$ on $\Si.$ As a result, the line bundle $E=E_{\sss \Si}$ and its degree $d=\deg(E),$ as well as the characteristic line bundle $L=L_{\sss \Si}=\kcan^*\otimes E^2$ of $\deg(L)=2-2g+2d,$ are also fixed. However, if $\deg(Y)\neq 0$, recall that the \spinc structure $\mfsyeclass$ depends on the equivalence class of $E;$ therefore, to find all the solutions of $\tn{SW}^3_\eta$ on $Y$ corresponding to the \spinc structure $\mfsyeclass,$ we need to account for all the solutions to the $\tn{SW}^2_\nu$ equations on $\Si$ corresponding to each line bundle in the equivalence class $[E].$ We denote the irreducible component of $\cM(Y,\mfsyeclass)$ corresponding to the line bundle $E$ by $\cM(Y,\mfsye).$

\medskip
\par\noindent
Another point to note in the following discussion is that there is a natural involution in Seiberg-Witten theory which sends the characteristic line bundle $L$ to its inverse. Therefore, it will be enough to focus on the case where $\deg(L)\leq0,$ that is, when $d\leq g-1.$

\medskip
\par\noindent
\textbf{Unperturbed case.} This is the case of $\SW^3$ equations $(\eta\!=\!0).$
Let us first consider the case of $\deg(L)<0.$ If the equations have an irreducible solution, then it will be the pullback of a solution $(A,\Psi)$ of $\SW^2$ on $\Si.$ If $\Psi_+=0,$ then integrating the curvature equation of $\tn{SW}^2$
\bEq{SW2Curvature_e}
F_A=\frac{|\Psi_+|^2-|\Psi_-|^2}{2}\,\mfi\om
\eEq
over $\Si$ yields $\deg(L)=2(1-g+d)\geq0.$ Thus, when $\deg(L)< 0,$ we conclude that $\Psi_+$ must be non-zero. 
If the spinor $\Psi_+$ is non-zero, it will be unique up to a constant multiple if its zero set $\tn{Div}(\Psi_+),$ which is an effective divisor on $\Si,$ is prescribed. Moreover, the norm of this constant multiple can be determined from the curvature equation above by integrating over $\Si.$ Therefore, by fixing $\tn{Div}(\Psi_+),$ the plus spinor $\Psi_+,$ and hence the pair $(A,\Psi)$ on $\Si,$ will be unique up to gauge equivalence.
We conclude that the solution $[(B,\Psi)]$ on $Y$ will be determined by the zero set of $\Psi_+$ on $\Si,$ which is a point in $\tn{Div}_d(\Si),$ the space of effective divisors of degree $d$ on $\Si.$ 
In other words, each monopole in $\cM(Y,\mfsye)$ can be identified with a point in the symmetric $d$-fold product $\tn{Sym}^d(\Si),$ where its plus spinor ``lands".
The process above defines a smooth {\it landing map}
\bEq{Landing_e}
\flat:\cM(Y,\mfsye)\lra\tn{Div}_{d}(\Si)\cong\tn{Sym}^{d}(\Si),
\eEq 
which turns out to be a diffeomorphism. The inverse $\flat^{-1}$ is constructed in the following way. Let $D$ be an arbitrary effective divisor in $\tn{Div}_d(\Si).$ Corresponding to this divisor, we will have a unique holomorphic line bundle $[D]$ on $\Si$ of degree $d,$ which is topologically the same as $E,$ together with a holomorphic section vanishing at $D,$ which is unique up to scalar multiplication. Using the induced holomorphic connection on $E$ and the holomorphic section, we obtain a solution $(A,\Psi)$ to $\SW^2,$ unique up to gauge equivalence. We can then pull $(A,\Psi)$ back to a solution for $\SW^3$ on $Y.$

\medskip
\par\noindent
Similarly, in the case of $\deg(L)>0,$ we can conclude in the same vein that $\Psi_-$ must be non-zero. In this case, instead of $\Psi_+,$ we will look at the holomorphic section $\ov{\Psi}_-$ of the Serre dual line bundle $\kcan\otimes E^*$ and 
$$d=\deg(\ub^+)=\deg(E)=(g-1)+\frac12\deg(L)$$ 
will be replaced by 
$$\deg((\ub^-)^*)=\deg(\kcan\otimes E^*)=2g-2-d=(g-1)-\frac12\deg(L)$$ 
as the degree of the effective divisor $\tn{Div}(\ov\Psi_-).$
This justifies our definition of $d(\sfrak)$ as in (\ref{DegS_e}) in the introduction.
The landing map $\flat\!\colon\!\cM(Y,\mfsye)\!\lra\!\tn{Sym}^{d(\sfrak)}(\Si)$ will similarly be defined as in (\ref{Landing_e}).
Meanwhile, note that the existence of a non-zero holomorphic section for $\ub^+$ (respectively $(\ub^-)^*$) implies that $d\geq 0$ (respectively $d\leq 2g-2$). Therefore, we have just showed that there are no irreducible solutions to the $\SW^3$ equations unless $0\leq d=\deg(E)\leq 2g-2.$
Moreover, when $d=0$ (respectively $d=2g-2$), the line bundle $\ub^+$ (respectively $\ub^-$) is trivial, so such a non-zero section will be constant. Therefore, if $d=0,2g-2,$ the irreducible solution to $\SW^3$ will be unique up to gauge equivalence, unless $g=1,$ in which case $\deg(L)=0$ and there will be no irreducible solution at all, as will be discussed below.
For consistency of notation throughout this paper, we take $\tn{Div}_{d}(\Si)\cong\tn{Sym}^{d}(\Si)$ to consist of a single point when $d=0.$
\medskip
\par\noindent
Finally, we consider the case of $\deg(L)=0,$ that is, when $L=\bbC_{\sss \Si},$ $E=\sqrt{\kcan}$ and $d=g-1.$
If $\SW^3$ has an irreducible solution, by integrating the curvature equation (\ref{SW2Curvature_e}), we conclude that the spinors $\Psi_+$ and $\Psi_-$ have the same $L^2$-norms on $\Si.$ Since the product of these spinors is zero, they both have to be identically zero. 
This means that all the solutions to $\SW^3,$ corresponding to $E=\sqrt{\kcan},$ are reducible and correspond to flat connections on the torsion line bundle $\piysi^*L$ on $Y.$ 
The space of such flat connections, up to gauge equivalence, is homeomorphic to the Jacobian torus $H^1(\Si,S^1)=H^1(\Si,\R)/H^1(\Si,\Z)=\T^{2g}$ when $\deg(Y)\neq 0.$
Moreover, if $\deg(L)$ is not a multiple of $\deg(Y),$ the reducible component will be non-degenerate in the sense defined in \cite[Def.~5.15, p.~717]{moy}.
If $\deg(L)$ is a multiple of $\deg(Y),$ this space can be identified with the theta divisor $W_{g-1}$ inside the Jacobian $J_{g-1}(\Si)$ (see Fact~3 in \cite[p.~95]{N1}, \cite{GH}). 
If $\deg(Y)=0,$ that is, when $Y=S^1\times\Si$ is a trivial circle bundle over $\Si,$ $\piysi^*L$ will be topologically trivial and its space of flat connections can be identified with the torus $\T^{2g+1}$ as the space of $S^1$-representations of the fundamental group of $Y$ (see \cite[p.~369]{N2}).

\medskip
\par\noindent
Now, let $\cJ$ denote the component of reducible solutions to $\SW^3,$ that is, $\cJ$ is the torus $\T^{2g}$ when $\deg(Y)\!\neq\!0,$ the torus $\T^{2g+1}$ when $\deg(Y)\!=\!0$ and $E\!=\!\sqrt{\kcan},$ and empty otherwise. 
The preceding discussion shows that we can write the moduli space of unperturbed $\SW^3$ equations as
$$\cM(Y,\mfsyeclass)=\mc{J}\cup\bigcup_{\sss E\in[E]}\cM(Y,\mfsye),$$
where $\cM(Y,\mfsye)\cong\tn{Sym}^{d(\mfsye)}(\Si)$ is the irreducible component\footnote{We are using a different notation here compared to the notation in the introduction. There, we denoted $\cM(Y,\mfsye)$ by $\cM(Y,\mfsy)_m,$ where $m=\deg(L)/2.$} 
for which $(\Psi_+,\ov{\Psi}_-)$ are holomorphic sections of $(E,\kcan\otimes E^*)$ and $\cJ$ is the reducible component as above.
\medskip
\par\noindent
When $\deg(Y)\!=\!0,$ the pullback $\spinc$ structure on $Y$ is uniquely identified by $E$ (i.e., $[E]=\{E\}$) and the relation above can be simplified. 
If $E\!\neq\!\sqrt{\kcan},$ we simply get an irreducible moduli space
$$\cM(Y,\mfsye)\cong\tn{Sym}^{d(\mfsye)}(\Si)$$
and if $E\!=\!\sqrt{\kcan},$ the moduli space is just the reducible component $\T^{2g+1}.$ 
Of course, a perturbation of the $\SW^3$ equations will help us get rid of these reducible solutions, as explained below.  

\medskip
\par\noindent
\textbf{Perturbed case.} This is the case of $\SW^3_\eta$ equations with $\eta\!=\!\mfi(*n)\!\neq\!0.$
As before, we conclude from Lemma~\ref{Nic-lm} 
that both spinors $\Psi_+$ and $\Psi_-$ are covariantly constant in the direction of the circle, $(B,\Psi)$ is a pullback of $(A,\Psi)$ from $\Si$, and $\dbar_{A_E}^*(\Psi_-)\!=\!\dbar_{A_E}(\Psi_+)\!=\!0.$ This time, however, none of $\Psi_\pm$ can be identically zero, so our solution is not reducible. Moreover, since $\Psi_+\ov{\Psi}_-=\mfi\nu,$ knowing one of the spinors $\Psi_+$ and $\Psi_-$ will uniquely identify the other.  Again, consider the zero set $\tn{Div}(\Psi_+)$ of the holomorphic spinor $\Psi_+$ on $\Si.$ This divisor will determine $\Psi_+$ up to a constant multiple.  The norm of this constant multiple can be determined from (\ref{SW2Curvature_e}) 
by integrating over $\Si,$ as in \cite[Cor.~5.5, p.~726]{mst}. 
Therefore, by fixing $\tn{Div}(\Psi_+),$ the plus spinor $\Psi_+,$ and hence the pair $(A,\Psi)$ on $\Si,$ will be uniquely determined up to gauge equivalence. As in the unperturbed case above, for degree reasons, the equations have no solution unless $0\leq \deg(E) \leq 2g-2$. We also have
$$
\tn{Div}(\Psi_+)+\tn{Div}(\ov{\Psi}_-)= \tn{Div}(\nu),
$$
so, unlike the unperturbed case above, $\tn{Div}(\Psi_+)$ is constrained to a sum of $d$ points in the zero set $\tn{Div}(\nu).$ Then $\tn{Div}(\ov\Psi_-)$ is the sum of the remaining $2g-2-d$ points in the complement of $\tn{Div}(\Psi_+)$ in $\tn{Div}(\nu)$. Therefore, $\cM_\eta(Y,\mfs_{\sss Y,E})$ is a finite set of regular points and the reducible component $\cJ$ is empty.
\subsection{SW equations on a four-manifold with a cylindrical end}\label{SW4cyend_ss}
%
Let $M$ be a smooth connected oriented riemannian 4-manifold with a cylindrical end\footnote{One may allow more than one cylindrical end. The results will readily generalize.}, that is, one which is orientation-preserving isometric to the cylinder $\R_+\!\times\! Y$ outside a compact submanifold $M_o$ with non-empty interior, where $Y$ is a closed oriented riemannian 3-manifold.
Moreover, the boundary $\partial M_o$ has a collar neighborhood diffeomorphic to $(-1,0]\times Y,$ where $\partial M_o$ corresponds to $\{0\}\times Y.$
Use the above isometry and diffeomorphism to identify the cylindrical end (also called the neck) and the collar with $\R_+\!\times\! Y$ and $(-1,0]\times Y,$ respectively, and let $\tcyend$ denote the longitudinal ``time coordinate" on the neck $\R_+\!\times\! Y.$ 
We can view $\tcyend$ as the projection map onto the first factor $\tcyend\colon\R_+\!\times\! Y\lra\R_+,$ 
so that $\tcyend(t,y)\!=\!t$ for $(t,y)\!\in\!\R_+\!\times\! Y,$
and smoothly extend it to the entire manifold $\tcyend\colon M\lra[-1,\infty)$ using the collar diffeomorphism 
such that $\tcyend^{-1}(\R_+)$ is the neck, 
$\tcyend^{-1}((-1,0])$ is the collar, 
and $\tcyend^{-1}(-1)$ is the rest of the manifold. 
We will have $M_o=\tcyend^{-1}([-1,0])$ and $\partial M_o=\tcyend^{-1}(0).$

\medskip
\par\noindent
We are going to consider the Seiberg-Witten equations on this cylindrical-end manifold $M,$ where we perturb the equations using an {\it adapted} perturbation term 
\bEq{AdaptedEta_e}
\etacyend=\etacyend(\eta_o,\eta)=\eta_o+\be(\tcyend)\eta^+,
\eEq 
in which 
\bIt
\item $\eta_o$ is a purely imaginary self-dual 2-form on $M$ with compact support, 
\item $\eta^+$ denotes the self-dual part of the pullback of a closed 2-form $\eta\in\Om^2(Y,\mfi\R)$ under the natural projection $\R_+\!\times\! Y\lra Y,$ 
\item $\be(\tcyend)$ is a smooth cutoff function on $M$ which is identically 1 on the neck, non-negative on the collar, and zero otherwise. It is defined as the composition of a smooth increasing function $\be\colon[-1,\infty)\lra[0,1]$ satisfying $\be(-1)=0$ and $\be(0)=1$ with the map $\tcyend\colon M\lra[-1,\infty).$
\eIt
An adapted perturbation will enable us to control how the Seiberg-Witten equations are perturbed on the compact piece and the neck, and to study the equations on the cylinder $\R_+\!\times\! Y$ (with the product metric) as a stand-alone problem. 
Fixing a \spinc structure on the cylinder will induce a \spinc structure $\syfrak$ on $Y$ by restriction,
where the spinor bundle $\sb_{\sss Y}$ is identified with both $\sb^+$ and $\sb^-,$ which are themselves identified via Clifford multiplication by $\nd t.$
Conversely, a \spinc structure on $Y$ induces one on the cylinder in an obvious way via pullback; we will thus use the notation $\syfrak$ interchangeably for both.
A connection on the cylinder is called {\it temporal} if it acts as $\partial_t$ in the time direction; this can always be arranged by a suitable gauge transformation.

\medskip
\par\noindent
Therefore, any solution $(A,\Phi)$ of $\SW^4_{\eta^+}$ on the cylinder $I\!\times\! Y,$ where $I\!\subset\!\R$ is an interval and $A$ is in temporal gauge, can be identified with a path $(B(t),\Psi(t))$ in the configuration space on $Y.$
Moreover, these paths are the (upward) gradient flow lines for the {\it Chern-Simons-Dirac} function on this configuration space, as defined below. Fix a background connection $B_0$ on $L_{\sss Y}$ and, for $$(B,\Psi)\in\cC(Y,\mfsy)=\mc{A}(L_{\sss Y})\times\Gamma(\sb_{\sss Y}),$$ 
define
$$
\tn{CSD}_\eta(B,\Psi)=\int_Y F_{B_0}\wedge b -\int_Y \eta\wedge b 
+\frac{1}{2}\int_Y b\wedge\nd b +\int_Y \langle\Psi,\dirac_B\Psi\rangle\nd\tn{vol},
$$
where $b=B-B_0$ is a 1-form on $Y.$ A different choice of $B_0$ would change $\tn{CSD}_\eta$ only by a constant and therefore the gradient on $\cC(Y,\syfrak)$ is independent of this choice. 
Now, the gradient flow lines of $\tn{CSD}_\eta$ identify with solutions to $\tn{SW}^4_{\eta^+}$ with temporal gauge and the stationary points correspond to solutions to $\tn{SW}^3_\eta.$ 
Note that changing the orientation on $Y$ reverses the direction of the flow lines
and that $\tn{CSD}_\eta$ is locally constant on $\cM_\eta(Y,\syfrak)$ as the zero locus of the gradient, hence constant on each connected component.

\medskip
\par\noindent
The function $\tn{CSD}_\eta$ is not necessarily invariant under the gauge transformation $u\colon Y\!\lra\! S^1$ though and will change by
$$\tn{CSD}_\eta(u\cdot(B,\Psi))-\tn{CSD}_\eta(B,\Psi)=4\uppi^2\Big([u]\cup\big(c_1(\mfsy)\!-\!\frac{\mfi}{2\uppi}[\eta]\big)\Big)[Y],$$
where $[u]$ denotes the homotopy class of $u\colon Y\!\lra\! S^1$ as a cohomology class in $H^1(Y,\Z).$
Thus, $\tn{CSD}_\eta$ descends to a real-valued function on the quotient if $c_1(\mfsy)\!=\!\frac{\mfi}{2\uppi}[\eta],$ 
and to a circle-valued function otherwise.
In the former case, the value of $\tn{CSD}_\eta$ is increasing along the flow lines in the quotient.

\medskip
\par\noindent
Now let $(A,\Phi)$ be a solution to $\tn{SW}^4_{\eta^+}$ on $I\!\times\! Y,$ where $A$ is in temporal gauge, and identify that with $(B(t),\Psi(t))$ in the configuration space on $Y,$ as above. The topological {\it energy} of $(A,\Phi)$ is defined by any of the following equivalent formulas:
$$
\label{energy_eq}
E(A,\Phi)=\int_I \|\dot{B}\|^2 +\|\dot{\Psi}\|^2
=\int_I \|\nabla \tn{CSD}_\eta(B(t),\Psi(t))\|^2
=\tn{CSD}_\eta(B(t_1),\Psi(t_1))-\tn{CSD}_\eta(B(t_0),\Psi(t_0)),
$$
where in the last equation the interval $I$ is assumed to be $I=[t_0,t_1].$
When working with four-manifolds $M$ with cylindrical ends, we will only consider solutions with finite energy on the ends and $\cM_{\etacyend}(M,\smfrak)$ denotes the moduli space of finite-energy monopoles on $M$ modulo gauge transformations. We will complete the relevant underlying spaces using appropriate {\it weighted} Sobolev norms, as in \cite{gluingbook} or \cite{N}.

\medskip
\par\noindent
If $[A,\Phi]=[B(t),\Psi(t)]$ is a finite-energy $\eta^+$-monopole on $\R_+\!\times\! Y,$ one can show that the limit
$[B,\Psi]=\lim_{t\to\infty}[B(t),\Psi(t)]$
exists and is an $\eta$-monopole on $Y,$ using the same methods established in \cite{mmr} and \cite{tl2}. 
Moreover, if the limit $[B,\Psi]$ is irreducible, the convergence to the asymptotic limit is exponential. 
This defines a {\it limiting map} 
\bEq{bMap_e}
\partial\colon\cM_{\etacyend}(M,\smfrak)\lra\cM_\eta(Y,\syfrak).
\eEq

\bRm{Generic-eta_rmk}
If $b^+(M)>0,$ a generic choice of the compact perturbation $\eta_o$ will ensure that all $\etacyend$-monopoles on $M$ are irreducible and {\it strongly regular}, in the sense defined in \cite[Def.~4.3.20, p.~383]{N}. This implies that, for a generic compact perturbation $\eta_o,$ the moduli space $\cM_{\etacyend}(M,\smfrak)$ is a smooth manifold \cite[Prop.~4.4.1, p.~405]{N} and the limiting map above is a submersion \cite[Rmk.~4.4.2]{N}.
\eRm

\medskip
\par\noindent
A finite-energy monopole on $\R\times Y$ is called a {\it tunneling}. 
To simplify notation in the following discussion, let us drop the reference to \spinc structures and assume $\eta=0.$
The real line $\R$ acts on $\cM(\R\times Y)$ by translations and we denote the quotient, that is, the moduli space of unparametrized tunnelings, by $\breve\cM(\R\times Y).$ Corresponding to the two ends of the cylinder $\R\times Y,$ there are two induced limiting maps $\partial_\pm\colon\breve\cM(\R\times Y)\lra\cM(Y)$ as $t\to\pm\infty,$ since the limiting maps are $\R$-invariant. If $C_-$ and $C_+$ are two components of $\cM(Y),$ we will denote the moduli space of such tunnelings from $C_-$ to $C_+$ by $\breve\cM(C_-,C_+)\subset\breve\cM(\R\times Y).$
If the \spinc structure $\syfrak$ is torsion, so that $\tn{CSD}$ is a real-valued function on the moduli space, then $\breve\cM(C_-,C_+)$ will be empty if $\tn{CSD}(C_-)>\tn{CSD}(C_+),$ but there may be non-trivial tunnelings when $\tn{CSD}(C_-)<\tn{CSD}(C_+).$
We will use these tunnelings to construct a suitable compactification for $\cM_{\eta_o}(M,\smfrak)$ 
in Section~\ref{Compact_ss}.

\medskip
\par\noindent
The moduli space $\cM_{\etacyend}(M,\smfrak)$ in general is not compact and has components of different dimensions. At any point of the moduli space, the expected real dimension can be calculated using the Atiyah-Patodi-Singer index theorem \cite[Thm.~3.10]{APS} and is given by 
\bEq{DimCyl4_eq}
\dsw=\frac{c_1(A)^2-2\chi(M)-3\si(M)}{4}+\xi,
\eEq
where $\xi$ stands for $\frac{1}{2}(\eta(0)+h),$ i.e., half of the sum of the APS $\eta$-invariant and the dimension of the kernel of $\SW^3_\eta$ operator on $Y$ at the given point.\footnote{We have opted to work with the notation convention $\xi$ of \cite{N1} instead to avoid confusion regarding the terms $h$ and $\eta$ we have used elsewhere in the paper.}
The term $c_1(A)^2$ in (\ref{DimCyl4_eq}) denotes the integral of the Pontrjagin form as in \cite[Thm.~4.2]{APS} and is given by 
\bEq{Integralc_1^2_e}
c_1(A)^2\defeq\frac{-1}{4\uppi^2}\int_M F_A\wedge F_A.
\eEq
The signature $\si(M)$ is defined as in \cite[Prop.~7.1]{AS} and equals $b^+_M-b^-_M,$ the difference of dimensions of maximal positive- and negative-definite subspaces in $H^2_c(M,\R)$ with respect to the Poincar\'e intersection pairing.
The dimension (\ref{DimCyl4_eq}) is not a topological quantity and changes as we move from one component of $\cM_{\etacyend}(M,\smfrak)$ to another. We will fix this problem on $M=X-\Si$ by choosing a \spinc structure for $TX(-\log\Si)$ over the entire $X.$ 

\subsection{Gluing monopoles}\label{Gluing_ss}

In this section, which is based on \cite{gluingthesis}, \cite{swgluing}, \cite{gluingbook}, we state a gluing theorem for SW monopoles on a closed 4-manifold $X$ if it can be ``decomposed" as a union of two cylindrical-end manifolds $X_\pm$, in the sense described below. 
Proof of the sum formula in Theorem~\ref{SW-Sum_thm2} follows from this gluing theorem and a convergence result. Such a gluing statement is also needed to prove that the compactified moduli space $\ov\cM_{\eta_o}(X-\Si,\mfs)$ is a $C^0$-manifold.

\medskip
\par\noindent
Let $Y$ be a smooth closed connected oriented riemannian 3-manifold.
Suppose $X_+$ and $X_-$ are two smooth connected oriented riemannian 4-manifolds with cylindrical ends $\R_+\!\times\!Y$ and $\R_-\!\times\!Y,$ respectively. 
Note that we can also think of $X_-$ as a cylindrical-end manifold with cylindrical end $\R_+\!\times\!Y$ at the expense of changing the orientation on $Y$. 
The collars of $X_\pm$, as in Section~\ref{SW4cyend_ss}, are identified with $(-1,0]\!\times\!Y$ and $[0,1)\!\times\!Y$, respectively.
For each $T\gg0$, let $X_{\sss T}$ denote the closed riemannian 4-manifold obtained 
by identifying
$$
X_{\sss T,+}=X_+-(T,\infty)\!\times\!Y \qquad \tn{and} \qquad 
X_{\sss T,-}=X_--(-\infty,-T)\!\times\!Y 
$$
along $$[0,T]\times Y \qquad \tn{and} \qquad [-T,0]\times Y,$$ respectively, via the shift map $t_+\longleftrightarrow t_-=t_+-T,$ i.e., 
\bEq{Neck_e}
[0,T]\times Y \ni (t_+,y)\sim (t_+-T,y)\in [-T,0]\times Y.
\eEq
When the choice of $T$ is not important, we write $X$ to denote any $X_{\sss T}$. 

\medskip
\par\noindent
With notation as in the previous section, suppose $\mfs_\pm$ are \spinc structures on $X_\pm$ that induce the same \spinc structure on $Y$ (or $\R_\pm\times Y$). By gluing $\mfs_\pm$ along the neck (\ref{Neck_e}) 
we obtain a \spinc structure $\mfs$ on $X$.  

\medskip
\par\noindent
Let $\beta_{\sss T}\colon X_{\sss T} \lra \R$ be a smooth non-negative function which identically equals $1$ on the neck and vanishes on
$$
\Big(X_+-(-1,\infty)\!\times\!Y\Big) \cup  
\Big(X_--(-\infty,1)\!\times\!Y\Big) \subset X_{\sss T}.
$$
We consider an {\it adapted} perturbation form on $X_{\sss T}$
\bEq{AdaptedEta_e2}
{\etacyend}_{\sss T}=\eta_{o,-}+\be_{\sss T}\cdot\eta^++\eta_{o,+},
\eEq 
in which 
\bIt
\item $\eta_{o,+}$ is a purely imaginary self-dual 2-form on $X_+-(-1,\infty)\!\times\!Y$ with compact support, 
\item $\eta_{o,-}$ is a purely imaginary self-dual 2-form on $X_--(-\infty,1)\!\times\!Y$ with compact support, 
\item $\eta^+$ denotes the self-dual part of the pullback of a closed 2-form $\eta\in\Om^2(Y,\mfi\R)$ under the natural projection $\R\!\times\! Y\lra Y.$
We can make sense of $\be_{\sss T}\cdot\eta^+$ as a 2-form on $X_{\sss T}$ in an obvious way. 
\eIt
We can meanwhile use these forms to define adapted perturbations $\etacyend_\pm$ on $X_\pm$ as in (\ref{AdaptedEta_e}). 
When restricted to $X_{\sss T,+}$ and $X_{\sss T,-}$ as subsets of $X_{\sss T},$ the perturbations $\etacyend_\pm$ coincide with $\etacyend_{\sss T},$ respectively.
The following is an immediate consequence of the gluing theorem in \cite{gluingthesis, gluingbook}.

\bTh{Gluing_thm}
Suppose $Y$ is a circle bundle over a Riemann surface, as in Section~\ref{SW3_ss}.
Let $\check{\cM}_\eta(Y,\mfsy)$ denote an irreducible component of $\cM_\eta(Y,\mfsy)$
and $\check{\cM}_{\etacyend_\pm}(X_\pm,\mfs_\pm)$ be the union of those components of $\cM_{\etacyend_\pm}(X_\pm,\mfs_\pm)$ that map into $\check{\cM}_\eta(Y,\mfsy)$ under the limiting map (\ref{bMap_e}).
Assume moreover that the perturbation terms are chosen generically so that $\check{\cM}_{\etacyend_\pm}(X_\pm,\mfs_\pm)$ are regular at every point and the limiting maps are submersions, therefore transversal, so the fiber product is a smooth manifold. Then, for sufficiently large $T$, there is a natural gluing map 
$$
\check{\cM}_{\etacyend_+}(X_+,\mfs_+)\, \times_{\check{\cM}_{\eta}(Y,\mfsy)} \check{\cM}_{\etacyend_-}(X_-,\mfs_-)\lra \cM_{{\etacyend}_{\sss T}}(X_{\sss T},\mfs),
$$
which is an embedding onto an open subset of $\cM_{{\etacyend}_{\sss T}}(X_{\sss T},\mfs).$
\eTh

\medskip
\par\noindent
As a special case, when the perturbation form $\eta$ on $Y$ is zero, we will have compact perturbations on the cylindrical-end pieces $X_\pm$ and the gluing map above is written as a map
\bEq{GCS_e}
\check{\cM}_{\eta_{o,+}}(X_+,\mfs_+)\, \times_{\check{\cM}(Y,\mfsy)} \check{\cM}_{\eta_{o,-}}(X_-,\mfs_-)\lra \cM_{{\eta}_{o,\sss T}}(X_{\sss T},\mfs),
\eEq
which is an embedding onto an open subset of $\cM_{{\eta}_{o, \sss T}}(X_{\sss T},\mfs)$, with $\eta_{o,\sss T}=\eta_{o,-}+\eta_{o,+}$.

\section{Relative Seiberg-Witten moduli spaces}\label{RelSW_s}
In Section~\ref{LogStr_ss}, we define the notion of logarithmic tangent bundles and \spinc structures. Then in Section~\ref{RelSWModuli_ss}, using the moduli spaces defined in~Section~\ref{setup_s}, we elaborate on the definition and properties of the relative moduli spaces $\cM_{\eta_o}(X-\Si,\mfs)$ and $\cM_{\etacyend}(X-\Si,\mfs)$. In Section~\ref{Dim_ss}, we use \cite{N1} to find topological formulae for the dimensions of the moduli spaces $\cM_{\eta_o}(X-\Si,\mfs)$ and $\cM_{\etacyend}(X-\Si,\mfs)$, as well as $\cM_{\eta_o}(X-\Si,\mfs_{\sss X-\Si},\mc{J})$. This involves a highly non-trivial calculation of the $\xi$-invariant in formula (\ref{DimCyl4_eq}), which is carried out in Section~\ref{PfPrp_ss}.
In Section~\ref{Tunneling_ss}, we give a holomorphic description of the tunneling moduli spaces. This plays a key role in proving that certain relative moduli spaces are indeed compact.
In Section~\ref{Compact_ss}, we describe the compactification of $\cM_{\eta_o}(X-\Si,\mfs)$.
In Section~\ref{Proofs_ss}, we show how the results from the earlier sections come together to prove Theorems~\ref{Relaive_thm}--\ref{SW-Sum_thm2}.

\subsection{Logarithmic tangent bundles and relative \spinc structures}\label{LogStr_ss}
Suppose $X$ is a closed oriented 4-manifold and $\Si$ is a closed oriented surface of genus $g$ in $X.$ Equip $\Si$ with a K\"ahler structure $(\mfj,\om)$  compatible with the orientation. The normal bundle 
$$
\pi\colon \cN=\frac{TX|_\Si}{T\Si}\lra \Si
$$
is an oriented vector bundle of real rank 2. Choose a complex structure $\mfi$ on $\cN$ compatible with the orientations and let $(\varrho,\nabla)$ be a hermitian metric and a compatible complex linear connection. By abuse of notation, we also let $\varrho(v)$ denote $|v|^2$ for every $v\!\in\!\cN$. The connection $\nabla$ gives rise to a decomposition of $T\cN$ into horizontal and vertical directions
\bEq{Decomp_e}
T\cN=T^{\tn{hor}}\cN\oplus T^{\tn{ver}}\cN \cong \pi^*T\Si\oplus\pi^*\cN.
\eEq
The tuple $(\mfi,\varrho,\nabla)$ also gives rise to a $1$-form $\alpha$ on $\cN-\Si$ such that $\alpha(\partial_\theta)\!=\!1$, $\alpha$ is zero on the horizontal sub-space, and 
$\nd\al=\pi^*F$ with $[F]=-2 \uppi c_1(\cN)$. The $2$-form
\bEq{KahlerN_e}
\om_\cN=\pi^*\om_\Si+\frac{1}{2}\nd(\varrho\alpha) =\pi^*\om_\Si+\frac{1}{2}\varrho~\pi^*F+\frac{1}{2}\nd\varrho\wedge\alpha
\eEq
on $\cN$ is closed and non-degenerate in a sufficiently small neighborhood $\cD$ of $\Si$ in $\cN$. Via decomposition (\ref{Decomp_e}), define $J_\cN$ to be the almost complex structure on $\cN$ given by $\pi^*\mfj$ on $\pi^*T\Si$ and $\pi^*\mfi$ on $\pi^*\cN$. The tuple 
$$
\big(\cD,\om_\cD\defeq \om_\cN|_{\cD} ,J_\cD\defeq J_{\cN}|_{\cD}\big)
$$ 
is K\"ahler. In other words, the $\dbar$ operator associated to $\nabla$,
\bEq{cNHolo_e}
\dbar_\na\ze= \nabla^{0,1}\ze= \frac{1}{2}\big (\nabla \ze + \mfi \nabla \ze \circ \mfj\big), \qquad \forall~\ze\in\Gamma(\Si,\cN),
\eEq
is integrable and defines a holomorphic structure on $\cN$. Locally around every point $p\!\in\!\Si$, if $z$ is a local holomorphic coordinate with $z(p)=0$ and $\ze$ is a holomorphic section with $\ze(p)\neq 0,$ then 
\bEq{LHC_e}
\cN\ni(z,w\ze(z))\lra (z,w)\in\C^2
\eEq
defines the holomorphic structure associated to the almost complex structure $J_\cN$.

\medskip
\par\noindent
Let 
$$
\Upsilon\colon \cD\lra X, \qquad \Upsilon|_\Si=\tn{id}, \qquad \nd\Upsilon|_{\Si}=\tn{id},
$$
be an identification of a neighborhood $\cD$ of $\Si$ in $X$ with a neighborhood of that in $\cN$. Via this identification, the previous paragraph shows that a neighborhood of $\Si$ in $X$ can be equipped with a K\"ahler structure of a standard form. If $(X,\om_X)$ is a symplectic manifold and $\Si$ is a symplectic submanifold with $\om\!=\!\om_X|_\Si$, by Symplectic Neighborhood Theorem \cite[Thm.~3.4.10]{MS17}, we can choose $\Upsilon$ to also satisfy
$$
\Upsilon^*\om_X=\om_\cD.
$$ 
Then, via the identification $\Upsilon$, we can extend the complex structure $J_\cD$ to a compatible almost complex structure on the entire $X$ compatible with $\om_X$.

\medskip
\par\noindent
If $X$ is a holomorphic surface and $\Si$ is a smooth holomorphic curve in $X$, the {\it log tangent sheaf} $\cT X(-\log \Si)$ is a sub-sheaf  of the holomorphic tangent sheaf $\cT X$ that is equal to $\cT X$ away from $\Si$ and is generated by
\bEq{LGT_e}
\partial^{\log}_{z}\defeq z\partial_{z} \quad \tn{and}\quad \partial_{w}
\eEq
in any local holomorphic chart $(z,w)\colon V\!\lra\!\C^2$ around a point of $\Si$ with $\Si\cap V \equiv (z=0)$.
Since $\cT X(-\log \Si)$ is locally free, it is the sheaf of holomorphic sections of a holomorphic vector bundle $TX(-\log \Si)$. The inclusion $\cT X(-\log \Si)\!\subset\!\cT X$ gives rise to a holomorphic homomorphism
\bEq{homo_e}
\iota\colon T X(-\log \Si)\lra T X
\eEq
which is an isomorphism away from $\Si$. 

\medskip
\par\noindent
Since $TX$ and $TX(-\log \Si)$ only differ in a neighborhood of $\Si$, the notion of logarithmtic tangent bundle can be defined for every pair $(X,\Si)$ of an oriented $4$-manifold and an oriented $2$-dimensional submanifold. 
For such $(X,\Si)$ and any choice of $(\Upsilon,\mfj,\mfi,\varrho,\nabla)$ as above, define
\bEq{LTdfn_e}
\aligned
&TX(-\log \Si)= \big(\Upsilon^{-1\,*} (\pi^* T\Si \oplus \C_\cD) \sqcup 
TX|_{X-\Si}\big)\big/\!\!\sim,\\
&\Upsilon^{-1\,*} (\pi^* T\Si\oplus \C_\cD)
  \ni\big(\Upsilon(v),u \oplus c\big)\sim\big(\Upsilon(v),\nd_v\Upsilon(u+cv)\big)\in TX|_{X-\Si},
  \endaligned
\eEq
where $\C_\cD$ is the trivial complex line bundle on $\cD$ and in the last equation we think of $u+cv$ as a tangent vector in $T_v\cN$ via the right-hand side of identification (\ref{Decomp_e}). 
Instead of the global definition in (\ref{LTdfn_e}), one may as well define $TX(-\log \Si)$ via local holomorphic charts (\ref{LHC_e}) in a neighborhood of $\Si$ as in (\ref{LGT_e}) and identify that with $TX|_{X-\Si}$ on the overlap as in the second line of (\ref{LTdfn_e}).

\medskip
\par\noindent
The riemannian metric $\om(\cdot,\mfj\cdot)$ on $T\Si$ and the standard metric on $\C_\cN$ give us a riemannian metric on $\pi^*T\Si\oplus\C_\cD.$ Via (\ref{LTdfn_e}), that can be extended to a metric $h_{\sss X}$ on the entire $TX(-\log\Si).$ The restriction of this metric to $X-\Si$ via the identification 
\bEq{Restriction_e}
TX(-\log \Si)|_{X-\Si}\cong T(X-\Si)
\eEq
defines a complete riemannian metric $h_{\sss X-\Si}$ on $X\!-\!\Si,$ giving it the structure of a riemannian 4-manifold with a cylindrical end $[0,\infty)\!\times\!Y,$ where $Y$ is the unit circle bundle in $\cN$ and $[0,\infty)$ is parametrized by $t\!=\!-\frac{1}{2}\log(\varrho).$ 
Note that in this paper, $Y$ gets the orientation induced by $X\!-\!\Si,$ that is, the one induced by its cylindrical end $[0,\infty)\!\times\!Y$ using the ``outward-normal-first" convention. 
This is opposite to the orientation $Y$ gets as $\partial\cD;$ therefore, in calculations related to the degree of the circle bundle $Y\!\lra\!\Si,$ we have
\bEq{degY_e}
\deg(Y)=-\Si\cdot\Si.
\eEq

\medskip
\par\noindent
The dual space
$$T^*X(\log \Si)\defeq TX(-\log \Si)^*$$ 
is the {\it logarithmic cotangent space}. 
Along $\Si,$ the space $\Omega^1_{\log}(X)$ of smooth sections of $T^*X(\log\Si)$ is locally generated by $\nd z/z$ and $\nd w$ and their complex conjugates.

\medskip
\par\noindent 
If $(X,\om_X)$ is a symplectic manifold, $\Upsilon$ is a symplectic identification, and $J$ is an almost complex structure on $X$ such that $\Upsilon^*J=J_\cD$, the complex structures $\pi^*\mfj\oplus\mfi_{\C}$ on $\pi^*T\Si\oplus\C_\cD$ and $J$ on $TX|_{X-\Si}$ match on the overlap region and define a complex structure on $TX(-\log \Si)$. Then, we have 
$$
c(TX(-\log \Si))=c(TX)/(1+\tn{PD}(\Si)).
$$
This identity follows from an isomorphism of vector bundles
$$
TX(-\log \Si) \oplus \cO(\Si)\cong TX\oplus \C_X,
$$
where $\cO(\Si)$ is the complex line bundle on $X$ with $c_1(\cO(\Si))=\tn{PD}(\Si)$ and $\C_X$ is the trivial complex line bundle.

\medskip
\par\noindent
In the symplectic case, the hermitian metric on $T\Si$ and the standard hermitian metric on $\C_\cD$ give us a hermitian metric on $\pi^*T\Si\oplus\C_\cD.$ Via (\ref{LTdfn_e}), the latter can be extended to a hermitian metric $h_{\sss X}$ on the entire $TX(-\log \Si)$. The restriction of this metric to $X-\Si$ defines a complete hermitian metric $h_{\sss X-\Si}$ on $X\!-\!\Si$. 

\medskip
\par\noindent
In order to define relative SW invariants, we use the moduli spaces $\cM_{\etacyend}(M,\smfrak)$, or indeed their compactification, where $M\!=\!X\!-\!\Si$ is equipped with the cylindrical-end metric $h_{\sss X-\Si}.$
As mentioned in Section~\ref{SW4cyend_ss}, the moduli space $\cM_{\etacyend}(M,\smfrak)$ often has components of different dimensions. In order to choose a component of fixed dimension, we use a $\spinc$ structure $\mfs_{\sss X}$ on $X$ that restricts to $\mfsm$ on $M.$ By \cite[p.~94, Fact~1]{N1}, if $\mfsm$ does not admit such an extension over $\cD$, then $\cM_{\etacyend}(M,\smfrak)$ is empty.

\medskip
\par\noindent
Let $\mfs_{\sss X}\!=\!(\sb_{\sss X},\clm_{\sss X})$ be a  $\tn{spin}^c$ structure on $TX(-\log \Si)$. We say $\mfs_{\sss X}$ is {\it relatively canonical} if $\Upsilon^* \mfs_{\sss X}$ is the canonical $\tn{spin}^c$ structure of the complex vector bundle 
$$
\Upsilon^* TX(-\log\Si)\cong \pi^* T\Si\oplus\C_\cD.
$$
Therefore, if $\mfs_{\sss X}$ is relatively canonical, by~(\ref{CanW_e}),
\bEq{LCanW_e}
\aligned
&\Upsilon^*\sb^+_{\sss X} = \Lambda^{0,0}(\pi^*T^*\Si\oplus\C_\cD) \oplus \Lambda^{0,2}(\pi^*T^*\Si\oplus\C_\cD) \cong \C_\cD\oplus\pi^*\Lambda^{0,1}\Si,\\
&\Upsilon^*\sb^-_{\sss X} = \Lambda^{0,1}(\pi^*T^*\Si\oplus\C_\cD) \cong \C_\cD\oplus\pi^*\Lambda^{0,1}\Si.\\
\endaligned
\eEq
Observe that $\Upsilon^*\sb^\pm_{\sss X}$ coincide with the pullback of the canonical spinor bundle $\sb_{\Si,\tn{can}}$ on $\Si$ in Section~\ref{SW2_ss} and that $d(\mfs_{\sss X})=0$ by (\ref{DegS_e}). 
If $X$ is a symplectic manifold equipped with an $(\om_X,\Upsilon)$-compatible almost complex structure $J$, then $TX(-\log \Si)$ has the structure of a complex vector bundle and the canonical $\tn{spin}^c$ structure is globally well-defined.

\subsection{The main component}\label{RelSWModuli_ss}
Fix a relatively canonical \spinc structure 
$\mfs_{\sss X,\tn{can}}=(\sb_{\sss X,\tn{can}},\clm_{\sss X,\tn{can}})$ on $TX(-\log\Si);$ this will be the globally canonical \spinc structure if $X$ is a symplectic manifold. 
Every other $\tn{spin}^c$ structure $\mfs\!=\!\mfs_{\sss X,E}\!=\!(\sb_{\sss X,E},\clm_{\sss X,E})$ can be obtained by tensoring with a hermitian line bundle $E$ on $X$. In the neighborhood $\cD$ of $\Si$, we may assume 
$$
E_\cD= \pi_{\sss \cD,\Si}^*E_{\sss\Si}.
$$
By (\ref{LCanW_e}),  we have
\bEq{CompWU_e}
\Upsilon^*\sb^+_{\sss X,E}\cong \Upsilon^*\sb^-_{\sss X,E} \cong \pi^*\big( E \oplus (\kcan^*\otimes E)\big)=\pi^*(\sb_{\sss \Si,E}),
\eEq
that is, $\Upsilon^*\sb^\pm_{\sss E}|_{\cD-\Si},$ where $\cD-\Si\!\cong\![0,\infty)\times Y,$ 
coincide with the pullback of the spinor bundle on $\Si$ corresponding to $E_{\sss\Si}$ in Sections~\ref{SW3_ss} and \ref{SW4cyend_ss}. In this case, (\ref{DegS_e}) gives us
$$
d(\mfs_{\sss X,E}) \!=\begin{cases}
\tn{deg}(E_{\sss \Si})& \quad \tn{if}\quad \tn{deg}(L_{\sss\Si})\leq 0,\\
\tn{deg}(K_{\sss \Si}-E_{\sss\Si})& \quad \tn{if}\quad \tn{deg}(L_{\sss \Si})\geq 0\,.
\end{cases}
$$

\noindent
Associated to $\mfs$, we are interested in two types of moduli spaces.

\medskip
\par\noindent
(1) Let $\cM_{\eta_o}(X-\Si,\sfrak_{\sss X-\Si})$ denote the moduli space of monopoles on $X-\Si$ with finite energy on the end with respect to the induced $\spinc$ structure $\sfrak|_{X-\Si}$ on $X-\Si$, where $\eta_o$ 
is a perturbation with compact support, as in Section~\ref{SW4cyend_ss}.
Different $\spinc$ structures $\mfs$ and $\mfs'$ on $TX(-\log \Si)$ may have the same restriction on $X-\Si$. Therefore, for such $\mfs$ and $\mfs'$,  $\cM_{\eta_o}(X-\Si,{\sfrak}_{\sss X-\Si})$ will be the same as $\cM_{\eta_o}(X-\Si,\sfrak'_{\sss X-\Si})$. As we described in the introduction, by looking at the limits of the monopoles over the cylindrical end, we choose a component $\cM_{\eta_o}(X-\Si,\sfrak)$ of $\cM_{\eta_o}(X-\Si,{\sfrak}_{\sss X-\Si})$ corresponding to $\mfs$; here we assume that $(\Si\cdot\Si,\deg(L))\!\neq\!(0,0)$, otherwise there is no such component and we will have to use an adapted perturbation term instead. This alternative perturbation is described below.

\medskip
\par\noindent
(2) Consider an adapted perturbation 
$$
\eta_\nu\defeq\etacyend(\eta_o,\eta)=\eta_o+\be(\tcyend)\eta^+
$$
as in (\ref{AdaptedEta_e}), where 
$$
\eta=\mfi (*_{\sss Y}n), \qquad n=\piysi^*(\nu+\ov{\nu})
$$ 
comes from a holomorphic 1-form $\nu$, as in (\ref{om-pert_e}). In fact, on the cylindrical end we have
\bEq{ExplicitEtaPlus_e}
\eta^+=\mfi(\nd t \wedge n + *_{\sss Y} n).
\eEq
Let $\cM_{\eta_\nu}(X-\Si,\sfrak_{\sss X-\Si})$ denote the moduli space of monopoles on $X-\Si$ with finite energy on the end with respect to the induced $\spinc$ structure $\sfrak|_{X-\Si}$ on $X-\Si$.
As before, by looking at the limits of the monopoles over the cylindrical end, we choose a component $\cM_{\eta_\nu}(X-\Si,\sfrak)$ of $\cM_{\eta_\nu}(X-\Si,{\sfrak}_{\sss X-\Si})$ corresponding to $\mfs$. In this case, all the monopoles are irreducible.

\medskip
\par\noindent
In both cases, there are limiting maps
$$
\cM_{\eta_o}(X-\Si,\sfrak)\lra\cM(Y,\mfsye), \qquad
\cM_{\eta_\nu}(X-\Si,\sfrak)\lra\cM_\eta(Y,\mfsye),
$$
i.e., if $(A,\Phi)$ is a finite-energy solution of $\tn{SW}^4_{\eta_o}$ (resp. $\tn{SW}^4_{\eta_\nu}$), then, restricted to the cylindrical part $\cD-\Si\cong[0,\infty)\times Y,$ this solution identifies with the gradient flow line $(B(t),\Psi(t))$ of $\tn{CSD}$ (resp. $\tn{CSD}_\eta$) such that the exponentially-converging limit
\bEq{local-Limit_e}
\pi_{\sss Y,\Si}^*[A_\infty,\Psi_\infty]=
[B,\Psi]=\lim_{t\lra\infty} [B(t),\Psi(t)]
\eEq
is induced by the pullback of an irreducible solution $(A_\infty,\Psi_\infty)$ of $\SW^2$ (resp. $\SW^2_\nu$) for $\mfs_{\sss \Si,E}$. 
Recall that an irreducible monopole on $Y$ can be identified, up to gauge equivalence, with a divisor on $\Si.$
In the compactly supported case (1), every divisor in $\tn{Div}_{d(\mfs)}(\Si)$ can arise as a limit, but in the adapted case (2), 
the limiting map takes values in $S_{d(\mfs)}(\nu)$. \\

\bRm{Questions}
Since the \spinc structure $\mfs$, and thus the characteristic line bundle $L$, is defined over the entire $X$, it is natural to ask whether 

\medskip
\par\noindent
(Q1) \textit{there exists a pair $(A_{\sss X},\Phi_{\sss X}),$ consisting of a connection on $L$ and a section of $\sb^+$, such that}
$$
(A_{\sss X},\Phi_{\sss X})|_{X-\Si}=(A,\Phi) \quad \tn{and} \quad (A_{\sss X}|_{T\Si},\Phi_{\sss X}|_{\Si})=(A_\infty,\Psi_\infty).
$$
(Q2) \textit{there exists a system of equations, defined over the entire $X$, whose solutions are gauge equivalent to such pairs $(A_{\sss X},\Phi_{\sss X})$.}

\medskip
\par\noindent
Note that we already know from (\ref{local-Limit_e}) that $\Psi_\infty$ and $\Phi$ define a continuous section $\Phi_{\sss X}$ of $\sb^+$. We will come back to these questions in Section~\ref{LSW_s}.
\eRm

\noindent
In the following subsections, we summarize the relevant conclusions from the introduction and Sections~\ref{SW3_ss} and \ref{SW4cyend_ss} and extend them to prove Theorems~\ref{Relaive_thm}--\ref{SW-Sum_thm2}. For simplicity, we assume that $\Si$ is connected. If $\Si$ is not connected, some of the related numbers and structures should be treated component-wise.

\subsection{Dimension}\label{Dim_ss}
As seen in Section~\ref{SW4cyend_ss}, and according to (\ref{DimCyl4_eq}), 
if $b^+_{\sss X-\Si}>0$, for generic $\eta_o$, $\cM_{\eta_o}(X-\Si,\sfrak)$ and $\cM_{\eta_\nu}(X-\Si,\sfrak)$ are smooth (but probably not closed) manifolds of real dimension 
\bEq{DimGeneral_e}
\frac{c_1(A)^2-2\chi(X-\Si)-3\si(X-\Si)}{4}+\xi.
\eEq
While $\chi(X-\Si)$ and $\si(X-\Si)$ are topological quantities, $c_1(A)^2$ and $\xi$ depend  (at least) on the asymptotic behavior of the solution. In order to (topologically) fix the asymptotic behaviour, we have considered a $\spinc$ structure on $TX(-\log \Si)$, instead of $T(X-\Si)$, resulting in a characteristic bundle $L$ that is defined over the entire $X$. Below, we show that $c_1(A)^2$ is equal to the topological quantity $c_1(L)^2$ and calculate $\xi$ for each of $\cM_{\eta_o}(X-\Si,\sfrak)$ and $\cM_{\eta_\nu}(X-\Si,\sfrak)$.
\bLm{dim-Calc_lm1}
We have
$$
\aligned
{\dsw}=\dim_\R \cM_{\eta_o}(X-\Si,\sfrak)&=\frac{c_1(L)^2-2\chi(X-\Si)-3\si(X-\Si)}{4}+\frac{\Si\cdot \Si- 3\tn{sign}(\Si\cdot \Si)}{4}+d(\mfs)\\
&=\frac{(c_1(L)+\Si)^2-2\chi(X)-3\si(X)}{4} 
\endaligned
$$ 
\eLm
\bPf
We use the formula (\ref{DimGeneral_e}), identify $c_1(A)^2$ with the topological quantity $c_1(L)^2,$ and use the explicit computation of $\xi$ in \cite[(3.27)]{N1}. 
\\
\par\noindent
\textbf{Claim}. The integral $c_1(A)^2$ in (\ref{Integralc_1^2_e}) is equal to $c_1(L)^2$.\\

\noindent
\tn{Proof}.
Fix a $U(1)$-connection $A_{\sss X}$ on $L$ whose restriction to $T\Si$ is the limiting connection $A_\infty$ defined in (\ref{local-Limit_e}). 
Since $L$ is the characteristic line bundle of the $\spinc$ structure $\mfs$ defined on the entire $X$, $c_1(L)^2$ is a well-defined topological quantity that coincides with the integral
$$
\frac{-1}{4\uppi^2}\int_X F_{A_{\sss X}}\wedge F_{A_{\sss X}}.
$$
Restricted to $X\!-\!\Si$, we have $\big(A_{\sss X} -A\big)|_{\sss X-\Si}=a$
for some imaginary-valued $1$-form $a$. Since $X-\Si$ is open and dense in $X$, we have 
\bEq{c1toint_e}
\frac{-1}{4\uppi^2}\int_X F_{A_{\sss X}}\wedge F_{A_{\sss X}}= 
\frac{-1}{4\uppi^2}\int_{X-\Si} F_{A_{\sss X}}\wedge F_{A_{\sss X}}= 
\frac{-1}{4\uppi^2}\int_{X-\Si} F_{A}\wedge F_{A}+ 2 \nd a \wedge F_A + \nd a\wedge \nd a.
\eEq
Restricted to the neck, we can write $a=a(t)+ f\, \nd t$, where $a(t)$ is a smooth $1$-parameter family of imaginary-valued $1$-forms on $Y$ and $f$ is a smooth imaginary-valued function. By (\ref{local-Limit_e}) and our assumption in the first line of the proof, $a(t)$ converges to $0$ exponentially fast. By Stokes' theorem, we have
$$
\int_{X-[t,\infty)\times Y} \nd a \wedge F_A = \int_{Y} a(t) \wedge F_A|_{(\{t\}\times Y)} \quad \tn{and} \quad  \int_{X-[t,\infty)\times Y} \nd a\wedge \nd a =  \int_{Y} a(t) \wedge \nd a(t).
$$
Both integrals converge to $0$ as $t$ converges to $\infty$. That finishes the proof of the claim.\qed

\noindent

\noindent
Back to the proof of Lemma~\ref{dim-Calc_lm1}, we use the calculation of $\xi$ in \cite{N1}. Let 
$$
\ve_{\sss \Si,X}\!=\!\tn{sign}(\Si\cdot \Si) \in \{-1,0,1\}.
$$ 
By the additivity of the signature, we have 
$$
\si(X-\Si)=\si(X)-\ve_{\sss \Si,X}.
$$
The number $\ell$ in \cite[(3.27)]{N1} is $\tn{deg}(Y)=-\Si\cdot\Si$, as explained in (\ref{degY_e}). Also, recall from Remark~\ref{DwN_rmk} that the twisting line bundle denoted by $L$ in \cite{N1} corresponds to $E_{\sss \Si}\otimes\kcan^{-1/2}$ in our setup. Thus the number $n$ in \cite[Rmk.~3.4]{N1} is $\deg(E_{\sss \Si})\!+\!(1-g)$ in our notation, which is equal to $\frac{1}{2}\deg(L_{\sss\Si})$. In conclusion, by \cite[(3.27)]{N1}, if $\tn{deg}(L_{\sss\Si})\leq 0$, we have 
\bEq{OldZeta_e}
\aligned
\xi= 
& \big(-\frac{\ve_{\sss \Si,X}+1}{2}+\tn{deg}(E_{\sss \Si})+(1-g)\big)+\big(g-\frac{1}{2}\big)+\frac{1}{4}\big(\Si\cdot \Si -\ve_{\sss \Si,X}\big)
\\
= & \deg(E_{\sss \Si}) +\frac{1}{4} (\Si\cdot \Si- 3\ve_{\sss \Si,X}).\\
\endaligned
\eEq
We conclude that, 
$$
\aligned
\mathsf{d} = &\frac{c_1(L)^2-2\chi(X)-3\si(X)}{4}+(1-g)+\frac{3}{4} \ve_{\sss \Si,X}-\frac{3}{4} \ve_{\sss \Si,X}+\frac{1}{4} \Si\cdot \Si+\tn{deg}(E_{\sss \Si})\\
=& \frac{(c_1(L)+\Si)^2-2\chi(X)-3\si(X)}{4}.
\endaligned
$$
Similar calculations give us the dimension formula for the case $\deg(L_{\sss\Si})\geq 0.$
\ePf

\noindent
For the component $\cM_{\eta_o}(X-\Si,\mc{J})$ ending at reducible solutions, \cite[(3.29)]{N1} yields the following.
\bLm{dim-Calc_lmRed}
If $(g-1)\not\equiv 0$ modulo $\ell\!=\!-\Si\cdot\Si$ and $\frac{1}{2}c_1(L_{\sss Y})$ is a torsion class $[k]\in H^2_{\tn{tor}}(Y,\Z)\cong \Z/\ell\Z$, with $0<k<|\ell|$, we have
\bEq{dJ_e}
\aligned
{\dsw}_{\mc{J}}=\dim_\R \cM_{\eta_o}(X-\Si,\cJ)=&\frac{c_1(A)^2-2\chi(X-\Si)-3\si(X-\Si)}{4}+\\
+&(g-\frac{1}{2})+\frac{1}{4} (\Si\cdot \Si-\ve_{\sss \Si,X})+\frac{k^2}{\Si\cdot \Si}-k\cdot \ve_{\sss \Si,X}.
\endaligned
\eEq
\eLm
\noindent
If $\Si$ is not connected, each component contributes by the factor in the second line.
If $(g-1)\equiv 0$ modulo $\ell$, then $k=0$ and the dimension formula will be slightly different. Note that, in the formula above, $A$ does not extend to a connection on $L$ over the entire $X$. If $\mfs_k$ is the logarithmic $\spinc$ structure on $TX(-\log \Si)$ for which $\frac{1}{2}\tn{deg}(L_{\mfs_k}|_{\Si})=k$, then \cite[(2.13)]{N1} and an integral formula similar\footnote{with $\lim_{t\to \infty} a(t)=\frac{2\mfi k}{\ell}\al$ instead of $0$.} to (\ref{c1toint_e}) implies
$$
c_1(A)^2-c_1(L_{\mfs_k})^2=-\frac{4k^2}{\Si\cdot \Si}.
$$ 
We can therefore rewrite (\ref{dJ_e}) as the topological formula
\bEq{dJ_e2}
\aligned
{\dsw}_{\mc{J}}=&\frac{c_1(L_{\mfs_k})^2-2\chi(X-\Si)-3\si(X-\Si)}{4}+\\
+&(g-\frac{1}{2})+\frac{1}{4} (\Si\cdot \Si-\ve_{\sss \Si,X})-k\cdot \ve_{\sss \Si,X}.
\endaligned
\eEq
We will use this in Section~\ref{Compact_ss} to explain why $\mc{J}$ does not contribute to the sum formula~(\ref{SW-Sum_e}).

\medskip
\par\noindent
Finally, we prove the following result for the dimension of $\cM_{\eta_\nu}(X-\Si,\sfrak).$
\bLm{dim-Calc_lm2}
We have
$$
\aligned
\wt\dsw\defeq \dim_\R \cM_{\eta_\nu}(X-\Si,\sfrak)=
&\frac{c_1(L)^2-2\chi(X-\Si)-3\si(X-\Si)}{4}+\frac{1}{4} (\Si\cdot \Si- 3\ve_{\sss \Si,X})=\\
=&\frac{(c_1(L)+\Si)^2-2\chi(X)-3\si(X)}{4}-\nd(\mfs).
\endaligned
$$ 
\eLm

\noindent
The explicit computation of $\xi$ in \cite[(3.27)]{N1} does not apply to monopoles in $\cM_{\eta_\nu}(X-\Si,\sfrak)$, because starting at \cite[Sec~3.3]{N1}, the author assumes $\psi_-$ is zero. In the following, we adapt the calculations in \cite{N1} to monopoles in $\cM_{\eta_\nu}(X-\Si,\sfrak)$ and prove the following proposition. 
With (\ref{Newzeta_e}) in place of (\ref{OldZeta_e}), the proof of Lemma~\ref{dim-Calc_lm2} is similar to the proof of Lemma~\ref{dim-Calc_lm1}.

\bPr{xi-non-compact_prp}
The $\xi$-invariant of the monopoles in $\cM_{\eta_\nu}(X-\Si,\sfrak)$ is equal to 
\bEq{Newzeta_e}
\xi=\frac{1}{4} (\Si\cdot \Si- 3\ve_{\sss \Si,X}).
\eEq
\ePr

\subsection{Proof of Proposition~\ref{xi-non-compact_prp}}\label{PfPrp_ss}
In this subsection, in order to prove (\ref{Newzeta_e}), we follow and modify the relevant details of \cite{N1}. To keep the proof short, and to enable the reader to directly compare with \cite{N1}, we use the same setup and notation in \cite{N1} with slight simplifications. 
In particular, our 3-manifold $Y$ is denoted by $N$ in this subsection.

\medskip
\par\noindent
The calculation of the $\xi$-invariant and thus the dimension formula \cite[(3.26)]{N1} involves contributions from several spectral flow calculations in \cite[Sec.~2.2-3.4]{N1}.
Up to \cite[Sec.~3.3]{N1}, the results apply to an arbitrary spinor $\phi=(\phi_-,\phi_+)$. Starting in \cite[Sec.~3.3]{N1}, the author assumes $\phi_-\!=\!0$; that will have major impacts on some of the calculations. Compared to \cite[(3.26)]{N1}, we get a different value for $\tn{SF}_+$ and the dimension of the asymptotic limit set is zero. In order to find $\tn{SF}_+$ in our case, we similarly want to calculate the spectral flow of 
\bEq{cHt_e}
\wt\cH_t=\wt\cH_0+t\cP_{\phi}\qquad t\in [0,1],
\eEq
but this time, $\phi\!=\!\phi_-\oplus\phi_+$ has two non-zero components in $\mc{K}^{-1/2}\otimes L$ and $\mc{K}^{1/2}\otimes L$ with $\ov\phi_-\phi_+=\nu\in H^0(\mc{K})$. Here $\mc{K}$ and $L$ are pullback to $N$ of the canonical line bundle and an arbitrary complex line bundle on $\Si$, respectively. The operator $\wt\cH_0$ is the same but $\cP_{\phi}=\cP_{\phi_+}+\cP_{\phi_-}$ is different. As a result, the resonance matrix $\cR_\phi$ used at the top of \cite[p.~97]{N1} is different, but it is defined on the same space $\tn{Ker}~\wt\cH_0$. Also, in this case $\tn{Ker}~\wt\cH_1\!=\!0$ since $(A,\phi)$ is a regular point. However, similarly to \cite[p.~97, STEP~1]{N1},  there is no spectral flow along the path $(t\lra \wt\cH_{t})_{t\in (0,1]}$. The only contribution is from $t\!=\!0$.

\medskip
\par\noindent
Before moving forward to find the generalization of \cite[Lemma~3.2]{N1} and calculate the spectral flow contribution $\tn{SF}_+$ at $t\!=\!0$, let us set up the notation, recall the definition of $\wt\cH_0$, and find the generalization of $\cP_{\phi}$. These operators are linear maps from 
\bEq{Conf-Space_e}
\Gamma(\mathbb{S}_{L}\oplus \mfi \Lambda^1 T^*N \oplus \Lambda^0 T^*N)
\eEq
to itself, where 
$$
\mathbb{S}_{L}= \mathbb{S}^-_{L}\oplus \mathbb{S}^+_{L}=\big(\mc{K}^{-1/2}\otimes L\big) \oplus \big(\mc{K}^{1/2}\otimes L\big).
$$
We have $\phi=\phi_-\oplus\phi_+\in \mathbb{S}_{L}$ and an arbitrary element in (\ref{Conf-Space_e}) is written\footnote{Compared to $\cite{N}$, we have simplified $\dot\psi$ and $\dot{a}$ to $\psi$ and $a$.} as 
$$
\Xi=(\psi_-\oplus \psi_+) \oplus \mfi a \oplus \mfi f.
$$
Furthermore, we write
$$
\mfi a= \mfi h \varphi + \frac{\om -\ov\om}{2}
$$
where $\{\varphi,\varphi_1,\varphi_2\}$  is the local orthonormal coframe of $T^*N$ in \cite[Sec~2.1]{N1}, and $\om$ is a complex multiple of the $(1,0)$ form $\frac{1}{\sqrt{2}}(\varphi_1+\mfi \varphi_2)$ as in \cite[Appendix~D]{N1}. At the monopole $(A,\phi)$, the operators $\wt\cH_0$ and $\cP_\phi$ are given by
$$
\wt\cH_0=
\begin{bmatrix} 
D_A & 0 & 0 \\
0& -*\nd & \nd\\
0 & \nd^* & 0 \\
\end{bmatrix}
\quad\tn{and}\quad 
\cP_\phi(\Xi)=
\begin{bmatrix} 
c(\mfi a)\phi- \mfi f \phi \\
\dot{q}(\phi,\psi)\\
\mfi \tn{Im}\ll \phi,\psi\rr \\
\end{bmatrix}.
$$
Using the same calculations as in \cite[(D.1)-(D.3)]{N1}, we have
$$
\mc{P}_\phi\Big(\Xi\Big)
=  (\Psi_-\oplus \Psi_+) \oplus \Big(\mfi H \varphi+   \frac{\Om -\ov\Om}{2}\Big)\oplus \mfi F,
$$
where
$$
\begin{bmatrix} 
\Psi_-\oplus \Psi_+  \\
\mfi H \varphi+  \frac{\Om -\ov\Om}{2}\\
\mfi\,F \\
\end{bmatrix}=
\begin{bmatrix} 
(-\ov{\om} \phi_+- (h+\mfi f) \phi_-) \oplus (-\om \phi_-+ (h-\mfi f )\phi_+)  \\
\mfi \tn{Re}(\ll \phi_+,\psi_+\rr-\ll \phi_-,\psi_-\rr)\varphi- \frac{1}{\sqrt{2}}\Big((\ov{\psi_-}\phi_++\ov{\phi_-}\psi_+)-({\psi_-}\ov{\phi_+}+{\phi_-}\ov{\psi_+})\Big) \\
\mfi\, \tn{Im}\Big(\ll \phi_+,\psi_+\rr+\ll \phi_-,\psi_-\rr\Big) \\
\end{bmatrix}.
$$

\noindent
The calculation of the spectral flow contribution at $t\!=\!0$ of the path $\{\wt\cH_t\}_{t\in [0,1]}$ is done in the following way in \cite[Step~2]{N1}. The spectral data of $\wt\cH_t$ can be organized in families depending analytically on $t$. Denote by $Z$ the set of all pairs $(\Xi(t),\la(t))$ where $\la(t)$ is a real eigenvalue of $\wt\cH_t$, $\Xi(t)$ is an eigenvector of unit length corresponding to $\la(t)$, $\la(0)=0$, and $t\longmapsto (\Xi(t),\la(t))$ is analytic. We have $\#Z=\tn{dim}_\R \tn{Ker}~\wt\cH_0$ and 
\bEq{Xit_e}
(\Xi(t);\la(t))=(\Xi_0+t\Xi_1+\ldots~; \la_mt^m+\la_{m+1}t^{m+1}+\ldots), 
\eEq
where the integer $m$ is the order of the pair $(\Xi(t),\la(t))$. A pair is called degenerate if $m>1$ and non-degenerate otherwise. Since $\tn{Ker}~\wt\cH_t=0$ for $t\!>\!0$, all pairs have finite order (so in our case $Z=Z^*$). The spectral flow $\tn{SF}_+$ is then determined by
$$
-\tn{SF}_+=\#\{(\Xi(t),\la(t)\colon~\la_m<0\}.
$$
\par\noindent
For $(\Xi(t),\la(t))$ in (\ref{Xit_e}), the equation $\wt\cH_t \Xi(t)=\la(t)$ expands to  
\bEq{Expansion_e}
\begin{cases}
\wt\cH_0(\Xi_0)&=0,\\
\wt\cH_0(\Xi_1)+\cP_\phi (\Xi_0)&=0,\\
&\vdots\\
\wt\cH_0(\Xi_{m-1})+\cP_\phi (\Xi_{m-2})&=0~,\\
\wt\cH_0(\Xi_m)+\cP_\phi (\Xi_{m-1})&=\la_m \Xi_0,\\
&\vdots
\end{cases}
\eEq
If $|\Xi_0|=1$, note that 
\bEq{laEq_e}
\la_m=\ll \la_m \Xi_0,\Xi_0\rr=\ll \cP_\phi (\Xi_{m-1}),\Xi_0\rr.
\eEq
We have
$$
\tn{Ker}~\wt\cH_0\cong \ov{H^0(K^{1/2}\otimes L^*)}\oplus H^0(K^{1/2}\otimes L) \oplus  H^0(\Si,\R)\oplus H^1(\Si,\R),
$$
where we have identified $H^0(K^{-1/2}\otimes L)$ with $\ov{H^0(K^{1/2}\otimes L^*)}$. We will also identify $H^0(\Si,\R)$ with $\R$ and $ H^1(\Si,\R)$ with the space of holomorphic $1$-forms $\om$ via 
$$
H^0(\Si,\R)\ni a\lra \om \quad \tn{s.t.}\quad \mfi a = \frac{\om-\ov\om}{2}.
$$
For each $\Xi_0\!\in\!\tn{Ker}~\wt\cH_0,$ let $\cR_\phi(\Xi_0)$ denote the projection to $\tn{Ker}~\wt\cH_0$ of $\cP_\phi(\Xi_0)$. The operator $\cR_\phi$ is called the \textit{resonance} form/matrix. For 
$$
\Xi_0= (\psi_-\oplus \psi_+)\oplus \frac{\om-\ov\om}{2}\oplus \mfi f\in \tn{Ker}~\wt\cH_0,
$$ 
the quadratic form
$$
Q(\Xi_0)\defeq \ll \mc{R}_{\phi} \Xi_0,\Xi_0\rr= \ll \mc{P}_{\phi} \Xi_0,\Xi_0\rr=\sqrt{2}\,f\,\tn{Im} \Big(\ll \phi_+,\psi_+\rr + \ll \phi_-,\psi_-\rr\Big)- 
\tn{Re}\Big( \big({\ov\psi}_-\phi_+  + {\psi}_+\ov \phi_-\big) \ov \om\Big)
$$
generalizes \cite[Lemma~3.2]{N1}. Non-degenerate pairs correspond to the case where $m=1$ in (\ref{Expansion_e}), whereas degenerate pairs correspond to the kernel of $\cR_\phi$ for which $m>1$. The contribution to the spectral flow of the non-degenerate pairs is equal to the number of negative eigenvalues of the quadratic form $Q$. In the following, we find the number of negative eigenvalues of $Q$ and the kernel of $\cR_\phi$.

\medskip
\par\noindent
Let 
$$
V_-=\ov{H^0(K^{1/2}\otimes L^*)},\quad V_+=H^0(K^{1/2}\otimes L), \quad  K=H^1(\Si,\R).
$$ 
Consider the multiplication maps
$$
\aligned
&m_+ \colon V_- \lra K, \quad  {\psi}_-\lra {\ov\psi}_-\phi_+ ,\\
&m_- \colon V_+ \lra K,\quad   \psi_+\lra {\psi}_+ \ov\phi_-\,.
\endaligned
$$
The second term in $Q$ takes the form
$$
(f,v^+,v^-,\om)\lra -\tn{Re} \big((m_+ ({\psi}_-)+m_-({\psi}_+))\ov\om \big).
$$
Let $R_\pm$ denote the range of $m_\pm$. Since $\phi_+\ov\phi_-=\nu$ ($\nu$ is the perturbation term), we get 
$$
R_+\cap R_-= \C\cdot \nu.
$$
Decompose $K$ into $\C\cdot \nu \oplus W_\pm \oplus W_0$, where $W_\pm$ is the intersection of orthogonal complement of $\C\cdot \nu$ with $R_++R_-$ and $W_0$ is the orthogonal complement of $R_++R_-$. 
As in \cite[p.~99-100]{N1}, we can choose complex bases $v^\pm_1,\ldots,v^\pm_{d_{\pm}}$ 
for $V_\pm$ such that
\bIt
\item $v_1^\pm=\phi_\pm$,
\item for $v^\pm =\sum_{i=1}^{d_\pm} z_i^\pm v^\pm_i$, the first term in $Q$ takes the form 
$$
(f,v^+,v^-,\om)\lra f \tn{Im}(z_1^++z_1^-),
$$
\item under $m_\mp$, $\sum_{i=2}^{d_\pm} z_i v^\pm_i$ has image in $W_\pm$.
\eIt
Then, for 
$$
\om = \la \nu \oplus \om_\pm\oplus \om_0\in \C\cdot \nu \oplus W_\pm \oplus W_0,
$$
we have 
$$
Q(f,v^+,v^-,\om)=  f \,\tn{Im}(z_1^++z_1^-) - \tn{Re} \big( (\ov{z_1^-}+z_1^+)\ov{\la} \big) |\nu|^2 - \tn{Re}\ll\cdot,\cdot\rr|_{W_\pm \times W_\pm},
$$
where the last term denotes the quadratic form $\tn{Re}\ll\cdot,\cdot\rr$ on $W_{\pm}$.
We conclude that 
$$
\aligned
\tn{dim}_\R \tn{Ker}~Q&= 1+2~\tn{dim}_\C W_0= 3+ 2(g-(d_+-1)-(d_--1)-1)=2(g-(d_++d_--1))+1 \\
\tn{dim}_-~Q&= 3 + 2~\tn{dim}_\C W_\pm= 2(d_-+d_+-1)+1.
\endaligned
$$
Therefore, (A) the spectral flow contribution of the non-degenerate pairs is 
$$
-2(d_-+d_+-1)-1, 
$$
and, (B) the number of degenerate pairs is equal to 
$$
2(g-(d_++d_--1))+1.
$$
Since $\dim_\R \tn{Ker}~\wt\cH_t = 0$ for $t>0$, there can be at most 
$$
x<\tn{dim}_\R \tn{Ker}~Q=2(g-(d_++d_--1))+1
$$
degenerate pairs contributing to the spectral flow. Thus, we get
\bEq{SF_e}
-\tn{SF}_+=x+2(d_-+d_+-1)+1.
\eEq
\bLm{SF_lm} 
We have
\bEq{xfactor_e}
x=\begin{cases}
(g-(d_++d_--1))+1 & \tn{if}\quad \ell>0;\\
(g-(d_++d_--1))+0 & \tn{if}\quad  \ell<0.
\end{cases}
\eEq
Therefore,
\bEq{SF_e2}
-\tn{SF}_+=\begin{cases}
g+d_++d_- +1& \tn{if}\quad \ell>0;\\
g+d_++d_- & \tn{if}\quad  \ell<0.
\end{cases}
\eEq
\eLm

\medskip
\par\noindent
First, we use Lemma~\ref{SF_lm} to finish the proof of Proposition~\ref{xi-non-compact_prp}.
Moving to \cite[Sec~3.4]{N1}, in (3.21), the first term involves $\wt\cH_0$ so it is as considered there, the second term is (\ref{SF_e}), and the third term (i.e., $\tn{SF}(\wt\cH_1\lra \cO_w)$) vanishes for the same reason\footnote{We get the same identity $\int_{N} |\nd f|^2+|f|^2|\phi|~\nd \tn{vol}_N=0$ as in \cite[p.~121]{N1}.} as in \cite[Cor.~D.3]{N1}. In \cite[(3.25)]{N1}, replacing $-\tn{SF}_+$ with (\ref{SF_e2}), since the dimension of the asymptotic moduli space is zero, we get
$$
\aligned
\xi&=-(d_-+d_ - +\tn{SF}_+)+ -g -\frac{1}{2}+\frac{1}{4}(\Si\cdot \Si- \ve_{\sss \Si,X})\\
&=g+\frac{1-\ve_{\sss \Si,X}}{2}+ -g -\frac{1}{2}+\frac{1}{4}(\Si\cdot \Si- \ve_{\sss \Si,X})\\
&=\frac{1}{4} (\Si\cdot \Si- 3\ve_{\sss \Si,X}).
\endaligned
$$
That finishes the proof of Proposition~\ref{xi-non-compact_prp}.\qed

\medskip
\par\noindent
It is just left to prove Lemma~\ref{SF_lm}, which comes next.
\bPf
Compared to the degenerate case in \cite[p.~101]{N1}, in our case, $\tn{Ker}~\wt\cH_1=0$, and $\dim_\R \tn{Ker}~\cR=2(g-(d_++d_--1))+1$;  $\tn{Ker}~\cR_\phi$ is the direct sum of 
\bEn
\item the complex linear subspace $W_0$, and 
\item the $1$-dimensional linear subspace generated by $e_0= \phi_+ \oplus -\phi_- $.
\eEn

\noindent
For $\Xi_0,\Xi_0'\!\in\!\tn{Ker}(\cR_\phi)$, the bilinear form
$$
\cB(\Xi_0,\Xi'_0)=\ll \Xi_0,\cP_{\phi}\Xi'_1\rr=-\ll \wt\cH_0\Xi_1,\Xi'_1\rr
$$
in \cite[p.~102]{N1} is well-defined and symmetric. These sub-spaces (1) and (2) above are orthogonal with respect to $\mc{B}$. Restricted to $W_0$, for $\Xi_0=\Xi_0(\om)= 0 \oplus \frac{\om-\ov\om}{2}\oplus 0$, we have 
$$
\mc{P}_{\phi}(\Xi_0)=\begin{bmatrix} 
-\ov{\om} \phi_+ \oplus -\om \phi_-  \\
0 \\
0 \\
\end{bmatrix}.
$$
Therefore, the $\psi$-component $\psi^\om=\psi_-^\om\oplus \psi_+^\om$ of the next term $\Xi_1$ in the expansion (\ref{Xit_e}) satisfies
$$
\dbar\psi_+^\om= \ov\om \phi_+ \quad\tn{and}\quad \partial\psi_-^\om=\om \phi_-,
$$
where $\dbar$ and $\partial$ are the operators at the bottom of \cite[p.~79]{N1}. We conclude that 
$$
\aligned
-\mc{B}(\Xi_0(\om_1),\Xi_0(\om_2))=\tn{Re} \ll \psi^{\om_1},D_A~\psi^{\om_2}\rr &=\tn{Re}\Big(\ll \psi_-^{\om_2}, \dbar\psi_+^{\om_1} \rr+\ll \psi_-^{\om_1}, \dbar\psi_+^{\om_2} \rr\Big)\\
&= \tn{Re}\Big(\ll \psi_-^{\om_2}, \ov\om_1 \phi_+\rr+\ll \psi_-^{\om_1}, \ov\om_2 \phi_+ \rr\Big).
\endaligned
$$
Fix a complex basis $\om_1,\ldots,\om_k$ for $W_0$. Consider the complex $k\times k$ symmetric matrix 
$$
X+\mfi Y=\Big(\ll \psi_-^{\om_j}, \ov\om_i \phi_+\rr+\ll \psi_-^{\om_i}, \ov\om_j \phi_+ \rr\Big)_{i,j \in \{1,\ldots, k\}}.
$$
Since 
$$
\psi_-^{\mfi\om}\oplus \psi_+^{\mfi\om}=\mfi \psi_-^\om\oplus -\mfi \psi_+^\om
$$
with respect to the real basis $\{e_0, \om_1, \mfi \om_1, \ldots, \om_k, \mfi\om_k\}$ for $\tn{Ker}~\cR_\phi$, the matrix $\cB$ is given by
\bEq{XY-matrix_e}
\begin{bmatrix} 
\al &0&0\\
0&X & -Y\\
0&-Y & -X
\end{bmatrix}.
\eEq
The real number $\al$ depends on $\ell$ in the following way. Similarly to \cite[(3.16)-(3.18)]{N1}, corresponding to $\Xi_0=e_0$ the solution $\Xi_1=\psi\oplus \mfi {a}\oplus \mfi f$ of the following set of equations 
$$
\begin{cases}
&D_A\psi=0\\
&-*d{a}+\nd f + (|\phi_+|^2+|\phi_-|^2)\varphi =0\\
&d^*{a}=0
\end{cases}
$$
satisfies $\wt\cH_0\Xi_1+ \cP_{\phi} \Xi_0=0$. The above equation has a unique solution for which
$$
{a}=\frac{1}{2\ell} (|\phi_-|^2+|\phi_+|^2)\varphi.
$$
Furthermore, 
$$
\mc{B}(\Xi_0)=-\frac{1}{2\ell} \int_{N} |\phi_-|^2+|\phi_+|^2 .
$$
We conclude that $\al$ is negative if and only if $\ell\!>\!0$. Also, if 
$\begin{bmatrix} 
0\\
u\\
v
\end{bmatrix}$ 
is an eigenvector of (\ref{XY-matrix_e}) with eigenvalue $\la$, then 
$\begin{bmatrix} 
0\\
-v\\
u
\end{bmatrix}$ 
is another eigenvector with eigenvalue $-\la$. 
Similarly to the non-degenerate (order 1) case, the negative eigenvalues of $\cB$ correspond to the spectral flow contributions of order-2 eigenvectors in (\ref{Expansion_e}). If $\mc{B}$ has no zero eigenvalue, this finishes the proof of Lemma~\ref{SF_lm}. Otherwise, we need to continue inductively in the following way.\\

\noindent
Inductively, starting with $\cR^1_\phi=\cR_\phi$, for $\Xi_0\in \tn{Ker}~\cR^k_\phi,$ let $\cR^{k+1}_\phi(\Xi_0) \in \tn{Ker}~\cR^k_\phi$ denote the orthogonal projection of  $\cP_\phi(\Xi_k)$ to $\tn{Ker}~\cR^k_\phi$. Eigenvectors of order $k+1$ in (\ref{Expansion_e}) correspond to negative eigenvalues of the symmetric bilinear form
\bEq{Orderkpairing_e}
\cB^k(\Xi_0,\Xi'_0)=\ll \Xi_0,\cP_{\phi}\Xi'_k\rr
\eEq
and eigenvectors of order $>\!k+1$ in (\ref{Expansion_e}) correspond to $\tn{Ker}~\cR^{k+1}_\phi$. For $k\!\geq\!2$, $\tn{Ker}~\cR^k_\phi$ is a complex linear subspace of $W_0$. For the same reason as above, the number of positive and negative eigenvalues of (\ref{Orderkpairing_e}) are the same. This inductive process ends in finite steps.
\ePf

\subsection{Tunneling spaces}\label{Tunneling_ss}
Suppose $\pi\colon\cN\lra \Si_-$ is a hermitian line bundle over a genus $g$ surface and $P$ is the projectivization of $\cN$. Let $\Si_+$ denote the section at infinity. We have
$$
  \ell\defeq\tn{deg}(\cN)=\Si_-^2=-\Si_+^2\quad  \tn{and}\quad c_1(K_P)= -\tn{PD}(\Si_-+\Si_+)+(2g-2) \tn{PD}(F),
$$
where $F$ is the fiber class, and the self-intersections numbers on the left are taken in $P$. The second homology of $P$ is generated by $\Si_\pm$ and $F$ satisfying the relation 
$$\Si_+= \Si_- - \ell F.$$
Let $\Si=\Si_-\cup\Si_+$. We have
$P-\Si=\R\times Y,$ where $Y$ is the unit circle bundle in $\cN.$ 
We may think of $P-\Si$ as a 4-manifold with two cylindrical ends modeled over the dual $S^1$-bundles $Y$ at $+\infty$ and $Y^*$ at $-\infty$. Since $P$ is a K\"ahler surface, every $\spinc$ structure $\mfs$ on $TP(-\log \Si)\cong \pi^* T\Si_{\pm}\oplus \C_{P}$ corresponds to a homology class 
\bEq{HomologyClass_e}
\sz= [a \Si_- + b_+ F]=[a \Si_+ + b_- F],\qquad b_--b_+=a\ell.
\eEq
The logarithmic canonical bundle of $(P,\Si)$ is 
$$
K_P(\log \Si)=K_P \otimes \cO(\Si)=\pi^*K_\Si.
$$

\noindent
The moduli space $\cM(P-\Si,\sfrak)$ has expected dimension
\bEq{d-Tunneling_e}
\dsw_{\sz}=\sz\cdot\sz-\sz\cdot K_P= (a+1)(b_-+b_+)+2a(1-g)
\eEq 
and admits a natural $\R\!\times \!S^1$-action. The quotient with respect to $\R$ is the unparametrized space $\breve\cM(C_-,C_+)$ of tunnelings between $C_-$ and $C_+$ described in Section~\ref{SW4cyend_ss}, where $C_\pm$ depends on the number $b_\pm$. The moduli space is empty unless
$$
0\leq \sz\cdot \Si_\pm = b_\pm \leq 2g-2 .
$$
If $a=0$, then $b_+=b_-$ and $\cM(X-\Si,\sfrak)\cong \tn{Sym}^{b_+}(\Si_-)\cong  \tn{Sym}^{b_-}(\Si_+)$ consists of only fiber-wise constant solutions. The action of $\R\times S^1$ is trivial in this case.\\

\noindent
Similarly, we are interested in perturbed moduli spaces $\cM_{\eta_\nu}(P-\Si,\sfrak).$
In this case, the expected dimension is 
\bEq{d-Tunneling_e2}
\wt\dsw_{\sz}=\sz\cdot \sz-\sz\cdot K_P(\log \Si)= a(b_-+b_+ +2(1-g)).
\eEq
The quotient with respect to $\R$ is the unparametrized space $\breve\cM(C_-,C_+)$ of tunnelings between two points $C_-\in S_{b_-}(\nu)$ and $C_+\in S_{b_+}(\nu)$. 

\begin{remark}
By Thom isomorphism, we have $H^2_c(P-\Si)\cong H^1(Y)$. It follows that $b^+_{\sss P-\Si}=0$. Therefore, the example above does not fit in the framework of Theorem~\ref{Relaive_thm}.
\end{remark}

\bPr{Tunneling_prp}
Suppose $(A,\Phi) \in \cM(P-\Si,\sfrak)$ or $\cM_{\eta_\nu}(P-\Si,\sfrak)$. Then the connection $A$ defines a holomorphic structure on the twisting complex line bundle $\cO(\sz)$, $\Phi_+$ is a holomorphic section of $\cO(\sz),$ $\ov\Phi_-$ is a holomorphic section of $K_P(\log \Si)\otimes \cO(-\sz)$, and $\Phi_+\ov\Phi_-=\mfi \pi^*\nu \in K_P(\log \Si)\cong \pi^* K_\Si$.
\ePr

\bPf
Comparing the notations in Section~\ref{SW3_ss} and \ref{LogStr_ss}, the connection $1$-form $\al$ on $Y$ is the same as $\al$ after (\ref{Decomp_e}). By (\ref{cNHolo_e}), it defines an integrable almost complex structure $J$ on $\cN$ and thus on $\R\times Y$. On $\R\times Y$, $J$ acts by the complex structure of $\Si$ on the horizontal subspace of every slice $\{t\}\times Y$, and $J\partial_t=\partial_\theta$. If $\ell\neq 0$, the cylindrical metric is not K\"ahler: the associated $(1,1)$-form is $\nd t\wedge \al +\pi^*\om$ and $\nd(\nd t\wedge \al +\pi^*\om)= -\nd t \wedge \nd \alpha \neq 0$ because $\mfi \nd \al$ is the curvature of the principal $U(1)$-bundle $Y$.

\medskip
\par\noindent
On $\R\times Y$, using the same notation as in the proof of Lemma~\ref{Nic-lm}, for 
$$
(A,\Phi)=\big(B(t),\Psi(t)\big)_{t\in\R}\, , \quad\tn{where}\quad \Psi(t)=\big(\Psi_+(t),\Psi_-(t)\big),
$$ 
the Dirac operator has the form\footnote{Note that $\la=0$ because we are using the adiabatic connection on $Y$.}
\bEq{DiracOnRxS1bundle_e}
\dirac_A=
\left[ 
  {\begin{array}{cc}
   \nabla_{-\partial_t-\mfi\partial_\theta} & \sqrt{2}\dbar_{B(t)}^* \medskip \\
   \sqrt{2}\dbar_{B(t)} & \nabla_{-\partial_t + \mfi\partial_\theta}  \\
  \end{array}}
\right].
\eEq
Therefore,
\bEq{SW3Dirac_e2}
\dirac_A\Phi =0 ~~\Leftrightarrow ~~\begin{cases}
\sqrt{2}\dbar_{B(t)}\Psi_+(t)-\nabla_{\partial_t-\mfi\partial_\theta}\Psi_-(t)=0\\
\sqrt{2}\dbar_{B(t)}^*\Psi_-(t)-\nabla_{\partial_t+\mfi\partial_\theta}\Psi_+(t)=0
\end{cases}
\eEq
By \cite[(4.9)]{KM}, the unperturbed curvature equation has the form
\bEq{CSD-Curvature_e}
*_{Y} \dot{B}(t)+ F_{B(t)}=(\Psi(t)\Psi(t)^*)_0.
\eEq
By (\ref{ExplicitEtaPlus_e}), the curvature equation, when perturbed using $\nu$, has the form
\bEq{Per-CSD-Curvature_e}
*_{Y} (\dot{B}(t)-\mfi n)+ (F_{B(t)}-*_{\sss Y} \mfi n)=(\Psi(t)\Psi(t)^*)_0.
\eEq
Note that $\dot{B}(t)$ is a family of imaginary valued $1$-forms $\beta(t)$ on $Y$. The curvature equation (\ref{Per-CSD-Curvature_e}) now reads:
\bEq{CurvEq_e2}
\aligned
-\mfi(F_{\al_1,\al_2,t}+\beta_{\al,t}) & = \frac{1}{2}(|\Psi_+(t)|^2-|\Psi_-(t)|^2) \\
\big(\frac{1}{2}(F_{\al_2,\al,t}-\mfi F_{\al_1,\al,t})+\beta_{(0,1),t}\big)\al_{0,1} -\mfi\pi^*\ov\nu& = \overline{\Psi}_+(t)\Psi_-(t),
\endaligned
\eEq
where $\al_{0,1}$ is defined as in (\ref{alhol_dfn}) and 
$$
\aligned
F_{B(t)}&= F_{\al_1,\al_2,t}~\al_1\wedge\al_2 + F_{\al_1,\al,t}~\al_1 \wedge \alpha+ F_{\al_2,\al,t}~\al_2 \wedge \alpha,\\
\beta(t)&= \beta_{\al,t} \al+ \beta_{(0,1),t}\, \al_{0,1} - \ov{\beta_{(0,1),t}}\, \al_{1,0}.
\endaligned
$$

\noindent
Applying $\dbar_{B(t)}$ to the second equation in (\ref{SW3Dirac_e2}) we get
\bEq{Start_Le2}
-\dbar_{B(t)}\nabla_{\partial_t+\mfi\partial_\theta}\Psi_+(t)+\sqrt{2}\dbar_{B(t)}\dbar_{B(t)}^*\Psi_-(t)=0.
\eEq
Commuting the differential operators $\nabla_{\partial_t+\mfi\partial_\theta}$ and $\dbar_{B(t)}$ produces a curvature term
$$
\nabla_{\partial_t+\mfi\partial_\theta} ~ \dbar_{B(t)}-
\dbar_{B(t)} ~ \nabla_{\partial_t+\mfi\partial_\theta}
=\big(2\beta_{(0,1),t}+(F_{\al_2,\al,t}-\mfi F_{\al_1,\al,t})\big)\al_{0,1}
$$
and using the second equation in (\ref{CurvEq_e2}), we can rewrite (\ref{Start_Le2}) as  
$$
-\nabla_{\partial_t+\mfi\partial_\theta}\,\dbar_{B(t)}\Psi_+ 
+2\big(\ov{\Psi}_+(t)\Psi_-(t)+\mfi\pi^*\ov\nu\big)\Psi_+(t)
+\sqrt{2}\dbar_{B(t)}\dbar_{B(t)}^*\Psi_-(t)=0.
$$
Applying the first equation in (\ref{SW3Dirac_e2}) to the first term and dividing by $\sqrt{2}$ we get
\bEq{LastEq_e}
-\frac{1}{2}\nabla_{\partial_t+\mfi\partial_\theta}\,\nabla_{\partial_t-\mfi\partial_\theta}\Psi_-(t)
+\sqrt{2}\big(\ov{\Psi}_+(t)\Psi_-(t)+\mfi\pi^*\ov\nu\big)\Psi_+(t)
+\dbar_{B(t)}\dbar_{B(t)}^*\Psi_-(t)=0.
\eEq
Similarly, we get
\bEq{LastEq_e2}
-\frac{1}{2}\nabla_{\partial_t-\mfi\partial_\theta}\,\nabla_{\partial_t+\mfi\partial_\theta}\Psi_+(t)
+\sqrt{2}\big(\Psi_+(t)\ov\Psi_-(t)-\mfi\pi^*\nu\big)\Psi_-(t)
+\dbar^*_{B(t)}\dbar_{B(t)}\Psi_+(t)=0.
\eEq
Define 
$$
\dbar^{\tn{ver}}=\nabla_{\partial_t+\mfi\partial_\theta}
$$
Note that the first terms in (\ref{LastEq_e}) and (\ref{LastEq_e2}) are 
$$
\frac{1}{2} \dbar^{\tn{ver}} \dbar^{\tn{ver}\,*}\qquad \tn{and}\qquad \frac{1}{2} \dbar^{\tn{ver}\,*}\,\dbar^{\tn{ver}},
$$
respectively. Also note that the $\dbar$-operator associated to the connection $A$ over the entire complex manifold $\R\times Y$ is
$$
\aligned
\dbar_A = &\big((\nd t -\mfi\al)\otimes \frac{1}{2}\nabla_{\partial_t+\mfi\partial_\theta} \big) + \dbar_{B}\\
=&\big((\nd t -\mfi\al)\otimes \frac{1}{2}\nabla_{\partial_t+\mfi\partial_\theta} \big) + 
\big(\al_{0,1} \otimes \frac{1}{2}\nabla_{\ze_1+\mfi\ze_2} \big)
\colon \Om^0(\cO(\sz))\lra \Om^{0,1}\big(\cO(\sz)\big).
\endaligned
$$ 
Therefore, if $\Psi_+$ satisfies $\dbar^{\tn{ver}}\Psi_+=0$ and $\dbar_B \Psi_+=0,$ then it satisfies $\dbar_A\Psi_+=0$.
We have a similar statement for $\Psi_-$.
As in the proof of Lemma~\ref{Nic-lm}, we want to take the inner product of, say (\ref{LastEq_e2}), with $\Phi_+$ over the entire $\R\times Y$ and conclude that
\bEq{ZeroNess_e}
\dbar^{\tn{ver}}\Phi_+=0,\quad \dbar_A \Phi_+=0,\quad \tn{and}\qquad \Phi_+\ov\Phi_-=\mfi\pi^*\nu.
\eEq
Towards this goal, we will work over $[-T_1,T_2]\times Y$ and then let $T_i\lra\infty.$
To find the boundary terms, consider the complex $(2,1)$-form
$$
\Om=-\mfi \ll \Phi_+, \dbar^{\tn{ver}}\Phi_+\rr (\nd t+ \mfi \al) \wedge \pi^*\om.
$$
Since $\R\times Y$ is a holomorphic manifold, we have 
$$
\nd \Om = \dbar \Om = 2 \Big(\,\ll \dbar^{\tn{ver}} \Phi_+, \dbar^{\tn{ver}}\Phi_+\rr-\ll  \Phi_+, \dbar^{\tn{ver}*}\dbar^{\tn{ver}}\Phi_+\rr \,\Big)~\nd t\wedge \nd\tn{vol}_Y.
$$
By the same reasoning as in (\ref{VanishingPart_e}), and integration by parts corresponding to the 3-form above, we get
\bEq{IntegralT_e}
\aligned
0=&\int_{[-T_1,T_2]\times Y} \ll\Phi_+, \frac{1}{2}\dbar^{\tn{ver}\,*}\dbar^{\tn{ver}}\Phi_+  +\dbar_{B(t)}^*\dbar_{B(t)}\Phi_+ +\sqrt{2}(\Phi_+\ov\Phi_- -\mfi\pi^*\nu)\Phi_-\rr  \nd t\wedge\nd\tn{vol}_{Y}\\
=&\int_{[-T_1,T_2]\times Y} \Big(\frac{1}{2}|\dbar^{\tn{ver}}\Phi_+|^2 +|\dbar_{B(t)}\Phi_+|^2 +\sqrt{2}|\Phi_+\ov\Phi_- -\mfi\pi^*\nu|^2\Big) \nd t\wedge\nd\tn{vol}_{Y}\\
&-\frac{1}{4}\int_{\{T_2\}\times Y} \ll\Phi_+, \dbar^{\tn{ver}}\Phi_+\rr_{t=T_2} \nd\tn{vol}_Y
+\frac{1}{4}\int_{\{-T_1\}\times Y} \ll\Phi_+, \dbar^{\tn{ver}}\Phi_+\rr_{t=-T_1} \nd\tn{vol}_Y.
\endaligned
\eEq
The integral in the second line is non-negative and is an increasing function of each $T_i$. The integrals in the third line converge to $0$ as $T_i\lra\infty.$ 
Thus the integral in the second line must be zero. 
We conclude that $\dbar^{\tn{ver}}\Phi_+\!=\!0,$ $\dbar_{B(t)}\Phi_+\!=\!0,$ and $\Phi_+\ov\Phi_-\!=\!\mfi\pi^*\nu,$ from which (\ref{ZeroNess_e}) follows. 
Moreover, the curvature of $A$ will be a $(1,1)$-form by (\ref{CurvEq_e2}), so $A$ defines a holomorphic structure. 
Note meanwhile that the boundary terms in the third line turn out to be zero, because $\dbar^{\tn{ver}}\Phi_+=0$.
\ePf

\bCr{Empty}
Unless $b_-=b_+$, the moduli space $\cM_{\eta_\nu}(P-\Si,\sfrak)$ is empty. If $b_+=b_-=b,$ which happens when $a=0,$ then $\cM_{\eta_\nu}(P-\Si,\sfrak)\cong S_{b}(\nu)$ consists of fiber-wise constant solutions. 
\eCr

\bPf
The zero sets of $\Phi_+$ and $\ov\Phi_-$ represent the homology classes 
$$
[a\Si_-+b_+F]\qquad\tn{and}\qquad [-a\Si_-+(2g-2-b_+)F],
$$ 
respectively. Therefore, we must have $a\!=\!0$.  If $a=0$, the restrictions of $\Phi_+$ and $\ov\Phi_-$ to each fiber of $P$ are constant sections. Therefore, both are pullbacks from $\Si$.
\ePf

\begin{remark}
In the case of $\cM(P-\Si,\sfrak)$, if $a\!>\!0$, then $\Phi_-=0$ and if $a<0$, then $\Phi_+=0$.
Furthermore, if $a>0$, then $0\leq b_{\pm} \leq g-1$, and if $a<0$, then $g-1\leq b_{\pm} \leq 2(g-1)$. The two cases are related by Serre duality.
\end{remark}

\subsection{Compactification}\label{Compact_ss}
As we mentioned in the introduction and Sec.~\ref{SW4cyend_ss}, the moduli spaces $\cM_*(X-\Si,\sfrak)$ (where $*$ means either a compact perturbation or an adapted perturbation) are not necessarily compact. Because of the tunneling phenomenon (see \cite[Sec.~16]{KM} or \cite[Sec.~4.4.2]{N}), a sequence of monopoles in $\cM(X-\Si,\sfrak)$ will, after passing to a sub-sequence, ``converge" to a finite ordered set of monopoles where the first one is a monopole on $X-\Si$ and the rest are non-trivial monopoles on the cylinder $\R\!\times\! Y$. 
For an adapted perturbation $\eta_\nu,$ with $\nu\!\neq\!0,$ $\cM_{\eta_\nu}(X-\Si,\sfrak)$ is compact by Corollary~\ref{Empty}. 
For a compact perturbation $\eta_o$, however, the relative moduli spaces $\cM_{\eta_o}(X-\Si,\sfrak)$ have the same associated tunneling spaces as in the unperturbed case and these tunneling spaces can be non-trivial.
\medskip
\par\noindent
In the following, we will assume that $\Si\cdot\Si\!\neq\!0$ to ensure that the \spinc structure on $Y$ is torsion and thus $\tn{CSD}$ is real-valued and increasing along the flow lines on $\cM(Y,\mfsy),$ as discussed in Sec.~\ref{SW4cyend_ss}. The case of $\Si\cdot\Si\!=\!0$ is slightly different and simpler and is already discussed in \cite{mst}.
\medskip
\par\noindent
Define $\ov\cM_{\eta_o}(X-\Si,\sfrak)$ to be the space of all level-$k$ broken trajectories for $\mfs$ in the following sense. A \textit{level-$k$} element of $\ov\cM_{\eta_o}(X-\Si,\sfrak)$ is a tuple 
\bEq{BT_e}
\Big([A_0,\Phi_0];
\llbracket A_1,\Phi_1\rrbracket;\ldots;
\llbracket A_k,\Phi_k \rrbracket \Big)
\eEq
in the fiber product
\bEq{LevelK_e}
\cM_{\eta_o}(X-\Si,\sfrak_{\sss X-\Si})\, {}_{\tn{ev}}\!\!\times_{\tn{ev}} \Big(\cM(\R\times Y,\sfrak_{\sss \R\times Y})/\C^*\Big) \ldots\, {}_{\tn{ev}}\!\!\times_{\tn{ev}} \Big(\cM(\R\times Y,\sfrak_{\sss \R\times Y})/\C^*\Big),
\eEq
such that 
\bEn
\item the $\C^*$-action on $\cM(\R\times Y,\sfrak_{\sss \R\times Y})$ is given by the translation symmetry on $\R$ and $S^1$-rotation symmetry on $Y$;
\item the evaluation map $\tn{ev}$ at the ends takes value in a divisor space $\tn{Div}_{m}(\Si)$ as in (\ref{Ev-Map_e}), if the limit is irreducible, or in $\cJ\cong \T^{2g}$ (or $\T^{2g+1}$), if the limit is reducible;
\item the limit at $+\infty$ of the last monopole belongs to $\tn{Div}_{d(\mfs)}(\Si)$.
\eEn

\par\noindent
It is known that $\ov\cM_{\eta_o}(X-\Si,\sfrak)$ has a natural sequential convergence topology that is Hausdorff and compact (e.g. see \cite[Sec.~16]{KM}). Compared to the compactification in \cite{KM}, the only difference here is that we take the quotient by the action of $\C^*$ instead of $\R$, therefore our compactification has no boundary.

\medskip
\par\noindent
For a broken trajectory as in (\ref{BT_e}), let $C_i$ denote the component of the limit at $+\infty$ of $[A_i,\Phi_i]$ for $i\!\in\!\{0,\ldots,k\}.$ Therefore, the limit at $-\infty$ of $[A_i,\Phi_i]$ belongs to $C_{i-1}$ for $i\!\in\!\{1,\ldots,k\}$ and $C_{k}=\tn{Div}_{d(\mfs)}(\Si)$. The $C_i$'s are in fact components of the decomposition
\bEq{DecompM_e}
\cM(Y,\sfrak_{\sss Y})=\mc{J}\cup \bigcup_{m} \cM(Y,\sfrak_{\sss Y})_{m}\, .
\eEq
Using the notation of Sec.~\ref{SW4cyend_ss}, for $i\!\in\!\{1,\ldots,k\}$, each $\llbracket A_i,\Phi_i\rrbracket$ is an element of 
$$\cM(C_{i-1},C_i)/\C^*=\breve\cM(C_{i-1},C_i)/S^1.$$
Recall from Sec.~\ref{SW4cyend_ss} that the components $C_0,\dots,C_k$ are ordered by the value of $\tn{CSD}$ in the sense that $\tn{CSD}(C_0)<\tn{CSD}(C_1)<\dots<\tn{CSD}(C_k).$

\medskip
\par\noindent
If none of these limits belongs to $\mc{J}$, it follows from Lemma~\ref{dim-Calc_lm1} that the expected dimension of such level-$k$ configurations is $2k$ lower than the expected dimension of  $\cM_{\eta_o}(X-\Si,\sfrak)$. 

\medskip
\par\noindent
On the other hand, depending on $\tn{sign}(\Si\cdot\Si)$, the reducible component $\mc{J}$ appears in two ways: 
\bIt
\item If $\Si\cdot\Si\!<\!0$, then $\tn{CSD}$ takes its minimum on $\mc{J}$;
\item If $\Si\cdot\Si\!>\!0$, then $\tn{CSD}$ takes its maximum on $\mc{J}$.
\eIt
To see this, first note that we can assume $\tn{CSD}(\mc{J})=0$ by taking our background connection to be a flat connection on $L_{\sss Y}$.
To evaluate $\tn{CSD}$ at an irreducible solution $(B,\Psi)$, recall that it is the pullback of a monopole on $\Si,$ so the holonomy of $B$ is $\exp(2\uppi\mfi c/\ell),$ where $c=\deg(L_\Si)$ and $\ell=\deg(Y)=-\Si.\Si.$ Therefore, for some flat connection $B_0$ on $L_{\sss Y},$ we have $B-B_0=(c/\ell)\mfi\al,$ where $\mfi\al$ is the connection 1-form associated to decomposition (\ref{Splitting_e}).
Since $\frac{\mfi}{2\uppi}\nd(\mfi\al)$ and $\ell.\piysi^*\om$ both represent the pullback of the (integral) Euler class of $Y,$ a straight-forward calculation shows that $\tn{CSD}(B,\Psi)$ is a positive multiple of $c^2/\ell$ and the two items above follow.
As a result, 
\bIt
\item If $\Si\cdot\Si\!<\!0$, then only $C_0$ can be equal to $\mc{J}$;
\item If $\Si\cdot\Si\!>\!0$, then none of $C_0,\ldots,C_k$ is equal to $\mc{J}$. 
\eIt
In the first case, where $C_0\!=\!\mc{J},$ the expected dimension of level-$k$ configurations is $2k+1$ lower than the expected dimension of $\cM_{\eta_o}(X-\Si,\sfrak)$ by lemmas~\ref{dim-Calc_lm1} and~\ref{dim-Calc_lmRed}. To summarize, we have the following theorem.

\bTh{Compactness}
The moduli space $\cM_{\eta_o}(X-\Si,\sfrak)$ naturally has a Hausdorff compactification $\ov\cM_{\eta_o}(X-\Si,\sfrak)$ with a ``boundary" of expected real codimension at least two. 
The moduli space $\cM_{\eta_\nu}(X-\Si,\sfrak)$ is compact.
\eTh

\begin{remark}
It may happen that (\ref{DecompM_e}) has only one irreducible component for topological reasons. If further $\Si\cdot\Si>0$, this implies that $\cM_{\eta_o}(X-\Si,\sfrak)$ is compact. For example, if $\Si\cdot\Si>2g-2$, since different $m$'s in (\ref{DecompM_e}) differ by a multiple of $\Si\cdot\Si$ and $\cM(Y,\mfsy)_{m}$ is empty unless $|m|\leq g-1,$ we conclude that $\cM_{\eta_o}(X-\Si,\sfrak)$ is compact. This point has been used in \cite{OS} to show that certain decompositions of 4-manifolds cannot happen.
\end{remark}

\begin{remark}
Subject to the regularity of the tunneling spaces, we get a gluing map 
$$
\aligned
\cM_{\eta_o}(X-\Si,\sfrak_{\sss X-\Si})\, {}_{\tn{ev}}\!\!\times_{\tn{ev}} \Big(\cM(\R\times Y,\sfrak_{\sss \R\times Y})/\C^*\Big) \ldots\, {}_{\tn{ev}}\!\!\times_{\tn{ev}}\Big(\cM(\R\times Y,\sfrak_{\sss \R\times Y})/\C^*\Big) \times_{(S^1)^k} \De^k& \\
\lra \cM_{\eta_o}(X-\Si&,\sfrak_{\sss X-\Si})
\endaligned
$$
as in (\ref{GCS_e}), where $\De\subset\C$ is a sufficiently small disk and the fiber product with the $i$-th copy of $\De$ is with respect to the $S^1$-actions on $\De$ and the $i$-th tunneling space. This gluing map gives the compactified moduli space $\ov\cM_{\eta_o}(X-\Si,\sfrak)$ the structure of a closed $C^0$-manifold in which level-$k$ configurations in (\ref{LevelK_e}) are embedded as a codimension $2k$ (or $2k+1$) submanifold. \qed
\end{remark}

\medskip
\par\noindent
We finish this section by some comments on the contribution of the component $\cM_{\eta_o}(X-\Si,\mfs_{\sss X-\Si},\cJ),$ consisting of monopoles ending at reducibles, to the relative invariants (\ref{SWI-Relative_e}) and the resulting sum formula (\ref{SW-Sum_e}).

\medskip
\par\noindent
By the discussion above, if $\Si\cdot\Si\!>\!0$, all the monopoles appearing in $\ov\cM_{\eta_o}(X-\Si,\sfrak)$ are automatically irreducible. Therefore, whenever the invariants are defined, the condition $b^{+}_{\sss X-\Si}\!>\!0$ is sufficient for the relative invariants $\SW^{\sss X,\Si}_{\eta_o}(\mfs;-)$ to be independent of the choice of generic $\eta_o$ and the cylindrical metric. However, the component $\cM_{\eta_o}(X-\Si,\mfs_{\sss X-\Si},\cJ)$ may appear in two ways in this paper.

\medskip
\par\noindent
Firstly, if $\Si\cdot\Si\!<\!0$, it contributes to the level-$k$ strata of the compactified moduli space $\ov\cM_{\eta_o}(X\!-\!\Si,\sfrak)$, where $[A_0,\Phi_0]\in\cM_{\eta_o}(X-\Si,\mfs_{\sss X-\Si},\cJ)$, $\llbracket A_1,\Phi_1 \rrbracket\in \cM(J,C_1)/\C^*$, and the rest of the tunnelings in (\ref{BT_e}) begin and end at irreducibles. The expected codimension of such a configuration is $2k+1\geq 3$. Therefore, in this case, reducibles happen in high codimension and do not contribute to (\ref{SWI-Relative_e}).

\medskip
\par\noindent
Secondly, the component $\cM_{\eta_o}(X-\Si,\mfs_{\sss X-\Si},\cJ)$ would appear in the proof of the alternative sum formula (\ref{SW-Sum_e}) in the following way. Proving (\ref{SW-Sum_e}) revolves around the following observation. 
Let $X$ be the connected sum of $X_1$ and $X_2$ along $\Si$ and $\mfs$ be a \spinc structure on $X.$ Consider the fiber product 
\bEq{Moduli-Sum_e}
\cM_{\eta_1}(X_1-\Si,\mfs|_{\sss X_1-\Si}) \times_{\cM(Y,\mfs_{\sss Y})} \cM_{\eta_2}(X_2-\Si,\mfs|_{\sss X_2-\Si})
\eEq
over $\cM(Y,\mfs_{\sss Y})$ via the limiting maps (\ref{Limit-Map_e}), where $\eta_1$ and $\eta_2$ are generic compact perturbations on $X_1-\Si$ and $X_2-\Si$, respectively, and $\eta=\eta_1+\eta_2$ is the resulting perturbation on $X.$
\footnote{Here, $X$ can be thought of as any of the glued manifolds $X_{\sss T}$ in Theorem~\ref{Gluing_thm} for some sufficiently large $T,$ with $\eta_T=\eta.$}
\medskip
\par\noindent
The fiber product (\ref{Moduli-Sum_e})
decomposes into a union of main components, indexed by various pairs of \spinc structures $(\mfs_1,\mfs_2)$ on $TX_1(-\log\Si)$ and $TX_2(-\log\Si)$, respectively, which restrict to $(\mfs|_{\sss X_1-\Si},\mfs|_{\sss X_2-\Si})$ and have the same degree on $\Si$, 
\bEq{Irr-Comp-of-Sum_e}
\bigcup_{\mfs_1\#\mfs_2=\mfs} \cM_{\eta_1}(X_1-\Si,\mfs_1) \times_{\tn{Div}_{d(\mfs_1)=d(\mfs_2)}(\Si)} \cM_{\eta_2}(X_2-\Si,\mfs_2)\, ,
\eEq
and the fiber product
\bEq{Red-Comp-of-Sum_e}
\cM_{\eta_1}(X_1-\Si,\mfs|_{\sss X_1-\Si},\mc{J}) \times_{\mc{J}} \cM_{\eta_2}(X_2-\Si,\mfs|_{\sss X_2-\Si},\mc{J}).
\eEq

\medskip
\par\noindent
The union in (\ref{Irr-Comp-of-Sum_e}) already embeds into $\cM_{\eta}(X,\mfs)$ by the Gluing Theorem~\ref{Gluing_thm}. Assuming that (\ref{Red-Comp-of-Sum_e}) also embeds into $\cM_{\eta}(X,\mfs)$, as in \cite{OS}, we can conclude that $\cM_{\eta}(X,\mfs)$ is the same as (\ref{Moduli-Sum_e}) except for a subset of real codimension 2 (consisting of broken trajectories involving non-trivial tunnelings).
This follows from a convergence result similar to the one used in the construction of the compactification $\ov\cM_{\eta_o}(X-\Si,\mfs)$ above; see further below. 

\medskip
\par\noindent
By Lemma~\ref{dim-Calc_lm1}, the expected dimension of each component in (\ref{Irr-Comp-of-Sum_e}) matches that of $\cM_{\eta}(X,\mfs).$ However, by Lemma~\ref{dim-Calc_lmRed}, the expected dimension of (\ref{Red-Comp-of-Sum_e}) is 1 less than the expected dimension of $\cM_{\eta}(X,\mfs).$ This explains why (\ref{Red-Comp-of-Sum_e}) does not contribute to the sum formula (\ref{SW-Sum_e}). This observation is used in \cite{OS} to prove their main theorem. 

\medskip
\par\noindent
With notation as above, in general, a sequence of monopoles $[A_{\sss T_i},\Phi_{\sss T_i}]$ in $\cM(X_{\sss T_i},\sfrak)$, where $\lim_{i\to \infty} T_i=\infty$, will, after passing to a sub-sequence, ``converge" to a finite chain of monopoles where the first one is a monopole on $X_1-\Si$, the last one is a monopole on $X_2-\Si$, and the rest are non-trivial monopoles on the cylinder $\R\!\times\! Y$. The proof is more or less identical to the proof of \cite[Thm. 16.1.3]{KM} in the following sense. 
Every compact sub-domain $K$ of $X_i-\Si$ can be identified with a sub-domain of $X_T$ for $T>T_{K}$. Restricted to $K$, after passing to a sub-sequence, $\{[A_{\sss T_i},\Phi_{\sss T_i}]\}_{T_i>T_K}$ converges to a monopole $[A_{1,K},\Phi_{1,K}]$. By considering a sequence of exhausting compact sets $K$, using unique continuity and a diagonal argument, we get a monopole $[A_1,\Phi_1]$ over $X_1-\Si$ whose restriction to any $K$ is $[A_{1,K},\Phi_{1,K}]$. Similarly, we get a monopole over $X_2-\Si$. However, as $T_i\lra \infty$, the energy of $[A_{\sss T_i},\Phi_{\sss T_i}]$ may concentrate at several locations along the expanding cylinder $[0,T_i]\times Y \subset X_{\sss T_i}$. The same proof as in \cite[Thm. 16.1.3]{KM} gives us the connecting tunneling monopoles between the monopoles on $X_1-\Si$ and $X_2-\Si$.

\medskip
\par\noindent
More precisely, define $\cM_{\eta_\infty}(X_{\infty},\sfrak)$ to be the space of all level-$k$ broken trajectories ($k\geq 0$) for $\mfs$ in the following sense. A level-$k$ element of $\cM_{\eta_\infty}(X_{\infty},\sfrak)$ is a tuple 
\bEq{BT_e2}
\Big([A_0,\Phi_0];
\llbracket A_1,\Phi_1\rrbracket;\ldots;
\llbracket A_k,\Phi_k \rrbracket; [A_{k+1},\Phi_{k+1}] \Big)
\eEq
in the fiber product
\bEq{LevelK_e2}
\cM_{\eta_1}(X_1-\Si,\sfrak_{\sss X_1-\Si})\, {}_{\tn{ev}}\!\!\times_{\tn{ev}} \Big(\frac{\cM(\R\times Y,\sfrak_{\sss \R\times Y})}{\C^*}\Big) \ldots\, 
\, {}_{\tn{ev}}\!\!\times_{\tn{ev}}
\cM_{\eta_2}(X_2-\Si,\sfrak_{\sss X_2-\Si}),
\eEq
such that 
\bEn
\item the $\C^*$-action on $\cM(\R\times Y,\sfrak_{\sss \R\times Y})$ is given by the translation symmetry on $\R$ and $S^1$-rotation symmetry on $Y$;
\item the evaluation map $\tn{ev}$ at the ends takes values in a divisor space $\tn{Div}_{m}(\Si)$ as in (\ref{Ev-Map_e}), if the limit is irreducible, or in $\cJ\cong \T^{2g}$ (or $\T^{2g+1}$), if the limit is reducible.
\eEn
An analogous construction exists for adapted perturbations $\eta_\nu,$ in which the intermediate tunneling spaces are empty and we simply have the fiber product of $\cM_{\eta_i}(X_i-\Si,\sfrak_{\sss X_i-\Si}).$
\medskip
\par\noindent
With notation as above, we have the following compactness theorem.

\bTh{convergence_thm}
The union 
$$
\bigcup_{T\in [0,\infty]} \cM_{\eta_T}(X_{T},\sfrak) 
$$
has a natural sequential convergence topology that is Hausdorff and compact.
\eTh

\noindent
Replacing $\C^*$ with $\R$, the same theorem holds for a decomposition of $X$ along an arbitrary $3$-dimensional manifold $Y$.

\subsection{Conclusion; proofs of Theorems~\ref{Relaive_thm}--\ref{SW-Sum_thm2}}\label{Proofs_ss}
Except for the orientability of the relative moduli spaces $\cM_*(X-\Si,\sfrak),$ so far we have discussed all the steps that go into proving Theorems~\ref{Relaive_thm}--\ref{SW-Sum_thm2}. In this section, we wrap up the proofs with some comments on the orientation problem. The SW moduli spaces of closed 4-manifolds are orientable; a choice of orientation on $H^1(X)$ and $H^2_+(X)$ determines an orientation on the moduli space. Unfortunately, a simple statement like this for the SW moduli spaces of cylindrical-end manifolds does not exist in the literature; we refer to \cite[Sec.~4.4.3]{N} and \cite[Sec.~20]{KM} for a rather lengthy discussion of the problem. 
The problem is on the cylindrical part, where it is not clear if all the operators $\wt\cH_t$ in (\ref{cHt_e}) are Fredholm.
In Section~\ref{LSW_s}, we introduce a different setup for constructing relative moduli spaces which, among other things, could provide a short answer to the orientation problem as well; see the discussion after Definition~\ref{LogSWEq_e}.

\newtheorem*{proofof-Relaive_thm}
{Proof of Theorem~\ref{Relaive_thm}}
\begin{proofof-Relaive_thm}
The fact that, for generic $\eta_o$, $\cM_{\eta_o}(X-\Si,\mfs)$ is a smooth manifold and the evaluation map (\ref{Ev-Map_e}) is a submersion is proved in \cite{N}; see Remark~\ref{Generic-eta_rmk}. The dimension formula (\ref{DimFormula_e}) is derived in Lemma~\ref{dim-Calc_lm1} by simplifying Nicolaescu's formula \cite[(3.27)]{N1}.
Since all of the non-standard analysis happens on the cylindrical part of $X-\Si$, the orientability of $\cM_{\eta_o}(X-\Si,\mfs)$ can be proved in the same way as in \cite[Cor. 20.4.1]{KM}. Here, for the sake of completeness, we present an ad~hoc way of proving orientability. We show that $\cM_{\eta_o}(X-\Si,\sfrak)$ embeds in the moduli space $\cM_{\eta_o}(X,\sfrak')$ of the closed-up manifold $X$ with a related classical \spinc structure $\sfrak'$ on $TX.$ Therefore, since $X$ is closed, $\cM_{\eta_o}(X,\sfrak')$ can be oriented by a choice of homology orientation on $H^1(X,\Z)\oplus H^2_+(X,\Z),$ which in turn induces an orientation on its submanifold $\cM_{\eta_o}(X-\Si,\sfrak)$ of the same dimension.
\medskip
\par\noindent
In the context of the gluing map (\ref{GCS_e}), consider the natural decomposition of $X=X_{\sss T}$ to the cylindrical-end manifolds $M_+\!=\!X-\Si$ and $M_-\!=\!\cN$, where $\cN$ is the normal bundle of $\Si\subset X,$ as defined in Sec.~\ref{LogStr_ss}. With notation as in Sec.~\ref{Tunneling_ss}, the cylindrical-end manifold $\cN$ can also be seen as $P-\Si_+$. Recall that $\cN$ has a canonical \spinc structure determined by its almost-complex structure, and any \spinc structure $\mfsn$ on $\cN$ can be obtained by twisting this canonical \spinc structure with $E,$ where $E$ is the pullback of a line bundle of degree $d$ on $\Si.$ We will be interested in the moduli space $\cM(\cN,\mfsn)$ of the cylindrical-end manifold $\cN$, where $\mfsn$ is determined by the line bundle $E$ with $d=d(\mfs).$ By gluing the pair $\mfs_+=\mfs$ on $M_+$ and $\mfs_-=\mfsn$ on $M_-=\cN$, we get a \spinc structure $\mfs'$ on $X$. If $\cM(\cN,\mfsn)$ is regular, the gluing map (\ref{GCS_e}) then gives an embedding 
$$
\cM_{\eta_{o}}(X-\Si,\mfs)\, \times_{\cM(Y,\mfsye)\cong\tn{Div}_d(\Si)} \cM(\cN,\mfsn)\lra \cM_{\eta_{\sss T}}(X,\mfs').
$$
If we show that $\cM(\cN,\mfsn)$ has a component isomorphic to $\tn{Div}_d(\Si)$, for that component, the left-hand side of the map above is simply $\cM_{\eta_{o}}(X-\Si,\mfs)$. Therefore, we get an embedding of equi-dimensional manifolds
$$
\cM_{\eta_{o}}(X-\Si,\mfs)\lra 
\cM_{\eta_{\sss T}}(X,\mfs').
$$
The following lemma finishes the proof.
\bLm{}
For the \spinc structure $\mfsn$ defined by $E$ as above, where $0<|d-g+1|\leq g-1,$ the moduli space $\cM(\cN,\mfs_{\sss\cN})$ is non-empty and has a regular component diffeomorphic to the irreducible component $\cM(Y,\mfsye)\cong\tn{Div}_d(\Si)$ of $\cM(Y,\mfs_{\sss Y}).$
\eLm

\bPf
As we have seen, if $0<|d-g+1|\leq g-1,$ then $\cM(Y,\mfsye)\cong\tn{Div}_d(\Si)$ is an irreducible component of $\cM(Y,\mfs_{\sss Y}),$ and if $d=g-1,$ we get the reducible component.
Let us start by constructing a family of metrics on $Y$ as follows. 
Recall that $TY$ decomposes as in (\ref{Splitting_e}), where the vertical sub-bundle is given as the kernel of a real 1-form $\al$ on $Y$ and the horizontal sub-bundle identifies with $\piysi^*T\Si.$
Let $\gmetric_{\sss\Si}$ denote the K\"ahler metric on $\Si$ corresponding to the K\"ahler form $\om.$ 
We can now use $\al$ and $\gmetric_{\sss\Si}$ to define a family of metrics $\gmetric_\la=(\la\al)\!\otimes\!(\la\al)+\piysi^*\gmetric_{\sss\Si}$ on $Y$ for any $\la>0.$
\\
\par\noindent
Observe meanwhile that, according to the descent process discussed in Section~\ref{SW3_ss}, if $(A,\Psi)$ is an irreducible monopole on $\Si,$ then its pullback $(B,\Psi)$ via $Y\!\lra\!\Si$ is also an irreducible monopole on $Y,$ regardless of which of the metrics $\gmetric_\la$ is used on $Y.$
\\
\par\noindent
Now, we can view $\cN$ topologically as a quotient of $[0,\infty)\!\times\! Y,$ where $\{0\}\!\times\! Y$ is identified with $\Si$ via the bundle projection $Y\!\overset{\piysi}{\lra}\!\Si.$
Using the family of metrics defined above, we can equip $\cN$ with a cylindrical-end riemannian metric as follows. Start by considering the following family of metrics $\gmetric_r$ on the slices $Y_r\!=\!\{r\}\!\times\! Y$ of $\cN,$ $r>0,$
$$\gmetric_r=\al_r\!\otimes\!\al_r+\piysi^*\gmetric_{\sss\Si}, \qquad \al_r=\be(r)\al,$$
where $\be\colon(0,\infty)\!\lra\![0,1]$ is a smooth increasing function such that $\be(r)\!=\! r$ for $r\!\leq\! 1/2$ and $\be(r)\!=\! 1$ for $r\!\geq\! 1.$
In other words, $\gmetric_r\!=\!\gmetric_\la$ for $\la\!=\!\be(r).$
As a result, this family of metrics has shrinking fibers on $Y$ as $r\!\to\! 0$ in one direction, and stabilizes at $\gmetric_r\!=\!\gmetric_1$ for $r\!\geq\! 1$ in the other direction.
Note by the way that the adiabatic connection on each slice $Y_r$ is compatible with the metric $\gmetric_r$ by construction.
Now define a ``radial" metric $\gmetric$ on $\cN\!-\!\Si$ by
\bEq{MetricN_e}
\gmetric=\nd r^2+\gmetric_r 
=\nd r \!\otimes\!\nd r+\al_r\!\otimes\!\al_r+\piysi^*\gmetric_{\sss\Si}.
\eEq
An elementary local calculation in polar coordinates \cite[Ch.~1, Sec.~3.4]{P} shows that $\gmetric$ has a limit as $r\!\to\! 0$ and smoothly extends to a (non-degenerate) riemannian metric on the entire $\cN.$ 
In a neighborhood of $r\approx0,$ we can formally write 
$\gmetric=\nd r^2+r^2\nd\theta^2+\gmetric_{\sss\Si},$
where the metric $\gmetric$ collapses to $\nd r^2+\gmetric_{\sss\Si}$ at $r=0$ as the slices $Y_r$ collapse to $\Si.$
The metric $\gmetric$ in a neighborhood of $\Si\subset\cN$ is the K\"ahler metric associated to the K\"ahler form (\ref{KahlerN_e}).
Moreover, $\gmetric$ is clearly cylindrical for $r\geq 1,$ where we can formally write $\gmetric=\nd r^2+\nd\theta^2+\gmetric_{\sss\Si}.$
\\
\par\noindent
We will now construct a monopole on $\cN$ from an arbitrary smooth irreducible monopole $(A,\Psi)$ on $\Si.$ To begin with, let us pull $(A,\Psi)$ back via $\piysi\colon Y\!\lra\!\Si$ to obtain a smooth irreducible monopole $(B,\Psi)$ on $Y,$ which is invariant under the circle action on $Y.$ We can extend this to a configuration on $[0,\infty)\!\times\! Y$ as a constant family $(B,\Psi)$ on each slice $Y_r,$ which in turn will descend to a smooth configuration $(\tilde{B},\tilde\Psi)$ on the quotient $\cN.$ We will show that $(\tilde{B},\tilde\Psi)$ satisfies the $\SW^4$ equations on $\cN.$ 
On $(0,\infty)\!\times\! Y,$ the spinor $\tilde\Psi$ is certainly harmonic, because it is harmonic on $Y$ and constant in the radial direction $(0,\infty).$ Moreover, the spinor bundles on the slices $Y_r$ are identified with the plus- and minus-spinor bundles on $(0,\infty)\!\times\! Y$ via Clifford multiplication by the unit cotangent vector $\nd r$ and, since we are using the adiabatic connection on each slice, the curvature $F_{\tilde B}$ has no component in the radial direction and the curvature equation in $\SW^4$ reduces to the curvature equation on $Y.$ We conclude that $(\tilde{B},\tilde\Psi)$ satisfies the $\SW^4$ equations on $\cN\!-\!\Si.$ Therefore, by smoothness, it will satisfy the equations over the entire $\cN.$
\\
\par\noindent
We have just showed that $\cM(\cN,\mfsn)$ is non-empty, so we can consider the limiting map
$$\partial\colon\cM(\cN,\mfsn)\lra\cM(Y,\mfsye)\cong\tn{Sym}^d(\Si)$$
for the cylindrical-end manifold $\cN.$
We have in fact constructed a right inverse to the limiting map in the previous paragraph, so $\tn{Sym}^d(\Si)$ embeds into $\cM(\cN,\mfsn).$
A calculation using Lemma~\ref{dim-Calc_lm1} shows that the expected dimension of $\cM(\cN,\mfsn)$ is equal to $2d,$ which matches the dimension of $\tn{Sym}^d(\Si).$ 
We conclude that at least a component of $\cM(\cN,\mfsn)$ is diffeomorphic to $\tn{Sym}^d(\Si)$ via the limiting map and the lemma follows.
\ePf

\end{proofof-Relaive_thm}

\newtheorem*{proofof-Relaive_thm2}
{Proof of Theorem~\ref{Relaive_thm2}}
\begin{proofof-Relaive_thm2}
As before, the fact that $\cM_{\eta_\nu}(X\!-\!\Si,\sfrak)$ is a smooth manifold for a generic choice of the compactly-supported part $\eta_o$ is proved in \cite{N}; see Remark~\ref{Generic-eta_rmk}.
The dimension formula (\ref{DimFormula_e2}) is proved in Lemma~\ref{dim-Calc_lm2}.
Compactness, as stated in Theorem~\ref{Compactness}, is derived from Corollary~\ref{Empty}. Orientability can be proved as above. 
The fact that, for different such $\eta_{\nu}$ and $\eta_{\nu'}$, $\cM_{\eta_\nu}(X\!-\!\Si,\sfrak)$ and $\cM_{\eta_{\nu'}}(X\!-\!\Si,\sfrak)$ are smoothly cobordant follows from considering a 1-parameter family of perturbations connecting $\eta_{\nu}$ and $\eta_{\nu'}$, as in the classic case.
\qed
\end{proofof-Relaive_thm2}

\newtheorem*{proofof-SW-Sum_thm2}
{Proof of Theorem~\ref{SW-Sum_thm2}}
\begin{proofof-SW-Sum_thm2}
The sum formula is a direct consequence of the gluing theorem and convergence/compactness argument at the end of Section~\ref{Compact_ss}. In other words, for sufficiently large $T$, the gluing map 
$$
\bigcup_{[\mfs]=\mfs_1\#\mfs_2} \cM_{\eta_{1,\nu}}(X_1-\Si,\mfs_1)\, \times_{S_{d(\mfs_1)}=S_{d(\mfs_2)}} \cM_{\eta_{2,\nu}}(X_2-\Si,\mfs_2)\lra \cM_{\eta_{{\sss T},\nu}}(X_{\sss T},\mfs)
$$
which is an embedding by Theorem~\ref{Gluing_thm}, is also onto by Theorem~\ref{convergence_thm}, thus 
is an identification of closed manifolds. Here $\eta_{1,\nu}$ and $\eta_{2,\nu}$ are adapted perturbations on $X_1-\Si$ and $X_2-\Si$ corresponding to a holomorphic $1$-form $\nu$ on $\Si$ and $\eta_{{\sss T},\nu}$ is the resulting perturbation term on $X_{\sss T}$.
\medskip
\par\noindent
To see the surjectivity of the gluing map, note that $X_{\sss T}$'s can be identified with each other in a natural way by re-scaling the neck in the time direction. Therefore, a given monopole in $\cM_{\eta_\nu}(X,\mfs)$ can be identified with a monopole in $\cM_{\eta_{{\sss T},\nu}}(X_{\sss T},\mfs)$ for any $T;$ these monopoles are compatible with each other and have the same energy on the neck. As $T\lra\infty,$ these monopoles converge to an element of the fiber product on the left, which is the desired inverse image of the original monopole on $X\simeq X_T$ under the gluing map above. 
\medskip
\par\noindent
As mentioned before, the $\pm$ signs $\ve_{\mfs_1,\mfs_2,q}\!\in\!\{\pm1\}$ in (\ref{SW-Sum_e2}) depend on the choices of orientations. 
For a meticulous discussion of how to orient the fiber products to be consistent with the gluing, see \cite[20.5]{KM}.
\qed
\end{proofof-SW-Sum_thm2}

\section{Logarithmic SW equations}\label{LSW_s}

In this section, we explain our idea for a direct construction of relative SW moduli spaces to bypass the issues related to working with non-closed manifolds. We will define logarithmic SW equations by replacing $TX$ with $TX(-\log \Si)$ and considering ``logarithmic connections". This construction can easily be generalized to the normal crossings case, where it is hard to work with $X-\Si$ as a cylindrical-end manifold. Through this construction, it should be possible to address the orientation problem more systematically.\\

\noindent
Let $(X,\Si)$ be as in the previous sections and $\mfs=(\sb,\clm)$ be a $\spinc$ structure on $TX(-\log \Si)$. 
If $W\!\lra\! X$ is a vector bundle over $X$, a connection $\nabla$ on $W$ is a bi-linear map
$$
\nabla\colon \Gamma(X,TX)\times \Gamma(X,W)\lra \Gamma(X,W), \qquad (\xi,\ze)\lra \nabla_{\xi}\ze,
$$
that is tensorial in the first input and satisfies the Leibniz rule in the second input. In the classical theory, the construction of SW moduli space involves a riemannian metric on $TX$, a hermitian metric on $\sb$, and compatible connections on $\sb$ and $TX.$ The latter is usually fixed to be the Levi-Civita connection.

\bDf{LogConnection_dfn}
Let $(X,\Si)$ be a pair of a closed oriented 4-manifold $X$ and a closed oriented 2-dimensional submanifold $\Si\!\subset\!X$. 
For any vector bundle $W\!\lra\! X$, a {\it logarithmic} connection $\nabla$ on $W$ is a bi-linear map
\bEq{LogNabla_e}
\nabla\colon \Gamma(X,TX(-\log \Si))\times \Gamma(X,W)\lra \Gamma(X,W), \qquad (\xi,\ze)\lra \nabla_{\xi}\ze,
\eEq
that is tensorial in the first input and satisfies the Leibniz rule in the second input.  \eDf
\noindent
As in the classical case, in any local trivialization, we have $\nabla= \nd_{\iota(\cdot)} +\Theta $, where  $\iota$ is the homomorphism in (\ref{homo_e}) and $\Theta$ is a matrix of logarithmic $1$-forms. Therefore, globally, every two logarithmic connections $\nabla$ and $\nabla'$ on $W$ differ by an $\tn{End}(W)$-valued logarithmic $1$-form 
$$
\Theta\in\Gamma\big(X, T^*X(\log \Si)\otimes \tn{End}(W)\big).
$$

\bDf{LogConMetric_dfn}
In the presence of a metric $\ll \cdot,\cdot\rr$ on $W$, we say $\nabla$ in  (\ref{LogNabla_e}) is compatible with the metric if 
$$
\nd_{\iota(\xi)}\ll \ze_1,\ze_2\rr= \ll \nabla_\xi \ze_1,\ze_2\rr +\ll \ze_1,\nabla_\xi\ze_2\rr.
$$
\eDf

\par\noindent
With the definition above, for $(X,\Si)$ and a \spinc structure $\mfs=(\sb,\clm)$ on $TX(-\log\Si),$ we need logarithmic connections $\nabla$ and $\nabla^\sb$ on $TX(-\log\Si)$ and $\sb$ which are compatible with the metrics on $TX(-\log\Si)$ and $W$, respectively, as well as with the Clifford multiplication:
$$
\nabla^\sb\big(\xi\cdot\Phi\big)=\xi\cdot\nabla^\sb\Phi+(\nabla\xi)\cdot\Phi.
$$
As in the classical case, a compatible connection $\nabla^{\sb}$ on $\sb$ is uniquely determined by $\nabla$ and a logarithmic connection $A$ on the characteristic line bundle $L_\sfrak$.

\medskip
\par\noindent
Associated to a logarithmic connection $\nabla^{\sb}$ as above we define the {\it logarithmic Dirac operator} to be 
\bEq{DiracLog_e}
\dirac^{\log}\colon \Ga(\sb^\pm)\lra \Ga(\sb^\mp), \qquad 
\dirac^{\log}\Phi=\sum_{{i=1}}^{4} e_i\cdot\nabla_{e_i}^\sb\Phi,
\eEq
where $e_1,\dots,e_4$ is an orthonormal basis for $T_x X(-\log\Si)$.
%
The metric considered on $TX(-\log\Si)$ is the one described before (\ref{Restriction_e}): on a neighborhood $\cD$ in $\cN,$ identified with a neighborhood of $\Si$ in $X$ using the map $\Upsilon,$ 
the metric is the direct sum of the pullback of the K\"ahler metric on $T\Si$ and the standard riemannian metric on $\C_\cD$ via the identification 
$$
\Upsilon^* TX(-\log\Si)=\pi^*T\Si\oplus\C_\cD.
$$
We take $\na$ to be the direct sum connection on $\cD$ and the Levi-Civita connection outside a larger neighborhood $(1+\ve)\cD$ and splice them in the middle using a convex combination with suitable smooth coefficients.
For each fixed-radius circle bundle $Y\subset \cD$, the restriction of $\nabla$ to $TY$ is the adiabatic connection mentioned in Section~\ref{SW3_ss}. Therefore, restricted to $X-\Si$, via the identification (\ref{Restriction_e}), $\nabla|_{T(X-\Si)}$ is the connection considered in the definition of relative moduli spaces $\cM(X-\Si,\mfs).$
\medskip
\par\noindent
As expected, the logarithmic Dirac operator (\ref{DiracLog_e}) is not elliptic on the entire $X$ in the classical sense of the word.  
The principal symbol of (\ref{DiracLog_e}), which is a function on the cotangent bundle $T^*X,$ is zero on the dual space $T\Si^\perp \subset T^*X|_{\Si}$ of $T\Si,$ while it is non-zero everywhere else. 
This is essentially due to the fact that the homomorphism $\iota$ maps $\partial_z^{\log}$ to $z\partial_z$ in $TX,$ which is zero along $\Si.$
We expect though that a {\it logarithmic} elliptic theory could be developed for $\dirac^{\log}$ that paves the way for working directly over $TX(-\log\Si)$ instead. 
In such a theory, the principal symbol of $\dirac^{\log}$ would rather be a function on the logarithmic cotangent bundle $T^*X(\log\Si),$ which, by analogy, is Clifford multiplication by that cotangent vector.

\medskip
\par\noindent
Next, we define the curvature of a logarithmic connection 
(\ref{LogConnection_dfn}). 

\bLm{Curvature_lm}
For $\xi_1,\xi_2\in \Gamma(X,TX(-\log \Si))$, there exists a unique
$$
\xi\defeq[\xi_1,\xi_2]\in \Gamma(X,TX(-\log \Si))
$$
such that 
\bEq{LogBracket_e}
\iota(\xi)=[\iota(\xi_1),\iota(\xi_2)].
\eEq
\eLm

\bPf
Fix local holomorphic coordinates $(z,w)\colon V\!\lra\!\C^2$ around a point of $\Si$ with $\Si\cap V \!=\! (z\equiv 0)$. 
For $a=1,2$, if 
$$
\xi_a= f_{a,1}\, \partial^{\log}_z+ 
f_{a,2}\, \partial^{\log}_{\ov{z}}+
f_{a,3}\, \partial_w+ 
f_{a,4}\, \partial_{\ov{w}}
$$
then 
$$
\iota(\xi_a)= f_{a,1}\, z\partial_z+ 
f_{a,2}\, \ov{z}\partial_{\ov{z}}+
f_{a,3}\, \partial_w+ 
f_{a,4}\, \partial_{\ov{w}}.
$$
We have
$$
\aligned
&[f_{1,1}\, z\partial_z,f_{2,1}\,z\partial_z]=  
\Big(f_{1,1}\frac{\partial (z f_{2,1})}{\partial z}-
f_{2,1}\frac{\partial (z f_{1,1})}{\partial z}\Big)\, z\partial_z, \\
&[f_{1,2}\, \ov{z}\partial_{\ov{z}},f_{2,1}\,z\partial_z]=
\Big(f_{1,2}\frac{\partial f_{2,1}}{\partial \ov{z}}\, \ov{z}\Big)\, z\partial_z
-\Big(f_{2,1}\frac{\partial f_{1,2}}{\partial z}\, z\Big)\, \ov{z}\partial_{\ov{z}},
\\
&[f_{1,3}\, \partial_w,f_{2,1}\,z\partial_z]=
\Big(f_{1,3}\frac{\partial f_{2,1}}{\partial w}\Big)\, z\partial_z
-\Big(f_{2,1}\frac{\partial f_{1,3}}{\partial z}\,z \Big)\, \partial_{w},
\\
&[f_{1,4}\, \ov{\partial_w},f_{2,1}\,z\partial_z]= 
\Big(f_{1,4}\frac{\partial f_{2,1}}{\partial \ov{w}}\Big)\, z\partial_z
-\Big(f_{2,1}\frac{\partial f_{1,4}}{\partial z}\,z \Big)\, \partial_{\ov{w}}.
\\
\endaligned
$$
Similarly, we see that (\ref{LogBracket_e}) holds for the rest of the terms.
Uniqueness follows from continuity and the fact that (\ref{LogBracket_e}) is the same as ordinary bracket away from $\Si.$
\ePf
\noindent
It follows from Lemma~\ref{Curvature_lm} that given a logarithmic connection $\nabla$ as in (\ref{LogNabla_e}), the curvature equation
$$
F_\nabla(\xi_1,\xi_2)\ze=
\nabla_{\xi_1} \nabla_{\xi_2} \ze-
\nabla_{\xi_2} \nabla_{\xi_1} \ze-
\nabla_{[\xi_1,\xi_2]} \ze
$$
is well-defined. Furthermore, $F_\nabla$ is tensorial; i.e., 
$$
F_\nabla\in \Gamma\big(X, \Lambda^{2}T^*X(\log \Si) \otimes \tn{End}(W)\big).
$$
Note that $T^*X(\log \Si)$ and thus $\Lambda^{2}T^*X(\log \Si)$ inherit metrics from $TX(-\log \Si)$. Therefore, the star operator $*$ on $\Lambda^{2}T^*X(\log \Si)$ is defined.
As in the classical case, let $\Lambda^{2,+}T^*X(\log \Si)$ denote the subspace of self-dual elements. As in (\ref{CM4_e}), the homomorphism
\bEq{LogIso_e}
\Lambda^{2,+}T^*X(\log \Si)\otimes_\R\C \lra \tn{End}_0(\sb^+), \qquad 
\sum_{i<j} c_{ij}e_i^*\wedge e_j^* \lra \sum_{i<j} c_{ij}\clm(e_i)\clm(e_j)
\eEq
is an isomorphism.

\bDf{LogSWEq_e}
Let $(X,\Si)$ be a pair of a closed oriented 4-manifold $X$ and a closed oriented 2-dimensional submanifold $\Si\!\subset\!X$. 
Let $\mfs=(\sb,\clm)$ be a $\spinc$ structure on $TX(-\log \Si)$. We define the (unperturbed) {\it Logarithmic Seiberg-Witten} (or LSW) monopole equations to be the system of equations 
$$
F_A^+=(\Phi\Phi^*)_0, \qquad \dirac^{\log}_A\Phi=0, \eqno{(\tn{SW}^4_{\log})}
$$
for a pair $(A,\Phi)$ in the configuration space $\cC(X,\sfrak)\!=\!\mc{A}_{\log}(L)\!\times\!\Gamma(\sb^+),$ where $A$ is a logarithmic connection on the characteristic line bundle of $\sfrak$, $\Phi$ is a plus-spinor, and the curvature equation is with respect to the identification (\ref{LogIso_e}).
\eDf

\medskip
\par\noindent
The Logarithmic Seiberg-Witten moduli space is the quotient
$$
\cM(X,\Si,\sfrak)\defeq
\{(A,\Phi)\in\mc{A}_{\log}(L)\times\Ga(\sb^+)~\tn{satisfying}~\tn{SW}^4_{\log}\}/\mc{G}. 
$$
The restriction of $\tn{SW}^4_{\log}$ to $X\!-\!\Si$ is the classical cylindrical-end $\SW^4$ equations considered in the previous sections. We expect that 
\bEq{Claim_e}
\cM(X,\Si,\sfrak)=\cM(X-\Si,\sfrak);
\eEq
i.e., the answer to questions (1) and (2) in Remark~\ref{Questions} is positive. Furthermore, local calculations suggest that
\bIt
\item the operator 
$$
\dirac^{\log}\colon \Ga_{2}^2(\sb^+)\lra \Ga_1^2(\sb^-)
$$
is Fredholm (even though it is not elliptic), where the Banach space completions are with respect to a classical metric on $X$ and the given metric on $W$;
\item the sequence 
\bEq{log-d_e}
0\longrightarrow
\Om^0(X;\R)
\stackrel{\nd}{\longrightarrow}
\Om^1_{\log}(X;\R)
\stackrel{\nd^+}{\longrightarrow}
\Om^{2,+}_{\log}(X;\R)
\longrightarrow 0
\eEq
has finite cohomology.
\eIt
If this is the case, then the index of $\cE_{\sss X,\Si}(A,\Phi)$ is the sum of the indices of $\dirac^{\log}$ and (\ref{log-d_e}), and orienting the cohomology classes of (\ref{log-d_e}) orients the moduli space $\cM(X,\Si,\mfs)$. We hope to shed more light on these questions in the near future.

\medskip
\par
\noindent
The equality (\ref{Claim_e}) extends to the perturbed moduli spaces considered in this paper. If $\eta_o$ is a self-dual $2$-form supported in $X-\Si$, the definition of $\cM_{\eta_o}(X,\Si,\mfs)$ is a straightforward generalization of $\cM(X,\Si,\mfs)$. Likewise, an adapted perturbation $\eta_\nu$ on $X-\Si$ extends to a logarithmic self-dual $2$-form on $X$ with a non-trivial residue ($1$-form) on $\Si$. Then, the definition of $\cM_{\eta_\nu}(X,\Si,\mfs)$ is a straightforward generalization of $\cM(X,\Si,\mfs)$ obtained by replacing the first equation in $\SW^4_{\log}$ with
$$
F_A^+=(\Phi\Phi^*)_0-\eta_\nu \in \Om^{2,+}_{\log}(X;\R).
$$


\end{document}